\documentclass[final,leqno]{siamltex}

\usepackage{amsmath,amssymb}
\usepackage[pdftex]{graphicx,color}
\usepackage{epstopdf}
\usepackage{multicol}
\usepackage{enumerate}
\usepackage{float}
\usepackage{hyperref} 
\hypersetup{
 colorlinks=true,
 linktocpage=true,
 bookmarksnumbered=true
}
\usepackage{algorithm,algorithmic}
\usepackage[all]{xy}

\newtheorem{thm}{Theorem}[section]
\newtheorem{rem}[thm]{Remark}

\newcommand{\seqof}[3]{(#1)_{#2}^{#3}}

\newcommand{\ssa}{\text{SSA}}

\newcommand{\MNRM}{\text{MNRM}}
\newcommand{\TL}{\text{TL}}
\newcommand{\ML}{\text{ML}}

\newcommand{\NSSAKone}{N_{\ssa,K1}}
\newcommand{\NSSAKtwo}{N_{\ssa,K2}}

\newcommand{\NMNRMKone}{N_{K1}}
\newcommand{\NMNRMKtwo}{N_{K2}}
\newcommand{\NMNRMKoneC}{N_{K1}^{(c)}}
\newcommand{\NMNRMKtwoC}{N_{K2}^{(c)}}

\newcommand{\NSSAP}{N_{{\ssa}^*}}
\newcommand{\NTL}{N_{\TL}}

\newcommand{\NTLC}{N_{\TL}^{(c)}}

\newcommand{\EE}{\mathcal{E}}

\newcommand{\WE}{\EE_I}

\newcommand{\WEH}[1]{\hat \EE_{I,#1}}

\newcommand{\avg}[2]{\mathcal{A}\left(#1;#2\right)}
\newcommand{\svar}[2]{\mathcal{S}^2\left(#1;#2\right)}

\newcommand{\dbar}[1]{\Bar{\Bar{#1}}}

\newcommand{\ie}{\emph{i.e.}}
\newcommand{\eg}{\emph{e.g.}}

\newcommand{\ud}{\mathrm{d}}
\newcommand{\prob}[1]{\mathrm{P}\left(#1\right)}
\newcommand{\expt}[1]{\mathrm{E}\left[#1\right]}

\newcommand{\var}[1]{\mathrm{Var}\left[#1\right]}
\newcommand{\cov}[2]{\mathrm{Cov}\left[#1, #2\right]}
\newcommand{\norm}[1]{\left\|#1\right\|}
\newcommand{\abs}[1]{\left|#1\right|}
\newcommand{\rset}{\mathbb{R}}

\newcommand{\zset}{\mathbb{Z}}
\newcommand{\JAC}{\mathbb{J}}

\newcommand{\Ordo}[1]{{\mathcal{O}}\left(#1\right)}

\newcommand{\ordo}[1]{{o}\left(#1\right)}

\newcommand{\red}[1]{{\color{red}#1}}

\newcommand{\latt}{\mbox{$\zset_+^d$}}

\newcommand{\indicator}[1]{\mathbf{1}_{#1}} 

\newcommand{\PERIOD}{.}
\newcommand{\COMMA}{,}

\newcommand{\LP}{\left(}
\newcommand{\RP}{\right)}

\newcommand{\SEP}{\, \big| \,}
\newcommand{\BX}{\bar X(t)}
\newcommand{\BXi}{\bar X_i(t)}

\newcommand{\sj}{\sum_{j=1}^J}

\newcommand{\Qi}{Q_i(t,\tau_i)}

\newcommand{\ldi}{\log(\delta_i)}
\newcommand{\taui}{\tau_i(s)}
\newcommand{\DEN}[1]{- a_0 \LP \bar X(t) \RP +  \sj a_j(\bar X(t))e^{-#1\, \nu_{ji}}}
\newcommand{\aj}{ a_j(\bar X(t))}

\newcommand{\tsi}{\tilde{s}_i}

\newcommand{\az}{a_0(\bar X(t))}

\newcommand{\mV}[1]{\mathcal{V}_{#1}}
\newcommand{\hmV}[1]{\hat{\mathcal{V}}_{#1}}

\title{Multilevel Hybrid Chernoff Tau-leap}

\author{Alvaro Moraes\thanks{Computer, Electrical and Mathematical Sciences and Engineering Division,
 King Abdullah University of Science and Technology (KAUST),
 Thuwal, Saudi Arabia ({\tt alvaro.moraesgutierrez@kaust.edu.sa}).}
        \and Raul Tempone\thanks{Computer, Electrical and Mathematical Sciences and Engineering Division,
 King Abdullah University of Science and Technology (KAUST),
 Thuwal, Saudi Arabia ({\tt raul.tempone@kaust.edu.sa}).}
\and Pedro Vilanova\thanks{Computer, Electrical and Mathematical Sciences and Engineering Division,
 King Abdullah University of Science and Technology (KAUST),
 Thuwal, Saudi Arabia ({\tt pedro.guerra@kaust.edu.sa}).}}

\date{Received: date / Accepted: date}

\numberwithin{equation}{section}
\numberwithin{figure}{section}
\numberwithin{table}{section}

\begin{document}
\maketitle

\begin{abstract}
In this work, we extend the hybrid Chernoff tau-leap method to the multilevel Monte Carlo (MLMC) setting. 
Inspired by the work of Anderson and Higham on the tau-leap MLMC method with uniform time steps, we develop a novel algorithm that is able to couple two hybrid Chernoff tau-leap paths at different levels. 
{Using dual-weighted residual expansion techniques, we also develop a new way to estimate the variance of the difference of two consecutive levels and the bias.
This is crucial because the computational work required to stabilize the coefficient of variation of the sample estimators of both quantities is 
 often unaffordable for the deepest levels of the MLMC hierarchy.}
Our method bounds the global computational error to be below a prescribed tolerance, $TOL$, within a given confidence level.
This is achieved with nearly optimal computational work. Indeed, the computational complexity of our method is of order $\mathcal{O}\left(TOL^{-2}\right)$, the same as with an exact method, but with a smaller constant. 
Our numerical examples show substantial gains with respect to the previous single-level approach and the Stochastic Simulation Algorithm.


\end{abstract}



\setcounter{tocdepth}{1}

\section{Introduction}
This work, inspired by {the multilevel discretization schemes introduced in} \cite{Anderson2012}, extends the {hybrid Chernoff tau-leap} method \cite{ourSL} to the {multilevel Monte Carlo} setting \cite{GilesMLMC}.
Consider a non-homogeneous Poisson process, $X$, taking values in the lattice of non-negative integers, $\zset_+^d$.
We want to estimate the expected value of a given observable, $g:\rset^d\to\rset$ of $X$, at a final time, $T$, \ie,
$\expt{g(X(T))}$.
For example, in a chemical reaction in thermal equilibrium, the $i$-th component of $X$, $X_i(t)$, could describe the number of particles of species $i$ present at time $t$. 
In the systems modeled here, different species undergo reactions at random times by changing the number of particles in at least one of the species. The probability of a single reaction happening in a small time interval is modeled by a non-negative {propensity function} that depends on the current state of the system. We present a formal description of the problem in Section \ref{sec:pjp}. 

Pathwise realizations of such pure jump processes (see, \eg, \cite{kurtzmp}) can be simulated exactly using the {Stochastic Simulation Algorithm} (SSA),
introduced by Gillespie in \cite{gillespie_ssa}, or the {Modified Next Reaction Method} (MNRM) introduced by Anderson in \cite{anderson2007modified}. 
Although these algorithms generate exact realizations for the Markov
process, $X$, they are computationally feasible for only relatively low propensities. 
  
For that reason, {Gillespie in \cite{gillespie_tau_leap} and  Aparicio and Solari in \cite{Aparicio2001} independently proposed} the {tau-leap method} to
{approximate} the SSA by evolving the process with fixed time
steps and by keeping the propensity fixed within each time step. In fact, the tau-leap method can be
seen as a forward Euler method for a stochastic differential equation driven by Poisson random measures (see, \eg, \cite{li:07}).

A drawback of the tau-leap method is that the simulated process may take negative values, which is an undesirable consequence of the approximation and not a qualitative feature of the original process. For this purpose, we proposed in \cite{ourSL} a Chernoff-type bound that controls the probability of reaching negative values
by adjusting the time steps. 
Also, to avoid extremely small time steps, we proposed  
switching adaptively between the tau-leap and an exact method, creating a {hybrid} tau-leap/exact 
method that combines the strengths of both methods. 

More specifically, let $\bar x$ be the state of the approximate process at time $t$, and let $\delta \in (0,1)$ be given. 
The main idea is to compute a time step, $\tau {=} \tau(\delta, \bar x)$, such that the probability that the approximate process reaches an unphysical negative value in $[t, t{+}\tau)$ is less than $\delta$. 
This allows us to control the probability that a entire hybrid path exits the lattice, $\latt$. In turn, this quantity leads to the definition of the global exit error, which is a global error component along with the time discretization error and the statistical error (see Section \ref{calibration} for details).

The multilevel Monte Carlo idea goes back at least to \cite{Heinrich, Hein98}. In that setting, the main goal was to solve high-dimensional, parameter-dependent integral equations and to conduct corresponding complexity analyses. Later, in \cite{GilesMLMC}, Giles developed and analyzed multilevel techniques that were used to reduce the computational work when estimating an expected value using Monte Carlo path simulations of a certain quantity of interest of a stochastic differential equation. {Independently, 
in \cite{speight09}}, Speight introduced a multilevel approach to control variates. Control variates are a widespread variance reduction technique with the main goal of increasing the precision of an estimator or reducing the computational effort. The main idea is as follows: to reduce the variance of the standard Monte Carlo estimator of  $\expt{X}$, 
$$\hat \mu_1 := \frac{1}{M} \sum_{m=1}^M X(\omega_m) \COMMA$$
we consider another unbiased estimator  of $\expt{X}$, 
$$\hat \mu_2 := \frac{1}{M} \sum_{m=1}^M \LP X(\omega_m)-(Y(\omega_m)-\expt{Y})\RP\COMMA$$
 where $Y$ is a random variable correlated with $X$ with known mean, $\expt{Y}$. 
The variable $Y$ is called a control variate. 
Since $\var{\hat \mu_2} {=} \var{\hat \mu_1} {+}  \var{Y} {-} 2\cov{X}{Y}$, whenever $\cov{X}{Y}{>}\var{Y}/2$, we have that $\var{\hat \mu_2}{\leq}\var{\hat \mu_1}$.
If we assume that the computational work of generating the pair $(X(\omega),Y(\omega))$ is less than twice the computational work of generating $X(\omega)$,
it is straightforward to conclude  that $\hat{\mu}_2$ is preferred when  
$\rho^2_{X,Y}{>}1/2$, where $\rho_{X,Y}$ is the correlation coefficient of the pair $(X,Y)$.    
We observe that $\hat \mu_2$ can be written as 
$$\hat \mu_2 = \expt{Y} +  \frac{1}{M} \sum_{m=1}^M \LP X-Y \RP (\omega_m)\PERIOD$$
In the case where $\expt{Y}$ is unknown and sampling from $Y$ is computationally less expensive than sampling from $X$, it is natural to estimate $\expt{Y}$ using Monte Carlo sampling to yield a two-level Monte Carlo estimator of $\expt{X}$ based on the control variate, $Y$, \ie,
$${\tilde \mu}_2 := \frac{1}{M_0} \sum_{m_0=1}^{M_0} Y(\omega_{m_0}) +  
\frac{1}{M_1} \sum_{m_1=1}^{M_1} \LP X-Y \RP(\omega_{m_1})\PERIOD$$
See Section \ref{MLMCintro} for details about the definition of levels in our context.

In this work, we apply Giles's multilevel control variates idea to the hybrid Chernoff  tau-leap approach to reduce the computational cost, which is measured as the amount of time needed for computing an estimate of $\expt{g(X(T))}$, within $TOL$, with a given level of confidence. 
We show that our hybrid MLMC method has the same computational complexity of the pure SSA, \ie, order $\Ordo{TOL^{-2}}$. 
From this perspective, our method can be seen as a variance reduction for the SSA since our MLMC method does not change the complexity; it just reduces the corresponding multiplicative constant.
We note in passing that in \cite{Anderson2013}, the authors show that the computational complexity for the pure MLMC tau-leap case has order $\Ordo{TOL^{-2} (\log(TOL))^2}$. {We note also that here our goal is to provide an estimate of $\expt{g(X(T))}$ in the probability sense and not in the mean square sense as in \cite{Anderson2012}}. 

The global error arising from our hybrid tau-leap MLMC method can naturally be decomposed into three components: the global exit error, the time discretization error and the statistical error. This global error should be less than a prescribed tolerance, $TOL$, with probability larger than a certain confidence level. The global exit error is controlled by the one-step exit probability bound, $\delta$ \cite{ourSL}. 
The time discretization error, inherent to the tau-leap method, is controlled through the size of the mesh, $\Delta t$, \cite{kt}. 
At this point, it is crucial to stress that, by controlling the exit probability of the set of hybrid paths,
we are indirectly turning this event into a rare event. 
Thus, direct sampling of exit paths is not an affordable way to estimate the probability of such an event.

{Motivated by the Central Limit results of Collier et al. \cite{Abdo2014} for the Multilevel Monte Carlo estimator (see appendix A, Theorem 1), we approximate the
statistical error with a Gaussian random variable with zero mean.} 
{In our numerical experiments, we tested this hypothesis by employing Q-Q plots and the Shapiro-Wilk test \cite{ShapiroWilk1965}. There, we did not reject the Gaussianity of the statistical error at the 1\% significance level. The variance of the statistical error} is a linear combination of the variance at the coarsest level and variances of the difference of two consecutive levels, which we sometimes call strong errors.  
{In Section \ref{sec:vl}, motivated by the fact that  sample variance and bias estimators are inaccurate on the deepest levels, we develop a novel dual-weighted residual expansion that allows us to estimate those quantities, cf. \eqref{eq:EILw}  
and \eqref{eq:varhatestimated}.}
We also control the statistical error
through the number of coupled hybrid paths, $\seqof{M_{\ell}}{\ell=0}{L}$, simulated at each level.

{We note that our use of duals in this work is different from the use  in \cite{kt}.
%
%
That earlier work  proposed an adaptive, single-level, tau-leap algorithm for error control, choosing the time steps non-uniformly to control the global weak error based on dual-weighted error estimators.
In this work, we do not have an adaptive time step based on  
dual-weighted error estimators as in \cite{kt}. We use instead dual-weighted error estimators to reduce the statistical error in our error estimates.}

\subsection{A Class of Markovian Pure Jump Processes}
\label{sec:pjp}
To describe the class of Markovian pure jump process, $X:[0,T]\times \Omega \to\zset_+^d$, that we use in this work, 
we consider a system of $d$ species interacting
through $J$ different reaction channels. For the sake of brevity, we write $X(t,\omega) {\equiv} X(t)$.  Let $X_i(t)$ be the
number of particles of species $i$ in the system at time $t$. We 
want to study the evolution of the state vector,
\begin{equation*}
  X(t) = (X_1(t), \ldots, X_d(t)) \in
  \zset_+^d,  
\end{equation*}
modelled as a continuous-time, discrete-space Markov chain starting at
some state, $X(0) \in \zset_+^d$. Each reaction can be described by the
vector $\nu_j\in \zset^d$, such that, for a state vector $x\in
\zset_+^d$, a single firing of reaction $j$ leads to the change
\begin{equation*}
  x \to x+ \nu_j.
\end{equation*}
The probability that reaction $j$ will occur during the small interval
$(t,t{+}\ud t)$ is then assumed to be
\begin{equation}\label{eq:prob}
  \prob{X(t+\ud t)=x+\nu_j | X(t) = x} =  a_j(x) \ud t + \ordo{\ud t}\COMMA
\end{equation}
with a given non-negative {polynomial propensity function}, $a_j:\rset^d \to
\rset$. We set $a_j(x){=}0$ for those $x$ such that $x{+}\nu_j\notin\latt$.
A process, $X$, that satisfies \eqref{eq:prob}, is a continuous-time, discrete-space
Markov chain that admits the following random time change representation \cite{kurtzmp}:
\begin{equation}
  \label{eq:exact_process}
  X(t) = X(0) + \sum_{j=1}^J  \nu_j Y_j \LP \int_0^t a_j(X(s))\, \ud s \RP  \COMMA
\end{equation}
where $Y_j:\rset_+ {\times} \Omega \to\zset_+$ are independent unit-rate Poisson processes.
Hence, $X$ is a non-homogeneous Poisson process.

In \cite{kt}, the authors assume that there exists a vector, $w\in\rset_+^d$, such that $(w,\nu_j)\leq 0 $, for any reaction $\nu_j$.
Therefore, every reaction, $\nu_j$, must have at least one negative component. This means that the species can be either transformed into other species or be consumed during the reaction. As a consequence, the space of states is contained in a simplex with vertices in the coordinate axis. This assumption excludes, for instance,  birth processes.
In our numerical examples, we allow the set of possible states of the system to be infinite, but we explicitly avoid cases in which one or more species grows exponentially fast or blows up in the time interval $[0,T]$.

\begin{rem}
In this setting, the solution of the following system of ordinary differential equations,
\begin{align*}
\left\{ \!
\begin{array}{lll}
\dot x (t)  &=&  \nu a(x(t)), \,\, t \in \rset_+ \\
x(0) &=& x_0 \in \rset_+ \COMMA
\end{array}
\right.
\end{align*}
is called mean field solution, where $\nu$ is the matrix with columns $\nu_j$ and $a(x)$ is the column vector of propensities. In Section \ref{EEC}, we use the mean field path for scaling and preprocessing constants associated with the computational work of the SSA and Chernoff tau-leap steps.
\end{rem}

%

\subsection{Description of the Modified Next Reaction Method (MNRM)}\label{sec:MNR}
The MNRM, introduced in \cite{anderson2007modified}, based on the Next
Reaction Method (NRM) \cite{GibBruck}, is an exact simulation algorithm like Gillespie's SSA that explicitly uses  representation \eqref{eq:exact_process} for simulating exact paths and generates only one exponential random variable per iteration. The reaction times are modeled with firing times of Poisson processes, $Y_j$, with
internal times given by the integrated propensity functions.
The randomness is now
separated from the state of the system and is encapsulated in the $Y_j$'s.
For each reaction, $j$, the internal time 
is defined as $R_j(t) {=} \int_0^t a_j(X(s))ds$. There are $J{+}1$ time frames in the system, the absolute one, $t$, and one for each Poisson
process, $Y_j$.
Computing the next reaction and its time is equivalent to computing how
much time passes before one of the Poisson processes, $Y_j$, fires, and
to identifying which process fires at that particular time, by taking the minimum of such times. The
NRM and MNRM make use of internal times to reduce the number of
simulated random variables by half.  In the following, we describe the MRNM and then we present its implementation in Algorithm \ref{alg:mnr}.

Given $t$, we have the propensity $a_j{=}a_j(X(t))$ and the internal
time $R_j{=}R_j(t)$. Now, let $\Delta R_j$ be the remaining
time for the reaction, $j$, to fire, assuming that $a_j$ stays constant
over the interval $[t,t{+}\Delta R_j)$. Then, $t{+}\Delta R_j$ is the time
when the next reaction, $j$, occurs. The next internal time at which the
reaction, $j$, fires is then given by $R_j {+} a_j \Delta R_j$. When
simulating the next step, the first reaction that fires occurs after
$\Delta {=} \min_j \Delta R_j$. We then update the state of the system according to
that reaction, add $\Delta$ to the global time, $t$, and then update the internal times by
adding $a_j \Delta$ to each $R_j$.  We are left to
determine the value of $\Delta R_j$, \ie, the amount of time until the
Poisson process, $Y_j$, fires, taking into account that $a_j$ remains
constant until the first reaction occurs. Denote by $R_j$ the first
firing time of $Y_j$ that is strictly larger than $R_j$, \ie,
$P_j {:=} \min  \{s{>}R_j:Y_j(s){>}Y_j(R_j)\}$ and finally $\Delta R_j {=} \frac{1}{a_j}(P_j-R_j)$.

\begin{algorithm}
\caption{The Modified Next Reaction Method. Inputs: the initial state, $X(0)$, the next grid point, $T_0> t_0$, the propensity functions, $\seqof{a_j}{j=1}{J}$, the stoichiometric vectors, $\seqof{\nu_j}{j=1}{J}$.
Outputs: the history of system states, $\seqof{X(t_k)}{k=0}{K}$. Here, we denote $S\equiv \seqof{S_j}{j=1}{J}$, $P\equiv \seqof{P_j}{j=1}{J}$, and $R\equiv \seqof{R_j}{j=1}{J}$.
}
\label{alg:mnr}
\begin{algorithmic}[1]
\STATE $k \leftarrow 0$, $t_k \leftarrow 0$, $X(t_k) \leftarrow X
(0)$ and $R \leftarrow \bf{0}$
\STATE Generate $J$ independent, uniform$(0,1)$ random numbers, $r_j$
\STATE $P \leftarrow \seqof{\log(1/r_j)}{j=1}{J}$
\WHILE {$t_k<T_0$}
\STATE $S \leftarrow \seqof{a_j(X(t_k))}{j=1}{J}$
\STATE $\seqof{\Delta R_j}{j=1}{J} \leftarrow \seqof{(P_j - R_j)/S_j}{j=1}{J}$
\STATE ${\mu} \leftarrow \text{argmin}_j \{\Delta R_j\}$
\STATE $\Delta \leftarrow \min_j\{ \Delta R_j \}$ 
\STATE $t_{k+1} \leftarrow t_k+\Delta$
\STATE $X(t_{k+1}) \leftarrow X(t_k) + \nu_\mu$
\STATE $R \leftarrow R + S \Delta$
\STATE $r \leftarrow$ uniform$(0,1)$
\STATE $P_{\mu} \leftarrow P_{\mu} + \log(1/r)$
\STATE $k \leftarrow k{+}1$
\ENDWHILE
\RETURN $\seqof{X(t_l)}{l=0}{k-1}$
\end{algorithmic}
\end{algorithm}


Among the advantages already mentioned, we can easily modify Algorithm 1 to generate paths in the cases where the rate functions depend on time and also when there are reactions delayed in time.
Finally, it is possible to simulate correlated exact/tau-leap paths using this algorithm as well as nested 
tau-leap/tau-leap paths. In \cite{Anderson2012}, this technique is used to develop a uniform-step, unbiased, multilevel Monte Carlo (MLMC) algorithm. In Section \ref{sec:chp}, we use this feature for coupling two exact paths.

\subsection{The Tau-Leap Approximation}\label{sec:num_approx}

In this section, we define $\bar X$, the tau-leap approximation of the process, $X$, which follows from applying the forward Euler approximation to the integral term in the following random time-change representation of $X$: 
\begin{equation*}
  \label{eq:exact_process_small}
  X(t+\tau) = X(t) + \sum_{j=1}^J  \nu_j Y_j \LP \int_t^{t+\tau} a_j(X(s))\, \ud s \RP  \PERIOD
\end{equation*}

The tau-leap method was proposed
in \cite{gillespie_tau_leap} to avoid the computational drawback of the exact methods, \ie,  when many reactions occur
during a short time interval. The tau-leap process, $\bar X$, starts from $X(0)$ at time $0$, and 
given that $\bar X(t) {=} \bar x$ and a
time step $\tau{>}0$, we have that
$\bar X$ at time $t{+}\tau$ is generated by
\begin{equation*}
  \bar X(t+\tau) = \bar x + \sum_{j=1}^J 
  \nu_j \mathcal{P}_j\left( a_j(\bar x)\tau \right) \COMMA
\end{equation*}
where $\{\mathcal{P}_j(\lambda_j)\}_{j=1}^J$ are independent
Poisson distributed random variables with parameter $\lambda_j$,
used to model the number of times that the reaction $j$ fires
during the $(t,t{+}\tau)$ interval.  Again, this is nothing other than a forward
Euler discretization of the stochastic differential equation formulation of the pure jump process
\eqref{eq:exact_process}, realized by the Poisson random measure with
state dependent intensity (see, \eg, \cite{li:07}).

In the limit, when $\tau$ tends to zero, the tau-leap method gives the same solution
as the exact methods. The total number of firings in each channel is 
a Poisson-distributed stochastic variable depending only on the
initial population, $\bar X(t)$.  The error thus comes from the
variation of $a(X(s))$ for $s\in(t,t{+}\tau)$. 

We observe that the computational work of a tau-leap step involves the generation of $J$ independent Poisson random variables. This is in contrast to the computational work of an exact step, which  involves only the work of generating two uniform random variables, in the case of the SSA, and only one in the case of MNRM. 

\subsection{The Chernoff-Based Pre-Leap Check}
\label{chernoff_onestep}
In \cite{ourSL}, we derived a Chernoff-type bound that allows us to guarantee that
the one-step exit probability in the tau-leap method is less than a
predefined quantity, $\delta{>}0$. 
We now briefly summarize the main idea.
Consider the following pre-leap check problem: find the largest possible
$\tau$ such that, with high probability, in the next step, the approximate process, $\bar{X}$, will take a value in the lattice, $\latt$, of non-negative integers.
%
The solution to that problem can be achieved by solving $d$ auxiliary
problems, one for each $x$-coordinate,
$i=1,2,\ldots,d$, as follows. Find the largest possible $\tau_i\geq 0$, such that
\begin{equation}
\label{eq:taui-deltai}
\prob{  \BXi + \sj \nu_{ji} \mathcal{P}_j\LP a_j\LP \bar X(t)\RP\tau_i\RP < 0 \Biggm| \BX}\leq \delta_i \COMMA
\end{equation}
where $\delta_i {=} \delta/d$, and $\nu_{ji}$ is the $i$-th coordinate of
the $j$-th reaction channel, $\nu_j$. Finally, we let $\tau{:=}\min\{\tau_i: i=1,2,\ldots,d\}$.
To find the largest time
steps, $\tau_i$, let $\Qi {:=} \sj (-\nu_{ji})\mathcal{P}_j\LP a_j\LP \bar X(t)\RP\tau_i\RP$. Then, for all $s{>}0$, 
%
we have the Chernoff
bound: 
\begin{equation*}
\prob{ \Qi{>}\BXi \SEP  \! \BX} \leq \inf_{s>0}  \exp \! \LP \! \! -s \BXi + \tau_i \sj \aj (e^{-s \nu_{ji}}{-}1) \! \RP \! .
\label{eq:chboundmulti}
\end{equation*}
Expressing $\tau_i$ as a function of $s$, we write
\begin{equation*}
\taui = \dfrac{\ldi + s \BXi}{\DEN{s}}=:\dfrac{R_i(s)}{D_i(s)}\COMMA
\end{equation*}
where 
\begin{equation*}
  \az := \sj \aj \PERIOD
\end{equation*}
We want to maximize $\tau_i$ while satisfying condition \eqref{eq:taui-deltai}. 
Let $\tau_i^*$ be this maximum.
We then have the following possibilities: If $\nu_{ji}\geq
0$, for all $j$, then naturally $\tau_i^*{=}+\infty$; otherwise, we have the
following three cases:
\begin{enumerate}
\label{summarizeChernoff}
\item \label{eq:tsi} $D_i(s_{i}){>}0$.  In this case,
  $\tau_i(s_{i}){=}0$ and $D_i(s)$ is positive and increasing as 
  $\forall s \geq s_{i}$. Therefore, $\taui$ is equal to the ratio of two
  positive increasing functions. The numerator, $R_i(s)$, is a linear
  function and the denominator, $D_i(s)$, grows exponentially
  fast. Then, there exist an upper bound, $\tau_i^*$, and a unique
  number, $\tsi$, which satisfy $\tau_i(\tsi){=}\tau_i^*$. We developed an algorithm in  \cite{ourSL} for approximating $\tsi$, using the relation
  $\tau'_i(\tsi){=}0$.
\item If $D_i(s_{i}){<} 0$, then  $\tau_i^*{=}+\infty$.
\item If $D_i(s_{i}){=}0$, then $\tau_i^*{=}\BXi/D'_i(s_{i})$.
\end{enumerate}
Here $s_{i} {:=}-\ldi/\BXi$.

\subsection{The Hybrid Algorithm for Single-Path Generation}
\label{hybrid_algo}

In this section, we briefly summarize our previous work, presented in \cite{ourSL}, on hybrid paths. 

The main idea behind the hybrid algorithm is the following. 
A path generated by an exact method (like SSA or MNRM) never exits the lattice, $\latt$, although the computational cost may be unaffordable due to many small inter-arrival times typically occurring when the process is ``far'' from the boundary. A tau-leap path, which may be cheaper than an exact one, could leave the lattice at any step. The probability of this event depends on the size of the next time step and the current state of the approximate process, $\BX$. This one-step exit probability could be large, especially when the approximate process is ``close'' to the boundary. We developed in \cite{ourSL} a Chernoff-type  bound to control the mentioned one-step exit probability. Even more, by construction, the probability that one hybrid path exits the lattice, $\latt$, can be estimated by
\begin{align*}
\prob{ A^c} \leq \expt{1-(1-\delta)^{N_{\TL}}} = \delta \expt{N_{\TL}} - \frac{\delta^2}{2}(\expt{N_{\TL}^2} - \expt{N_{\TL}}) + o(\delta^2),
\end{align*}
where  $\bar{\omega} \in A$ if and only if the whole hybrid path, $\seqof{\bar X(t_k,\bar \omega)}{k=0}{K(\bar{\omega})}$, belongs to the lattice, $\latt$, $\delta {>}0$ is the one-step exit probability bound, and $N_{\TL}(\bar \omega) {\equiv} N_{\TL}$ is the number of tau-leap steps in a hybrid path. Here, $A^c$ is the complement of the set $A$.

To simulate a hybrid exact/Chernoff tau-leap path, we first developed a one-step switching rule 
that, given the current state of the approximate process, $\BX$, adaptively determines whether to use an exact or an approximated method for the next step. This decision is based on the relative computational cost of taking an exact step (MNRM) versus the cost of taking a Chernoff tau-leap step.
We show the switching rule in Algorithm \ref{alg:sel}.
\begin{algorithm}
\caption{The one-step switching rule. Inputs: the current state  of the approximate process, $\BX$, the current time, $t$, the values of the propensity functions evaluated at $\BX$, $\seqof{a_j(\BX)}{j=1}{J}$, the one-step exit probability bound $\delta$, and the next grid point, $T_0$. Outputs: method and $\tau$. Notes:
based on $\expt{\tau_{\ssa}(\BX)\SEP \BX}=1/a_0(\BX)$ and $\tau_{Ch}(\BX,\delta)$, this algorithm adaptively selects between  MNRM and Chernoff tau-leap (TL).  
We denote by $\tau_{\MNRM}$  ($\tau_{Ch}$) the step size when the decision is to use the $\MNRM$  (tau-leap) method. 
}
\label{alg:sel}
\begin{algorithmic}[1]
	\REQUIRE $a_0 \leftarrow \sum_{j=1}^J a_j > 0$
	\IF{$K_1/a_0 < T_0-t$} 
		\STATE  $\tau_{Ch}\leftarrow $ compute Chernoff step size (see Section 2.2 in \cite{ourSL} )
		\IF{$\tau_{Ch} < K_2(\bar{X}(t),\delta)/a_0 $}
			\RETURN $(\MNRM, \tau_{\MNRM})$
		\ELSE
			\RETURN $(\TL, \tau_{Ch})$
		\ENDIF
	\ELSE
		\RETURN $(\MNRM, \tau_{\MNRM})$
	\ENDIF
\end{algorithmic}
\end{algorithm}
To compare the mentioned computational costs, we define $K_1$ as the ratio between the cost of computing $\tau_{Ch}$ and the cost of computing one step using the MNRM method, and $K_2{=}K_2(\BX,\delta)$ is defined as the cost of taking a Chernoff tau-leap step, divided by the cost of taking a MNRM step plus the cost of computing $\tau_{Ch}$. 
For further details on the switching rule, we refer to \cite{ourSL}.

\subsection{The Multilevel Monte Carlo Setting}
\label{MLMCintro}
In this subsection, we briefly summarize the control variates idea developed by Giles in \cite{GilesMLMC}. 
Let $\{\bar X_{\ell}(t)\}_{t\in[0,T]}$ be a hybrid Chernoff tau-leap process with a time mesh of size $\Delta t_{\ell}$ and a one-step exit probability bound, $\delta$. We can simulate paths of  $\{\bar X_{\ell}(t)\}_{t\in[0,T]}$ by using  Algorithm \red{4}  in \cite{ourSL}.
Let $g_{\ell}{:=} g({\bar X}_{\ell}(T))$.

Consider a hierarchy of nested meshes of the time interval $[0,T]$, indexed by $\ell = 0,1,\ldots,L$.
Let $\Delta t_0$ be the size of the coarsest time mesh that corresponds to the level $\ell{=}0$.
The size of the time mesh at  level $\ell \geq 1$ is given by $\Delta t_{\ell} {=} R^{-\ell} \Delta t_0$, where $R{>}1$ is a given integer constant. 

Assume that we are interested in estimating $\expt{g_{L}}$, and we are able to simulate correlated pairs, $(g_{\ell},g_{\ell{-}1})$ for $\ell=1,\ldots,L$. Then, the following unbiased Monte Carlo estimator of $\expt{g_{L}}$ uses $g_{L-1}$ as a control variate:
\begin{align*}
\tilde \mu_L &:= \frac{1}{M_L} \sum_{m_L=1}^{M_L} \LP g_{L}(\omega_{m_L}) - (g_{L{-}1}(\omega_{m_L}) - \expt{g_{L{-}1}})\RP\\
	   &\, =\expt{g_{L{-}1}} +  \frac{1}{M_L} \sum_{m_L=1}^{M_L}  (g_{L}- g_{L{-}1})(\omega_{m_L})\PERIOD
\end{align*}
Applying this idea recursively and taking into account the following telescopic decomposition: $\expt{g_L} = \expt{g_0} + \sum_{{\ell}=1}^{L} \expt{g_{\ell}-g_{\ell -1}}$,
we arrive at the multilevel Monte Carlo estimator of $\expt{g_L}$: 
\begin{align}\label{eq:mul}
\hat{\mu}_L := \frac{1}{M_0} \sum_{m_0=1}^{M_0}g_0(\omega_{m_0}) + \sum_{\ell=1}^L \frac{1}{M_{\ell}} \sum_{{m_{\ell}}=1}^{M_{\ell}} (g_{\ell} - g_{\ell-1})(\omega_{m_{\ell}}) \PERIOD
\end{align}

We have that $\hat{\mu}_L$ is unbiased, since $\expt{\hat{\mu}_L} {=} \expt{g_{L}}$.
The variance of $\hat{\mu}_L$ is given by 
$\var{\hat{\mu}_L} = \frac{\var{g_0}}{M_{0}} + \sum_{\ell=1}^L \frac{\var{g_{\ell}{-}{g_{\ell{-}1}}}}{M_{\ell}}$. Here, we are assuming independence among the batches between  levels.
For highly correlated pairs, $(g_{\ell},g_{\ell{-}1})$, we can expect, for the same computational work, that $\var{\hat{\mu}_L}$ is much less than the variance of the standard Monte Carlo estimator of $\expt{g_L}$.



\subsection{The Large Kurtosis Problem}
\label{sec:highk}
{Let us give a close examination of the problem of estimating $\var{g_{\ell}{-}{g_{\ell{-}1}}}$
for highly correlated pairs, $(g_{\ell},g_{\ell{-}1})$. This estimation is required to solve the optimization problem \eqref{eq:minworkprob}, that indicates how to choose the simulation parameters, particularly the number of simulated coupled paths for each pair of consecutive levels, $\seqof{M_{\ell}}{\ell=0}{L}$.

When $\ell$ becomes large, due to our coupling strategy developed in Section \ref{sec:couplingpaths}, we expect to obtain $g_{\ell}= g_{\ell{-}1}$ in most of our simulations, while observing differences only in a very small proportion of the simulated coupled paths.  

For the sake of illustration, let us assume that the random variable $\chi_{\ell}{:=}{g_{\ell}{-}{g_{\ell{-}1}}}$ takes values in the set $\{-1,0,1\}$, with respective probabilities $\{p_{\ell},1-2p_{\ell},p_{\ell}\}$, where $p_{\ell}$ goes to zero. The kurtosis of $\chi_{\ell}$ is by definition
\begin{align*}
\frac{\expt{\LP\chi_{\ell} - \expt{\chi_{\ell}} \RP^4}}{\LP\expt{\LP \chi_{\ell} - \expt{\chi_{\ell}}\RP^2} \RP^2}-3\PERIOD
\end{align*}
Simple calculations show that the kurtosis of $\chi_{\ell}$ is $(2p_{\ell})^{-1}$, and we
observe that $\chi_{\ell}^2 \sim \text{Bernoulli}(2p_{\ell})$. The maximum likelihood estimator of $2p_{\ell}$, $\hat{\theta}_{\ell}$, is the sample average of $M_{\ell}$ 
independent and identically distributed
(iid) values of $\chi_{\ell}^2$. The coefficient of variation of $\hat{\theta}_{\ell}$, defined as 
$({\var{\hat{\theta}_{\ell}}})^{1/2} (\expt{\hat{\theta}_{\ell}})^{-1}$, 
is $(2p_{\ell}M_{\ell})^{-1}$. Therefore, an accurate estimation of $p_{\ell}$ requires a sample of size 
\begin{align*}
M_{\ell}{\gg}(2p_{\ell})^{-1}{\to}\infty\PERIOD
\end{align*}
This lower bound on $M_{\ell}$ goes strongly against the spirit of the Multilevel Monte Carlo method, where $M_{\ell}$ should be a decreasing function of $\ell$. 

To overcome this difficulty, in Section \ref{sec:vl}, we developed a formula based on dual-weighted residuals. The technique of dual-weighted residuals can be motivated as follows: 
consider a process $\dbar X$, such that its position at time $s$, having departed from the state $x$, at a previous time $t$, is denoted as $\dbar X(s; t,x)$.
Notice that for $t{<}s{<}T$, we have that $\dbar X(T;t,x) = \dbar X(T;s,\dbar X(s;t,x))$. 
Let us define an auxiliary function $U(t,x):=g(\dbar X(T;t,x))$, where $g$ is an observable scalar function of the final state of the process $\dbar X$ that started from the state $x$ at the initial time, $t$. If $\bar X$ is a process approximating $\dbar X$, we want to have a computable approximation for $g(\bar X(T;0,x_0)) - g(\dbar X(T;0,x_0))$.
Consider a time mesh, $\{0{=}t_0,t_1,\ldots,t_N{=}T\}$, and define 
$\bar X_{t_{n}}{:=}\bar X(t_{n}; 0,x_0)$,  $\dbar X_{t_{n+1}}{:=}\dbar X(t_{n+1}; t_n,\bar X_{t_n})$ and $e_{n+1}:= \bar X_{t_{n+1}} {-} \dbar X_{t_{n+1}} $.
Observe that 
\begin{align*}
g(\bar X(T;0,x_0)) - g(\dbar X(T;0,x_0)) &= U(T,\bar X(T;0,x_0))-U(0,x_0)\\
&=\sum_{n=0}^{N-1} \LP U(t_{n+1},\bar X_{t_{n+1}}){-}U(t_{n},\bar X_{t_{n}})\RP\\
&=\sum_{n=0}^{N-1} \LP U(t_{n+1},\bar X_{t_{n+1}}){-}U(t_{n+1},\dbar X(t_{n+1}; t_n,\bar X_{t_n}))\RP\\
&=\sum_{n=0}^{N-1}\LP e_{n+1} \cdot \int_0^1 
\nabla_x U(t_{n+1}, \bar X_{t_{n+1}} {-} s e_{n+1}) ds \RP\\
&=\sum_{n=0}^{N-1} \LP e_{n+1} \cdot \nabla_x U(t_{n+1}, \bar X_{t_{n+1}}) + \norm{\nabla^2 U} \norm{e_{n+1}}^2+h.o.t.\RP\PERIOD
\end{align*}

We can now write a backward recurrence for the dual weights, $\seqof{\phi_{n}}{n=1}{N}$:
\begin{align*}
\phi_{n}&{:=}\nabla_x U(t_{n}, \bar X_{t_{n}}) = \partial_{ \bar X_{t_{n}}} g(\dbar X(T; t_{n}, \bar X_{t_{n}}))\\
&= \partial_{ \bar X_{t_{n}}} g(\dbar X(T; t_{n+1}, \bar X_{t_{n+1}}))\\
&= \partial_{ \bar X_{t_{n+1}}} g(\dbar X(T; t_{n+1}, \bar X_{t_{n+1}})) 
\displaystyle{\frac{\partial_{ \bar X_{t_{n+1}}}}{\partial_{ \bar X_{t_{n}}}}}\\
&= \nabla_x U(t_{n+1}, \bar X_{t_{n+1}})
\displaystyle{\frac{\partial_{ \bar X_{t_{n+1}}}}{\partial_{ \bar X_{t_{n}}}}}\\
&= \phi_{n+1}
\displaystyle{\frac{\partial_{ \bar X_{t_{n+1}}}}{\partial_{ \bar X_{t_{n}}}}}\\
\phi_N&{:=}\nabla g(\bar X(T;0,x_0))\PERIOD
\end{align*}

This reasoning evidently works for processes $\dbar X$ that are pathwise differentiable with respect to the initial condition. Our space state is in general a subset of the lattice, $\latt$, and for that reason, we can not directly apply this technique. 
In \cite{kkst}, the authors show how this dual-weighted residual technique can be adapted to the tau-leap case in regimes close to the mean field or to the Stochastic Langevin limit. 
In more general regimes, the formula \eqref{eq:varhatestimated}, which provides accurate estimates of $\var{g_{\ell}{-}{g_{\ell{-}1}}}$ in our numerical examples (see for instance Figure \ref{fig:effdec2} in Section \ref{sec:examples}), is promising but  more research is needed in this direction.
Specifically, in Section \ref{sec:vl}, the formula \eqref{eq:varhatestimated} is deduced from the conditional distribution of the local errors, $e_{n+1}|\mathcal{F}$, conditional on a sigma-algebra, $\mathcal{F}$, generated by the sequence, 
$\seqof{\bar X_{t_n}}{n=1}{N}$, and applying the tower properties of conditional expectation and conditional variance.
Similar comments apply to Formula \eqref{eq:EILw} regarding the weak error, $\expt{g(X(T))-g_L}$.}



\subsection{Outline of this Work}
In Section \ref{sec:couplingpaths}, we first show the main idea for coupling two tau-leap paths, which comes from a construction by Kurtz \cite{Kurtz82} for coupling two Poisson random variables. Then, inspired by the ideas of Anderson and Higham in \cite{Anderson2012}, we propose an algorithm for coupling two hybrid Chernoff tau-leap paths (see \cite{ourSL}). This algorithm uses four building blocks that result from the combination of the MNRM and the tau-leap methods.
In Section \ref{sec:MLMC}, we propose a novel hybrid MLMC estimator. Next, we introduce a global error decomposition; and finally, we develop formulae to efficiently estimate the variance of the difference of two consecutive levels {and to estimate the bias} based on dual-weighted residuals. 
{These estimates are particularly useful to addressing {the} large kurtosis problem, {described in Section \ref{sec:highk}}, that appears at the deeper levels and makes standard sample estimators too costly.} 
Next, in Section \ref{EEC}, we show how to control the three error components of the global error and how to obtain the parameters needed for computing the hybrid MLMC estimator to achieve a given tolerance with nearly optimal computational work. 
We also show that the computational complexity of our method is of order $\Ordo{TOL^{-2}}$.
In Section \ref{sec:examples}, the numerical examples illustrate the advantages of the hybrid MLMC method over the single-level approach presented in \cite{ourSL} and to the SSA. Section \ref{sec:conclusions} presents our conclusions and suggestions for future work.



\section{Generating Coupled Hybrid Paths}
\label{sec:couplingpaths}
In this section, we present an algorithm that generates coupled hybrid Chernoff tau-leap paths, which is an essential ingredient for the multilevel Monte Carlo estimator.
We first show how to couple two Poisson random variables and then we explain how we make use of the two algorithms presented in \cite{Anderson2012} as Algorithms 2 and 3 and two additional algorithms we developed to create an algorithm that generates coupled hybrid paths.

\subsection{Coupling Two Poisson Random Variables} 
\label{sec:couplingpoiss}
We motivate our coupling algorithm (Algorithm \ref{alg:coupled}) by first describing how to couple two Poisson random variables.
In our context, `coupling' means that we want to induce a correlation between them that is as strong as possible.
This construction was first proposed by Kurtz in \cite{Kurtz82}.
Suppose that we want to couple $\mathcal{P}_1(\lambda_1)$ and $\mathcal{P}_2(\lambda_2)$,
two Poisson random variables, with rates $\lambda_1$ and $\lambda_2$, respectively. 
Consider the following decompositions, 
\begin{align*}
\mathcal{P}_1(\lambda_1) &:= \mathcal{P}^*(\lambda_1 \wedge \lambda_2) + \mathcal{Q}_1(\lambda_1 -\lambda_1 \wedge \lambda_2) \\
\mathcal{P}_2(\lambda_2) &:= \mathcal{P}^*( \lambda_1 \wedge \lambda_2) + \mathcal{Q}_2(\lambda_2 - \lambda_1 \wedge \lambda_2)\COMMA
\end{align*}
where  $ \mathcal{P}^*(\lambda_1 \wedge \lambda_2)$, $\mathcal{Q}_1(\lambda_1 -\lambda_1 \wedge \lambda_2)$ and $\mathcal{Q}_2(\lambda_2 - \lambda_1 \wedge \lambda_2)$ are three independent Poisson random variables. Here, $\lambda_1 \wedge \lambda_2 :=\min\{\lambda_1,\lambda_2\}$.
Observe that at least one of the following vanishes: $\mathcal{Q}_1(\lambda_1 - \lambda_1 \wedge \lambda_2)$ and $\mathcal{Q}_2(\lambda_2 - \lambda_1 \wedge \lambda_2)$. This is because at least one of the rates is zero. Algorithm \ref{alg:coupled} implements these ideas.
Finally, note that, by construction, we have
\begin{align*}
\var{\mathcal{P}_1(\lambda_1) - \mathcal{P}_2(\lambda_2)} &= \var{\mathcal{Q}_1(\lambda_1 -\lambda_1 \wedge \lambda_2)  - \mathcal{Q}_2(\lambda_2 - \lambda_1 \wedge \lambda_2) }\\ 
&= \abs{\lambda_1 - \lambda_2} \nonumber \PERIOD
\end{align*}
However, if instead we consider making $\mathcal{P}_1(\lambda_1)$ and $\mathcal{P}_2(\lambda_2)$ independent,  then 
$$\var{\mathcal{P}_1(\lambda_1) {-} \mathcal{P}_2(\lambda_2)} = \lambda_1 + \lambda_2 \COMMA$$ which may be a large value even when $\lambda_1$ and $\lambda_2$ are close.

\subsection{Coupling Two Hybrid Paths}\label{sec:chp}
In this section, we describe how to generate two coupled hybrid Chernoff tau-leap paths, $\bar X$ and $\bar{\bar X}$, corresponding to two nested time discretizations, called coarse and fine, respectively.
Assume that the current time is $t$, and we know the states, $\bar X(t)$ and $\bar{\bar X}(t)$.
Based on this knowledge, we have to determine a method for each level. 
This method can be either the MNRM or the tau-leap one, determining 
 four possible combinations leading to four algorithms, B1, B2, B3 and B4, that we use as building blocks. Table \ref{tab:algs} summarizes them.

\begin{table}[h!]
\centering
\begin{tabular}{lcll}
  &Algorithm &at coarser mesh & at fine mesh \\ \noalign{\smallskip} \hline\noalign{\smallskip} 
B1 &(part of Algorithm \ref{alg:coupled})&TL & TL \\
B2& (Algorithm \ref{alg:auxilco}) &TL & MNRM \\
B3& " &MNRM & TL \\
B4& " &MNRM & MNRM
\end{tabular}
\bigskip
\caption{Building blocks for simulating two coupled hybrid Chernoff tau-leap paths.  
Algorithms B1 and B2 are presented as Algorithms 2 and 3 in \cite{Anderson2012}. 
Algorithm B3 can be directly obtained from Algorithm B2. 
Algorithm B4 is also based on Algorithm B2, but to produce MNRM steps, we update the propensities at the coarse level at the beginning of each time interval defined by the fine level.}
\label{tab:algs}
\end{table}

We note that the only case in which we use a Poisson random variates generator for the tau-leap method is in Algorithm B1. 
In Algorithms B2 and B3, the Poisson random variables are simulated by adding independent exponential random variables with the same rate, $\lambda$, until a given time final time $T$ is exceeded. The rate, $\lambda$, is obtained by freezing the propensity functions, $a$, at time $t$. More specifically, the Poisson random variates are obtained by using the MNRM repeatedly without updating the intensity.

We now briefly describe the Chernoff hybrid coupling algorithm, \ie,  Algorithm \ref{alg:coupled}. Given the current time, $t$, and the current state of the process at the coarse level, $\bar{X}(t)$, and the fine level, $\dbar{X}(t)$, this algorithm determines the next time point at which we run the algorithm (called time ``horizon''). 
To fix the idea, let us assume that, based on $\bar{X}(t)$, the one-step switching rule, \ie, Algorithm \ref{alg:sel}, chooses the tau-leap method at the coarse level, with the corresponding Chernoff step size, $\bar{\tau}$. As we mentioned, this $\bar{\tau}$ is the largest step size such that the probability that the process, in the next time step, takes a value outside $\latt$, is less than $\bar{\delta}$. This step size plus the current time, $t$, cannot be greater than the final time, $T$, and also cannot be greater than the next time discretization grid point in the coarse grid, $\bar{t}$, because the discretization error must be controlled. Taking the minimum of all those values, we obtain the next time horizon at the coarse grid, $\bar{H}$. Note that, if the chosen method is MNRM instead of tau-leap, we do not need to take into account the grid, and the next time horizon will be the minimum between the next reaction time and the final time, $T$.

\begin{figure}[h!]
\centering
\includegraphics[scale=1.00]{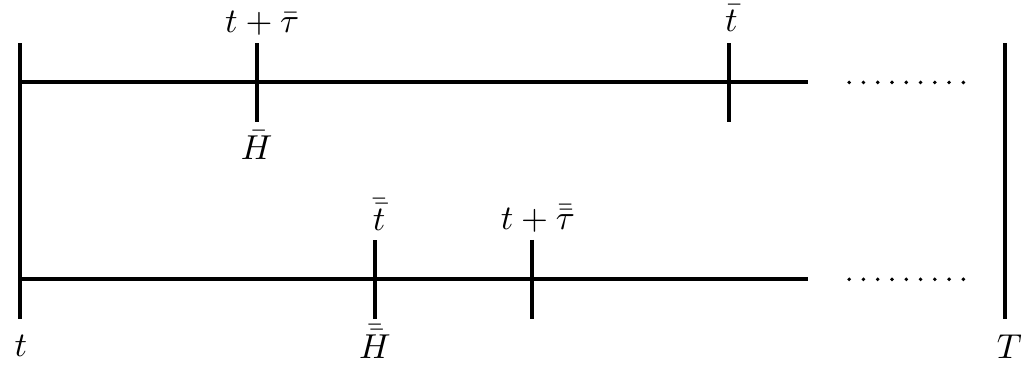}
\caption{This figure depicts a particular instance of the Chernoff hybrid coupling algorithm (Algorithm \ref{alg:coupled}), where $\bar{\tau}<\dbar{\tau}$. The synchronization horizon $H$, defined as $H{:=}\min\{\bar{H},\dbar{H}\}$, is equal to $\bar{H}$ in this case. Notice that $\bar{H}{:=}\min\{\bar{t},t+\bar{\tau},T\}$ and $\dbar{H}{:=}\min\{\dbar{t},t+\dbar{\tau},T\}$  }
\label{fig:hH}
\end{figure}

We now explain algorithm B1 (TL-TL). Assume that tau-leap is chosen at the coarse and at the fine level. We thus obtain two time horizons, one for the coarse level, $\bar{H}$, and another for the fine level, $\dbar{H}$. In this case, the global time horizon will be $H{:=}\min\{\bar{H},\dbar{H}\}$. Since the chosen method in both grid levels is tau-leap, we need to freeze the propensities at the beginning of the corresponding intervals. In the coarse case, during the interval $[t,\bar{H})$ (the propensities  are equal to $a(\bar{X}(t)) {=:}\bar{a}$), and in the fine case during the interval $[t,\dbar{H})$ (the propensities  are equal to $a(\dbar{X}(t)) {=:}\dbar{a}$). 
Suppose that $\bar{H}<\dbar{H}$ (see Figure \ref{fig:hH}). 

Then, we couple two Poisson random variables at time $t {=} \bar{H}$, using the idea described in Section \ref{sec:couplingpoiss}. 
When time reaches $\bar{H}$, the decision between which method to use (and the corresponding step size) at the coarse level must be made again. Note that the propensities of the process at the fine grid will be kept frozen until $\dbar{H}$. The case when $\bar{H}>\dbar{H}$ is analogous to the one we described, but the decisions on the method and step size are made at the finer level, when time reaches $\dbar{H}$. It can also be possible that $\bar{H}=\dbar{H}$. In that case, the decision between which method to use (and the corresponding step size) must be made at the coarse and at the fine level.

In the case of algorithm B2 (TL-MNRM), we assume that tau-leap is chosen at the coarse level, and MNRM at the fine level, obtaining two time horizons, one for the coarse level, $\bar{H}$, and another for the fine level, $\dbar{H}$. The only difference in how we determine the time horizons between algorithms B1 and B2 is that the time discretization grid points in the fine grid are not taken into account to determine $\dbar{H}$. Algorithm B2 is then applied until the simulation reaches $H{:=}\min\{\bar{H},\dbar{H}\}$. Suppose that $\dbar{H}<\bar{H}$. In this case, the process $\dbar{X}$ could take more than one step to reach $\dbar{H}$. At each step, the propensity functions $a(\dbar{X}(\cdot))$ are computed, but not the propensities for the coarse level, because in that case the tau-leap method is used. Note that the decision between which algorithm to use (B2 or another) is not made at those steps, but only when time reaches $\dbar{H}$. When time reaches $\dbar{H}$, the decision of which method to use (and the corresponding step size) at the fine level must be made again. In this case, the propensities at the coarse grid will be kept frozen until $\bar{H}$. The reasoning for the cases $\dbar{H}>\bar{H}$ and $\dbar{H}=\bar{H}$ are similar to before.

The other two cases, that is, B3 and B4, are the same as B2. The only difference resides is when to update the propensity values, $\bar{a}$ and $\dbar{a}$. See Algorithm \ref{alg:coupled} for more details.
As made clear in the preceding paragraphs,  the decision on which algorithm to use for a certain time interval is made only at the horizon points. 

\begin{rem}\label{rem:telescoping}[About telescoping]
{To ensure the telescoping sum property, 
the probability law of the hybrid process at level $\ell$ should be the same disregarding whether level $\ell$ is the finer in the pair $(\bar X_{\ell-1},\dbar X_{\ell})$ or the coarser in the pair $(\bar X_{\ell},\dbar X_{\ell+1})$.
For that reason, each process has its own next horizon as its decision points. See Figure \ref{fig:hH} showing the time horizons scheme and Figures \ref{fig:box2} and \ref{fig:box3} in Section \ref{sec:examples} to see that the telescoping sum property is satisfied by our hybrid coupling sampling scheme.}
\end{rem}

\section{Multilevel Monte Carlo Estimator and Global Error Decomposition}
\label{sec:MLMC}
In this section, we present the multilevel Monte Carlo estimator.
We first show the estimator and its properties and then we analyze and control the computational global error, which is decomposed into three error components: the discretization error, the global exit error, and the Monte Carlo statistical error. We give upper bounds for each one of the  three components. 

\subsection{The MLMC Estimator}\label{sec:nontel}
In this section, we discuss and implement a variation of the multilevel Monte Carlo estimator \eqref{eq:mul}  for the hybrid Chernoff tau-leap case. 
The main ingredient of this section is Algorithm \ref{alg:coupled}, which generates coupled hybrid paths at levels $\ell{-}1$ and $\ell$.
Let us now introduce some notation. 
Let $A_\ell$ be the event in which the $\bar X_{{\ell}}$-path arrived at the final time, $T$, without exiting the state space of $X$. Let $\indicator{A}$, be the indicator function of an arbitrary set, $A$.
Finally, $g_{\ell}:= g({\bar X}_{\ell}(T))$ was defined in Section \ref{MLMCintro}.

Consider the following telescopic decomposition:
\begin{align*}
\expt{g_L \indicator{A_L}} = \expt{g_0 \indicator{A_0}} {+} \sum_{\ell=1}^L \expt{g_\ell\indicator{A_\ell}-g_{\ell-1}\indicator{A_{\ell-1}}} \COMMA
\end{align*} 
which motivates the definition of our MLMC estimator of $\expt{g(X(T))}$,
\begin{align}\label{MLMCest}
\mathcal{M}_L := \frac{1}{M_0} \sum_{m=1}^{M_0}g_0\indicator{A_0}(\omega_{m,0}) + \sum_{\ell=1}^L \frac{1}{M_{\ell}} \sum_{m=1}^{M_{\ell}} [g_{\ell} \indicator{A_\ell} - g_{\ell-1} \indicator{A_{\ell-1}}](\omega_{m,\ell}) \PERIOD
\end{align}


\subsection{Global Error Decomposition} 
\label{calibration}
In this section, we define the computational global error, $\EE_L$, and show how it can be naturally decomposed into three components: the discretization error, $\EE_{I,L}$, and the exit error, $\EE_{E,L}$, both coming from the tau-leap part of the hybrid method and the Monte Carlo statistical error, $\EE_{S,L}$. Next, we show how to model and control the global error, $\EE_L$, giving upper bounds for each one of the  three components. 
We define the computational global error, $\EE_L$, as
\begin{equation*}
\EE_L := \expt{g(X(T))} - \mathcal{M}_L\PERIOD
\end{equation*}
Now, consider the following decomposition of $\EE_L$:
\begin{align*}
\expt{g(X(T))} - \mathcal{M}_L &= \expt{g(X(T))(\indicator{A_L}+\indicator{A^c_L})} - \expt{g_L \indicator{A_L}} + \expt{g_L \indicator{A_L}} - \mathcal{M}_L\nonumber \\
&= \underbrace{\expt{g(X(T)) \indicator{A^c_L}}}_{=:\EE_{E,L}} 
 + \underbrace{\expt{ \LP g(X(T)) {-} g_L \RP \indicator{A_L}}}_{=:\EE_{I,L}} + \underbrace{\expt{g_L \indicator{A_L}}{-}\mathcal{M}_L}_{=:\EE_{S,L}}\PERIOD 
\end{align*}

We show in \cite{ourSL} that by choosing adequately the one-step exit probability bound, $\delta$, the exit error, $\EE_{E,L}$, satisfies 
$|\EE_{E,L}| \leq |\expt{g(X(T))}|\, \prob{A^c_L}\leq TOL^2$. 
An efficient procedure for accurately estimating $\EE_{I,L}$ in the context of the tau-leap method is described in \cite{kt}. 
We adapt this method in Algorithm \ref{alg:weakerror} for estimating the weak error in the hybrid context. A brief description follows. For each hybrid path,  
$\seqof{\bar X_{\ell}(t_{n,\ell},\bar \omega)}{n=0}{N(\bar{\omega})}$, we define  the sequence of dual weights $\seqof{\varphi_{n,\ell}(\bar \omega)}{n=1}{N(\bar{\omega})}$ backwards as follows (see Section \ref{sec:highk}):
\begin{align}\label{eq:phisback}
\varphi_{N(\bar{\omega}),\ell} &:= \nabla g(\bar X_{\ell}(t_{N(\bar{\omega}),\ell},\bar \omega))\\
\varphi_{n,\ell} &:= \LP Id +  {\Delta t_{n,\ell}} \JAC_a^T(\bar X_{\ell}(t_{n,\ell},\bar \omega)) \,\nu^T\RP \varphi_{n+1,\ell}, \quad n=N(\bar{\omega}){-}1,\ldots,1\COMMA\nonumber
\end{align}
where $\Delta t_{n,\ell}{:=}t_{n+1,\ell}{-}t_{n,\ell}$, $\nabla$ is the gradient operator and $\JAC_a(\bar X_{\ell}(t_{n,\ell},\bar \omega)){\equiv}[\partial_i a_j(\bar X_{\ell}(t_{n,\ell},\bar \omega))]_{j,i}$ is the Jacobian matrix of the propensity function, $a_j$, for $j{=}1\ldots J$ and $i{=}1\ldots d$. 
According to this method, $\EE_{I,L}$ is approximated by $\avg{\EE_{I,L}(\bar\omega)}{\cdot}$, where 
{
\begin{equation}\label{eq:EILw}
\EE_{I,L}(\bar \omega) := \sum_{n=1}^{{N(\bar{\omega})}}  \LP\frac {\Delta t_{n,L}}{2}  \indicator{TL}(n)  \sj ( \varphi_{n,L} \cdot \nu_j) \Delta a_{j,n}\RP(\bar \omega)\COMMA
\end{equation}}
$\avg{X}{M}{:=}\frac 1 M \sum_{m=1}^M X(\omega_m)$, and, $\svar{X}{M}{:=}\avg{X^2}{M}- \avg{X}{M}^2$ denote
the sample mean and the sample variance of the random  variable, $X$, respectively.
Here,  $\Delta a_{j,n}(\bar \omega){:=} a_j(\bar X_{L}(t_{n+1,\ell},\bar \omega)){-}a_j(\bar X_{L}(t_{n,\ell},\bar \omega))$, $\indicator{TL}(n){=}1$ if and only if, at time $t_{n,\ell}$, the tau-leap method was used, and we denote by $Id$ the $d\times d$ identity matrix. 
{
\begin{rem}[Computational cost of dual computations]
It is easy to see that the computational cost per path of the dual computations in  \eqref{eq:phisback} is comparable, and possibly smaller than the hybrid path.
Indeed,  no new random variables, especially Poisson ones, which are the most computationally expensive in the forward simulation, need to be sampled and no coupling between levels is needed.
Moreover, we use \eqref{eq:phisback} only to determine the discretisation parameters for the actual run; so \eqref{eq:phisback} is thus used only in a fraction of the realisations.
\end{rem}
}

The variance of the statistical error, $\EE_{S,L}$, is given by $\sum_{\ell=0}^L \frac{\mV{\ell}}{M_{\ell}}$,
where $\mV{0}:=\var{g_{0}\indicator{A_0}}$ and 
$\mV{\ell}:=\var{ g_{\ell} \indicator{A_{\ell}} -g_{\ell-1} \indicator{A_{\ell{-}1}}},\,\,\ell \geq 1$. 
In the next subsection, we show how to estimate $\mV{\ell}$ {efficiently using the duals from \eqref{eq:phisback}.}



\subsection{Dual-weighted Residual Estimation of $\mV{\ell}$}\label{sec:vl}

Here, we derive the formula \eqref{eq:varhatestimated} for estimating the variance, $\mV{\ell},\,\ell\geq 1$. 
It is based on dual-weighted local errors arising from two consecutive tau-leap approximations of the process, $X$. 
For {each} level $\ell \geq 1$, the formula estimates $\mV {\ell}$ with much smaller statistical error than the standard sample estimator, which is 
seriously affected by the large kurtosis present at the deepest levels (see Section \ref{sec:highk}). 

Let us introduce some notation:
\begin{align*} 
f_{j,n}&:=(\varphi_{n+1}\cdot \nu_j)\COMMA \\
\mu_{j,n}&:=\frac{\Delta t_n}{2} \sum_i(\nabla a_j(x_n)\cdot \nu_i)a_i(x_n)\COMMA \\ 
\bar{\mu}_{j,n}&:=\frac{\Delta t_n}{2} \sum_i|(\nabla a_j(x_n)\cdot \nu_i)|a_i(x_n)\COMMA \\
\sigma^2_{j,n}&:=\frac{\Delta t_n}{2} \sum_i(\nabla a_j(x_n)\cdot \nu_i)^2 a_i(x_n)\COMMA \\
m_{j,n}&:=\min\{ \bar{\mu}_{j,n}, \sqrt{\mu_{j,n}^2+\sigma^2_{j,n}}\}\COMMA \\
q_{j,n}&:=\frac{\mu_{j,n}}{\sigma_{j,n}}\COMMA \\
p_{j,n}&:=\Phi(-q_{j,n})\COMMA \\
\tilde{\mu}_{j,n}&:=\mu_{j,n}(1{-}2p_{j,n})\COMMA \\
\tilde{\sigma}_{j,n}&:=\sqrt{\frac{2}{\pi}}\sigma_{j,n} \exp({-q_{j,n}^2/2}) \PERIOD
\end{align*}
Here, $\Phi(x)$ is the cumulative distribution function of a standard Gaussian random variable. 
We define our dual-weighted estimator of $\mV{\ell}$ as
\begin{align}
\label{eq:varhatestimated}
\hat{\mathcal{V}}_\ell&:=
\svar{\sum_n \indicator{TL}(n) \frac{\Delta t_n}{2} \sum_j f_{j,n} \mu_{j,n}} {M_{\ell}} \\\nonumber
&+ \avg{\sum_n \indicator{TL}(n) \frac{(\Delta t_n)^3}{8} \sum_{j,j'}f_{j,n}f_{j',n}
\sum_i (\nabla a_j(x_n)\cdot \nu_i)(\nabla a_{j'}(x_n)\cdot \nu_i)a_i(x_n)}{M_{\ell}}\\\nonumber
&+ \mathcal{A}\left(\sum_n \indicator{TL}(n) \frac{\Delta t_n}{2} \sum_{j}f_{j,n}^2 
\left( \indicator{G_n}(\tilde{\mu}_{j,n}+\tilde{\sigma}_{j,n}) + \indicator{G_n^c}\,m_{j,n} \right),{M_{\ell}}\right)\nonumber\COMMA
\end{align}
where $\indicator{G_n}{=}1$ if and only if $a_j(x_n) \frac{\Delta t_n}{2}{>}c$ for all $j \in \{1,\ldots,J\}$, where $c$ is a positive user-defined constant.

First, notice that $\mV {\ell}$ could be a very small positive number. In fact, in our numerical experiments, we observe that the standard Monte Carlo sample estimation of this quantity turns out to be computationally infeasible due to the huge number of simulations required to stabilize its coefficient of variation. For this reason, we initially consider the following dual-weighted approximations:
\begin{align}\label{eq:moments}
\expt{\ g_{\ell}-g_{\ell-1} } &\approx \expt{\sum_{n} \varphi_{n+1,\ell-1} \cdot e_{n+1,\ell-1}} \COMMA\\
\var{\ g_{\ell}-g_{\ell-1} } &\approx \var{\sum_{n} \varphi_{n+1,\ell-1} \cdot e_{n+1,\ell-1}}\COMMA\nonumber
\end{align}
where  $\seqof{\varphi_{n+1,\ell-1}}{n=0}{N(\bar \omega)-1}$,  defined in \eqref{eq:phisback}, is a sequence of dual weights computed backwards from a simulated path, $\seqof{\bar X_{\ell}(t_{n,\ell-1})}{n=1}{N(\bar \omega)}$, and 
{the sequence of local errors, $\seqof{e_{n+1,\ell-1}}{n=0}{N(\bar \omega)-1}$, defined in \eqref{eq:locerrRP}, is the subject of the next subsection}.

\subsubsection*{Defining the Sequence of Local Errors}\label{sec:hypTL}
For simplicity of analysis, we make two assumptions:  
i) the time mesh associated with the level, ${\ell}$, is obtained by halving the intervals of the level $\ell{-}1$;
ii) we perform the tau-leap at both levels without considering the Chernoff bounds described in Section \ref{chernoff_onestep}.

Let $\bar X$ and $\bar{\bar X}$ be two tau-leap approximations of $X$ based on two consecutive grid levels, for instance, 
$\bar X {:=} \bar X_{\ell{-}1}$ and $\bar{\bar X} {:=} \bar X_{\ell}$. 
Consider two consecutive time-mesh points for $\bar X$, $\{t_n, t_{n+1}\}$, and three consecutive time-mesh points for $\bar{\bar X}$, $\{t_n, (t_n{+}t_{n+1})/2, t_{n+1}\}$. 
Let $\bar X$ and $\bar{\bar X}$ start from $x_n$ at time $t_n$.  

The first step for coupling $\bar X$ and $\bar{\bar X}$ is to define 

\begin{align}
\label{eq:dxxb} \bar X_{n+1} &:= x_n + \sum_j \nu_j \mathcal{Y}_{j,n}(a_j(x_n)\Delta t_n)\COMMA\\
\label{eq:dxxbZ} Z_{n+1} &:= x_n + \sum_j \nu_j \mathcal{Q}_{j,n}(a_j(x_n)\frac{\Delta t_n}2)\COMMA\\
\bar{\bar X}_{n+1} &:= Z_{n+1} + \sum_j \nu_j \mathcal{R}_{j,n}(a_j(Z_{n+1})\frac{\Delta t_n}2) \nonumber\COMMA
\end{align}
where $\{\mathcal{Y}_{j,n}\}_{j=1}^J \cup \{\mathcal{Q}_{j,n}\}_{j=1}^J \cup \{\mathcal{R}_{j,n}\}_{j=1}^J$ are 
Poisson random variables.
To couple the  $\bar X$ and $\bar{\bar X}$ processes, we first decompose $\mathcal{Y}_{j,n}(a_j(x_n)\Delta t_n)$ as 
the sum of two independent Poisson random variables, $\mathcal{Q}_{j,n}(a_j(x_n)\frac{\Delta t_n}2)+\mathcal{Q'}_{j,n}(a_j(x_n)\frac{\Delta t_n}2)$.
As a consequence, $\bar X$ and $\bar{\bar X}$ coincide in the closed interval $[t_n,(t_n{+}t_{n+1})/2]$.
By applying this decomposition in \eqref{eq:dxxb}, we obtain
\begin{align}\label{eq:qq}
\bar X_{n+1} &= x_n + \sum_j \nu_j\mathcal{Q}_{j,n}(a_j(x_n)\frac{\Delta t_n}2)+  \sum_j \nu_j\mathcal{Q}_{j,n}'(a_j(x_n)\frac{\Delta t_n}2)\COMMA\\
\bar{\bar X}_{n+1}&= x_n + \sum_j \nu_j\mathcal{Q}_{j,n}(a_j(x_n)\frac{\Delta t_n}2) + \sum_j \nu_j\mathcal{R}_{j,n}(a_j(Z_{n+1})\frac{\Delta t_n}2)\PERIOD\nonumber
\end{align}

The second step for coupling $\bar X$ and $\bar{\bar X}$, according to \cite{Anderson2012}, is as follows:
let $m_j := \min\{a_j(x_n),a_j(Z_{n+1})\}$, $c_j := a_j(x_n) - m_j$ and $f_j := a_j(Z_{n+1}) - m_j$.
Notice that for each $j$, either $c_j$ or $f_j$ is zero (or both).

Now, consider the following decompositions:
\begin{align}\label{eq:qr}
\mathcal{Q}_{j,n}'(a_j(x_n)\frac{\Delta t_n}2) &= \mathcal{P}_{j,n}'(m_j\frac{\Delta t_n}2) + \mathcal{P}_{j,n}''(c_j\frac{\Delta t_n}2)\COMMA\\
\mathcal{R}_{j,n}(a_j(x_n)\frac{\Delta t_n}2)  &= \mathcal{P}_{j,n}'(m_j\frac{\Delta t_n}2) + \mathcal{R}_{j,n}'(f_j\frac{\Delta t_n}2)\COMMA\nonumber
\end{align}
where  $\mathcal{P}'$, $\mathcal{P}''$ and $\mathcal{R}'$ are independent Poisson random variables.

By substituting \eqref{eq:qr} into \eqref{eq:qq}, we define the local error, $e_{n+1,\ell-1}$, as
\begin{align}\label{eq:locerrRP}
e_{n+1,\ell-1} &:= \bar{\bar X}_{n+1} - \bar X_{n+1}\\
              &= \sum_j \nu_j \LP\mathcal{R}_{j,n}'(f_j\frac{\Delta t_n}2) -\mathcal{P}_{j,n}''(c_j\frac{\Delta t_n}2)\RP\nonumber\\
              &= \sum_j \nu_j \LP\mathcal{R}_{j,n}'(\Delta a_{j,n}\frac{\Delta t_n}2)\indicator{\{\Delta a_{j,n}>0\}} -\mathcal{P}_{j,n}''(-\Delta a_{j,n}\frac{\Delta t_n}2)\indicator{\{\Delta a_{j,n}<0\}}\RP\COMMA\nonumber
\end{align}
where $\Delta a_{j,n}:= a_j(Z_{n+1}){-}a_j(x_n)$ and $Z_{n+1}$ is defined in \eqref{eq:dxxbZ}.
Note that in \eqref{eq:locerrRP} not only are $\mathcal{R}_{j,n}'$ and $\mathcal{P}_{j,n}''$  random variables, but  
$\Delta a_{j,n}$ is also random because it depends on the random variables $\seqof{\mathcal{Q}_{j,n}}{j=1}{J}$. Also note that all the mentioned random variables are independent.
%

\subsubsection*{Conditioning}
At this moment, it is convenient to recall the tower properties of the conditional expectation and the conditional variance: given a random variable, $X$, and a sigma algebra, $\mathcal{F}$, defined over the same probability space, we have 
\begin{align}
\expt{X}&=\expt{\expt{X\SEP\mathcal{F}}}\COMMA\nonumber\\
\label{eq:towersv}
\var{X}&=\var{\expt{X\SEP\mathcal{F}}}+\expt{\var{X\SEP\mathcal{F}}}\PERIOD
\end{align}

Hereafter, we fix $\ell$ and, for the sake of brevity, omit it as a subindex.

Applying  \eqref{eq:towersv} to $\sum_n \varphi_{n+1} \cdot e_{n+1}$ and conditioning on $\mathcal{F}$,  we obtain
\begin{align*}
\var{\sum_n \varphi_{n+1} \cdot e_{n+1}}&=\var{\expt{\sum_n \varphi_{n+1} \cdot e_{n+1}\Big|\mathcal{F}}}+ \expt{\var{\sum_n \varphi_{n+1} \cdot e_{n+1}\Big|\mathcal{F}}}\nonumber\\
&= \var{\sum_n \expt{ \varphi_{n+1} \cdot e_{n+1}\SEP\mathcal{F}}}+ \expt{\sum_n \var{ \varphi_{n+1} \cdot e_{n+1}\SEP\mathcal{F}}}\PERIOD
\end{align*}
The main idea is to generate $M_{\ell}$ Monte Carlo paths, $\seqof{\bar X_{\ell}(t_n;\bar \omega)}{n=1}{N(\bar \omega)}$, and to estimate 
$\var{\sum_n \varphi_{n+1} \cdot e_{n+1}}$ using
\begin{align}
\label{eq:varhat}
\hat{\mathcal{V}}_\ell&:=\svar{\underbrace{\sum_n \expt{ \varphi_{n+1} \cdot e_{n+1}\SEP\mathcal{F}}(\bar \omega)}_{S_e(\bar\omega)}}{M_{\ell}}+ \avg{\underbrace{\sum_n \var{ \varphi_{n+1} \cdot e_{n+1}\SEP\mathcal{F}}(\bar \omega)}_{S_v(\bar\omega)}}{M_{\ell}}\PERIOD
\end{align}

To avoid nested Monte Carlo calculations, we develop exact and approximate formulas for computing 
$\expt{ \varphi_{n+1} \cdot e_{n+1}\SEP\mathcal{F}}$ and $\var{ \varphi_{n+1} \cdot e_{n+1}\SEP\mathcal{F}}$. To derive those formulas, we consider
 a sigma-algebra, $\mathcal{F}$, such that $\seqof{\varphi_n(\bar \omega)}{n=1}{N(\bar \omega)}$, conditioned on $\mathcal{F}$, is deterministic, \ie,  $\seqof{\varphi_n(\bar \omega)}{n=1}{N(\bar \omega)}$ is measurable with respect to $\mathcal{F}$. In this way, the only randomness in $\expt{ \varphi_{n+1} \cdot e_{n+1}\SEP\mathcal{F}}$ and $\var{ \varphi_{n+1} \cdot e_{n+1}\SEP\mathcal{F}}$ comes from the local errors, $\seqof{e_n}{n=1}{N(\bar \omega)}$. 

\subsubsection*{Conditional Local Error Representation} 
In this section, we derive a local error representation that takes into account the fact that 
the dual is computed backwards and the distribution of the local errors that is relevant to our calculations is therefore not exactly the one given by \eqref{eq:locerrRP}, but the distribution given by \eqref{eq:tclerr}. 

Consider the sequence $\seqof{\bar X_n}{n=0}{N(\bar \omega)}$ defined in \eqref{eq:dxxb}.
For fixed $n$, define $\mathcal{F}_n$ as the sigma-algebra  
\begin{align*}
\mathcal{F}_n := \sigma \LP \seqof{\mathcal{Y}_{j,k}(a_j(x_k)\Delta t_k)}{j=1,\dots,J,\, k=1,\ldots,n}{} \RP\COMMA
\end{align*}
\ie, the information we obtain by observing the randomness used to generate $\bar X_{n+1}$ from $x_0$. 
Motivated by dual-weighted expansions \eqref{eq:moments}, we want to express the local error representation \eqref{eq:locerrRP} conditional on $\mathcal{F}{:=}\mathcal{F}_{N(\bar \omega)}$. 

At this point, it is convenient to remember a key result for building Poissonian bridges. If $X_1$ and $X_2$ are two independent Poisson random variables with parameters $\lambda_1$ and $\lambda_2$, respectively, we have that $X_1\SEP X_1+X_2{=}k$ is a binomial random variable with parameters $k$ and $\lambda_1/(\lambda_1{+}\lambda_2)$. 

Applying this observation to 
the decomposition $\mathcal{Y}_{j,n}(a_j(x_n)\Delta t_n) {=} \mathcal{Q}_{j,n}(a_j(x_n)\frac{\Delta t_n}2)+\mathcal{Q'}_{j,n}(a_j(x_n)\frac{\Delta t_n}2)$, we conclude that the conditional distribution of $\mathcal{Q}_{j,n}(a_j(x_n)\frac{\Delta t_n}2)$ given $\mathcal{F}_n$, \ie, $\mathcal{Q}_{j,n}(a_j(x_n)\frac{\Delta t_n}2)\SEP\mathcal{F}_n$, is binomial with parameters $\mathcal{Y}_{j,n}$ and $1/2$.

Define now the sigma-algebra, $\mathcal{G}_n$, as
\begin{align*}
\mathcal{G}_n{:=}\sigma\LP\seqof{\mathcal{Q}_{j,n}(a_j(x_n)\frac{\Delta t_n}2)\SEP \mathcal{F}_n}{j=1}J\RP\PERIOD
\end{align*}
Applying the same argument to $\mathcal{P}_{j,n}''$, defined in \eqref{eq:qr}, we conclude that 
\begin{align*}
\mathcal{P}_{j,n}''\SEP\{\mathcal{F}_n, \mathcal{G}_n\} \sim \text{binomial}\LP\mathcal{Y}_{j,n}{-}\mathcal{Q}_{j,n},\frac{c_j}{a_j(x_n)}\RP\PERIOD
\end{align*}
From the definition of $Z_{n+1}{=} x_n {+} \sum_j \nu_j \mathcal{Q}_{j,n}$ in \eqref{eq:dxxbZ},  we conclude that
\begin{align*}
\mathcal{R}_{j,n}'\SEP \mathcal{G}_n \sim \text{Poisson}\LP(a_j(Z_{n+1})-m_j)\frac{\Delta t_n}{2}\RP\PERIOD
\end{align*}
Notice that, by construction, $\mathcal{P}_{j,n}''\SEP\{\mathcal{F}_n, \mathcal{G}_n\}$ and 
$\mathcal{R}_{j,n}'\SEP \mathcal{G}_n$ are independent random variables.
Since $c_j{=}-\Delta a_{j,n}\indicator{\{\Delta a_{j,n}<0\}}$ and 
$a_j(Z_{n+1}){-}m_j{=}\Delta a_{j,n}\indicator{\{\Delta a_{j,n}\geq0\}}$, we can express the conditional local error as
\begin{align}\label{eq:tclerr}
e_{n+1}&\SEP\{\mathcal{F}_n, \mathcal{G}_n\}=\\
 &\sum_j \nu_j \LP 
\mathcal{R}_{j,n}'\LP\Delta a_{j,n}\frac{\Delta t_n}{2}\RP\indicator{\{\Delta a_{j,n}\geq 0\}} {-}
\mathcal{P}_{j,n}''\LP\mathcal{Y}_{j,n}{-}\mathcal{Q}_{j,n}, \frac{{-}\Delta a_{j,n}}{a_j(x_n)}\RP\indicator{\{\Delta a_{j,n}<0\}}
\RP\nonumber
\end{align}
in the distribution sense.
For instance, we can easily compute the expectation of $e_{n+1}\SEP\{\mathcal{F}_n, \mathcal{G}_n\}$ as follows:
\begin{align*}
\expt{e_{n+1}\SEP\{\mathcal{F}_n, \mathcal{G}_n\}} = 
\sum_j \nu_j  \Delta a_{j,n}\LP 
 \frac{\Delta t_n}{2}\indicator{\{\Delta a_{j,n}\geq 0\}} +
 \frac{\mathcal{Y}_{j,n}{-}\mathcal{Q}_{j,n}}{a_j(x_n)}\indicator{\{\Delta a_{j,n}<0\}}\RP\PERIOD
\end{align*}
Taking into account that the joint distribution of $\seqof{\mathcal{Q}_{j,n}}{j=1}{J}\SEP \mathcal{F}_n$ is given by
\begin{align*}
\prob{\seqof{\mathcal{Q}_{j,n}=q_{j,n}}{j=1}{J}\SEP \mathcal{F}_n} = 2^{-\sum_j \mathcal{Y}_{j,n}}\prod_{j=1}^J \frac{{\mathcal{Y}_{j,n}}!}{ {q_{j,n}}! (\mathcal{Y}_{j,n} {-} q_{j,n})!},\quad 0\leq q_{j,n}\leq \mathcal{Y}_{j,n}\COMMA
\end{align*}
we can exactly compute  the expected value and the variance of $v_{n+1} \cdot e_{n+1}\SEP\mathcal{F}_n$ for any given deterministic vector, $v_{n+1}$. 
Notice that given $\mathcal{F}$, the sequence $\seqof{\bar X_n}{n=0}{N(\bar \omega)}$ is deterministic and, as a consequence, the sequence 
$\seqof{\varphi_n}{n=1}{N(\bar \omega)} \SEP \mathcal{F}$ is also a deterministic sequence of vectors.
We can thus compute 
\begin{align}\label{eq:thetwo}
\expt{\sum_n \varphi_{n+1}\cdot e_{n+1}\SEP\mathcal{F}} \text{ and } \var{\sum_n \varphi_{n+1}\cdot e_{n+1}\SEP\mathcal{F}}
\end{align}
exactly and proceed as stated at the beginning of this section.
However, trying to develop computable expressions from \eqref{eq:tclerr} has two main disadvantages: 
i) it may lead to computationally demanding procedures, especially for systems with many reaction channels or in regimes with high activity;
ii) it may be affected by the variance associated with the randomness in $\mathcal{F}_n$ and  $\mathcal{G}_n$.

\subsubsection*{Deriving a Formula for $\hmV{\ell}$}
In this section, we derive the formula \eqref{eq:varhatestimated}. 
Our goal is to find computable approximations of \eqref{eq:thetwo}, where  the underlying sigma-algebra, $\mathcal{F}$, is just the information gathered by observing the coarse path, $\bar X$.
This means that our formula should not depend explicitly on the knowledge of the random variables that generate $\mathcal{F}_n$ and  $\mathcal{G}_n$.   
At this point, it is important to recall the comments in Section \ref{sec:vl}; that is, the sequence 
$\seqof{\varphi_n(\bar \omega)}{n=1}{N(\bar \omega)}$ is measurable with respect to $\mathcal{F}$.
This implies that, for all $n$, $\varphi_{n+1}$ is independent of $\mathcal{G}_n$.
Hereafter, for notational convenience, we omit writing explicitly the conditioning on $\mathcal{F}$ in our formulae. 

It turns out that the leading order terms of the conditional moments obtained from \eqref{eq:tclerr} 
are essentially the same as those computed from \eqref{eq:locerrRP}. We will then derive \eqref{eq:varhatestimated} from \eqref{eq:locerrRP}. Using the notation from Section \ref{sec:vl}, we have that
\begin{align*}
(\varphi_{n+1} \cdot e_{n+1}) =  \sum_j f_{j,n}\LP \mathcal{R}_{j,n}'(\Delta a_{j,n}\frac{\Delta t_n}2)\indicator{\{\Delta a_{j,n}>0\}} -\mathcal{P}_{j,n}''(-\Delta a_{j,n}\frac{\Delta t_n}2)\indicator{\{\Delta a_{j,n}<0\}}\RP\PERIOD
\end{align*}
By the tower property, we obtain 
\begin{align*}
\expt{(\varphi_{n+1} \cdot e_{n+1})}=\expt{\expt{(\varphi_{n+1} \cdot e_{n+1})\SEP \mathcal{G}_n}} =  \frac{\Delta t_n}{2}\sum_j f_{j,n} \expt{\Delta a_{j,n}}\PERIOD
\end{align*}
Now let us consider the first-order Taylor expansion: 
\begin{align*}
\Delta a_{j,n} &:= a_j(x_n + \sum_i \nu_i\,\mathcal{Q}_{i,n}(a_i(x_n)\Delta t_n/2))-a_j(x_n)\\\nonumber
               &\approx (\nabla a_j(x_n)\cdot \sum_i \nu_i\,\mathcal{Q}_{i,n}(a_i(x_n)\Delta t_n/2))\\\nonumber
               &= \sum_i (\nabla a_j(x_n)\cdot \nu_i)\mathcal{Q}_{i,n}(a_i(x_n)\Delta t_n/2)\PERIOD
\end{align*}
Since $\mathcal{Q}_{i,n}(a_i(x_n)\Delta t_n/2) \sim \text{Poisson}(a_i(x_n)\Delta t_n/2)$, we have that
$\expt{\Delta a_{j,n}} = \mu_{j,n}$ and  $\var{\Delta a_{j,n}} = \sigma^2_{j,n}$.
Thus, 
\begin{align*}
\expt{(\varphi_{n+1} \cdot e_{n+1})}\approx \frac{\Delta t_n}{2} \sum_j f_{j,n}\,\mu_{j,n}\PERIOD
\end{align*}
Now, we use again the tower property for the variance:
\begin{align*}
\var{(\varphi_{n+1} \cdot e_{n+1})} = \var{\expt{(\varphi_{n+1} \cdot e_{n+1})\SEP \mathcal{G}_n}} + \expt{\var{(\varphi_{n+1} \cdot e_{n+1})\SEP \mathcal{G}_n}}\PERIOD
\end{align*}
We then immediately obtain
\begin{align*}
\var{\expt{(\varphi_{n+1} \cdot e_{n+1})\SEP \mathcal{G}_n}} &\approx \frac{(\Delta t_n)^3}{8}\sum_{j,j'} f_{j,n}f_{j',n}
\sum_i(\nabla a_j(x_n)\cdot \nu_i)(\nabla a_{j'}(x_n)\cdot \nu_i)a_i(x_n)\COMMA\\\nonumber
\expt{\var{(\varphi_{n+1} \cdot e_{n+1})\SEP \mathcal{G}_n}}&\approx \frac{\Delta t_n}{2} \sum_j f^2_{j,n}\expt{\Delta a_{j,n}\, \text{sgn}(\Delta a_{j,n})}\PERIOD
\end{align*}

Let us consider the case where $a_i(x_n)\Delta t_n/2$ is large enough for all $i$. 
It is well known that a Poisson random variable, $\mathcal{Q}(\lambda)$, is well approximated by a Gaussian random variable, $N(\lambda,\lambda)$, for moderate values of $\lambda$, say $\lambda{>}10$.
Since $\mathcal{Q}_{i,n}(a_i(x_n)\Delta t_n/2) \sim \text{Poisson}(a_i(x_n)\Delta t_n/2)$,
we have that, when $a_i(x_n)\Delta t_n/2$ is large enough for all $i$,  $\Delta a_{j,n} \approx N(\mu_{j,n},\sigma^2_{j,n})$.  
Consider a Gaussian random variable $Z$ with parameters $\mu$ and $\sigma^2{>0}$. Then,
\begin{align}\label{eq:dedu}
\expt{(\mu+\sigma Z) \indicator{\{\mu+\sigma Z>0\}}} &= \mu \prob{\mu+\sigma Z>0}+\frac{\sigma}{\sqrt{2\pi}} \int_{-\mu/\sigma}^{+\infty} z\exp{\left(-z^2/2\right)}\,dz\\ \nonumber
&= \mu(1-\Phi(-\mu/\sigma))+ \frac{\sigma}{\sqrt{2\pi}}\exp{\left(-(\mu/\sigma)^2/2\right)}\PERIOD
\end{align}
From \eqref{eq:dedu}, we immediately get
\begin{align}\label{eq:Dpos}
\expt{\Delta a_{j,n} \indicator{\{\Delta a_{j,n}>0\}}} &\approx \mu_{j,n}\LP 1-p_{j,n}\RP + \frac{\sigma_{j,n}}{\sqrt{2 \pi}} \exp\LP-\frac{q^2_{j,n}}{2}\RP\COMMA\\\nonumber
\expt{\Delta a_{j,n} \indicator{\{\Delta a_{j,n}<0\}}} &\approx \mu_{j,n}\,p_{j,n} - \frac{\sigma_{j,n}}{\sqrt{2 \pi}} \exp\LP-\frac{q^2_{j,n}}{2}\RP\PERIOD
\end{align}
By subtracting the expressions in \eqref{eq:Dpos}, we obtain 
\begin{align}\label{eq:Evarap}
\expt{\var{(\varphi_{n+1} \cdot e_{n+1})\SEP \mathcal{G}_n}}&\approx \frac{\Delta t_n}{2} \sum_j f^2_{j,n}
\left( \tilde{\mu}_{j,n}+\tilde{\sigma}_{j,n} \right)\PERIOD
\end{align}

Let us now consider the case where $a_i(x_n)\Delta t_n/2$ is close to zero for some $i$. 
We can bound the expression $\expt{\Delta a_{j,n}\, \text{sgn}(\Delta a_{j,n})}$ by  
$\expt{|\Delta a_{j,n}|}$ and  also $\sqrt{\expt{(\Delta a_{j,n})^2}}$.
It is easy to see that $\expt{|\Delta a_{j,n}|}\leq{\bar{\mu}_{j,n}}$. 
Regarding $\expt{(\Delta a_{j,n})^2}$, it can be approximated by
\begin{align*}
\expt{\sum_{i,i'} (\nabla a_j(x_n)\cdot \nu_i)\,(\nabla a_j(x_n)\cdot \nu_{i'})\,\mathcal{Q}_i\,\mathcal{Q}_{i'}} = 
\sum_{i,i'} (\nabla a_j(x_n)\cdot \nu_i)\,(\nabla a_j(x_n)\cdot \nu_{i'})\,\expt{\mathcal{Q}_i\,\mathcal{Q}_{i'}}\PERIOD
\end{align*}
Since 
\begin{align}\label{eq:eeqq}
\expt{\mathcal{Q}_i\,\mathcal{Q}_{i'}} =  \frac{(\Delta t_n)^2}{4}a_{i}(x_n)a_{i'}(x_n)\indicator{i\neq i'} + \LP a_i(x_n)\frac{\Delta t_n}{2} + \LP a_i(x_n)\frac{\Delta t_n}{2}\RP^2\RP\indicator{i= i'}\COMMA
\end{align}
we can rearrange terms and approximate $\expt{(\Delta a_{j,n})^2}$ by $\mu_{j,n}^2+\sigma^2_{j,n}$.

We conclude that $\expt{\Delta a_{j,n}\, \text{sgn}(\Delta a_{j,n})}$ can be bounded by $m_{j,n}$, which has been defined as $\min\{\bar{\mu}_{j,n}, \sqrt{\mu_{j,n}^2+\sigma^2_{j,n}}\}$.

\begin{rem}
Formula  \eqref{eq:varhatestimated} can be considered as an initial, relatively successful attempt to estimate $\mV{\ell}$, but there is still room for improvement. The main problem is the lack of sharp concentration inequalities for linear combinations of independent Poisson random variables. 
With the numerical examples, we show that the efficiency index of the formula is acceptable for our estimation purposes. 
\end{rem}
\begin{rem}
We are assuming that only tau-leap steps are taken, but in our hybrid algorithms, some steps can be exact, and, hence, do not contribute to the local error. 
For that reason, we include the indicator function of the tau-leap step, $\indicator{\TL}$, in the estimator, $\hat{\mathcal{V}}_{\ell}$. 
\end{rem}
\begin{rem}
The dual-weighted residual approach makes the estimation of $\mV{\ell}$  feasible.
{In our numerical experiments, we found that, using the same number of simulated coupled hybrid paths, the variance of  $\hat{\mathcal{V}}_{\ell}$ is much smaller than the variance of $\var{g_{\ell}{-}g_{\ell-1}}$, estimated by a standard Monte Carlo. Note that $\hat{\mathcal{V}}_{\ell}$ can be computed using only single-level hybrid paths at level $\ell{-}1$. In the upper right panel of Figure \ref{fig:gtt-diag}, we can see that due to the hybrid nature of the simulated paths, it is not possible to predict where the variance of $g_{\ell}{-}g_{\ell-1}$ will enter into a superlinear regime. Thus, by extrapolating the $\var{g_{\ell}{-}g_{\ell-1}}$ from the coarser levels, we may overestimate the values of $\var{g_{\ell}{-}g_{\ell-1}}$ for the deepest levels.}
\end{rem}

\section{Estimation Procedure}
\label{EEC}

In this section, we present a procedure that estimates $\expt{g(X(T))}$ within a given prescribed relative tolerance, $TOL{>}0$, with high probability. The process contains three phases:
\begin{description}
\item[Phase I] Calibration of virtual machine-dependent quantities.
\item[Phase II] Solution of the work optimization problem:  we obtain 
the total number of levels, $L$, and the sequences $\seqof{\delta_\ell}{\ell=0}{L}$ and $\seqof{M_\ell}{\ell=0}{L}$, \ie, the one-step exit probability bounds and the required number of simulations at each level. We recall that in Section \ref{MLMCintro}, we defined $\Delta t_\ell:= \Delta t_0 {R^{-\ell}}$, where $R>1$ is a  given  integer constant.
For that reason, to define the whole sequence of meshes, $\seqof{\Delta t_\ell}{\ell=0}{L}$, we simply need to define the size of the coarsest mesh, $\Delta t_0$.
\item[Phase III] Estimation of $\expt{g(X(T))}$. 
\end{description} 

\subsection{Phase I}

In this section, we describe the estimation of several constants, 
$C_1$, $C_2$, $C_3$ and $K_1$, and functions, $C_P$ and $K_2$,
that allow us to model the expected computational work (or just work), measured in terms of the runtime of hybrid paths, see definitions \eqref{eq:cost} and \eqref{eq:costc}. Those quantities are virtual machine dependent; that is, they are dependent on the computer system used for running the simulations and also on  the implementation language. Those quantities are also off-line estimated; that is, we need to estimate them only once for each virtual machine on which we want to run the hybrid method.

Constants $C_1$, $C_2$, and $C_3$ reflect the average execution times of each logical path of Algorithm \ref{alg:sel}. We have that $C_1$ and $C_2$ reflect the work associated with the two different types of steps in the MNRM. Constant $C_3$ reflects the work needed for computing the Chernoff tau-leap size, $\tau_{Ch}$. Finally, when we perform a tau-leap step, we have the work needed for simulating Poisson random variates, which is modeled by the function $C_P$ \cite{ourSL}. This function has two constants that are also virtual machine dependent.


The constant, $K_1$, and the function, $K_2\equiv K_2(x,\delta)$, defined through $C_1$, $C_2$, and $C_3$, were introduced in Section \ref{hybrid_algo}.

\subsection{Phase II}
In this section, we set and solve the work optimization problem.
Our objective function is the expected total work of the MLMC estimator, $\mathcal{M}_L$, defined in \eqref{MLMCest}, \ie,
\begin{equation*}
  \sum_{\ell =0}^{L}\psi_\ell M_{\ell}\COMMA
\end{equation*}
where $L$ is the maximum level (deepest level), $\psi_0$ is the expected work of a single-level path at level 0, and $\psi_\ell$, for $\ell\geq 1$, is the expected computational work of two coupled paths at levels $\ell{-1}$ and $\ell$. Finally, $M_0$ is the number of single-level paths at level 0, and $M_\ell$, for $\ell\geq 1$, is the number of coupled paths at levels $\ell{-1}$ and $\ell$. 

Let us now describe in detail the quantities, $\seqof{\psi_\ell}{\ell=0}{L}$.
For $\ell{=}0$, Algorithm \ref{alg:pathdualscostSL} generates a single hybrid path. The building block of a single hybrid path is Algorithm \ref{alg:sel}, which adaptively determines whether to use an MNRM step or a tau-leap one. According to this algorithm, there are two ways of taking an MNRM step, depending on the logical conditions, $K_1/a_0(x){>}T_0{-}t$ and $K_2/a_0(x){>} \tau_{Ch}$.
Given one particular hybrid path, let $\NMNRMKone(\Delta t_0,\delta_0)$ be the number of MNRM steps such that $K_1/a_0(x){>}T_0{-}t$ is true, and let $\NMNRMKtwo(\Delta t_0,\delta_0)$ be the number of MNRM steps such that $K_1/a_0(x) {>}T_0{-}t$  is false and $K_2/a_0(x){>} \tau_{Ch}$ is true. 
When a Chernoff tau-leap step is taken, we have constant work, $C_3$, and  variable work computed with the aid of $C_P$.
Then, the expected work of a single hybrid path, at level $\ell=0$, is
\begin{align}\label{eq:cost}
\psi_0 &:= C_1 \expt{\NMNRMKone(\Delta t_0,\delta_0)} + C_2 \expt{\NMNRMKtwo(\Delta t_0,\delta_0)} + C_3 \expt{\NTL (\Delta t_0,\delta_0)}\\ &+ \sum_{j=1}^J \expt{\int_{[0,T]} C_P(a_j(\bar X_0(s))\tau_{Ch}(\bar X_0(s),\delta_0))\indicator{TL}(\bar X_0(s)) ds} \nonumber
\COMMA
\end{align} 
where $\Delta t_0$ is the size of the time mesh at level 0 and $\delta_0$ is the exit probability bound at level 0.
Therefore, the expected work at level 0 is $\psi_0 M_0$, 
where $M_0$ is the total number of single hybrid paths.

For $\ell \geq 1$, we use Algorithm \ref{alg:coupled} to generate $M_{\ell}$-coupled paths that couple the $\ell{-}1$ and $\ell$ levels.
Given two coupled paths, let $\NMNRMKone(\Delta t_{\ell{-}1},\delta_{\ell-1})$ and $\NMNRMKone(\Delta t_{\ell},\delta_{\ell})$ be the number of exact steps for level ${\ell{-}1}$ (coarse mesh) and $\ell$ (fine mesh), respectively, with associated work $C_1$. We define $\NMNRMKtwo(\Delta t_{\ell{-}1},\delta_{\ell-1})$ and $\NMNRMKtwo(\Delta t_{\ell},\delta_{\ell})$ analogously. 
Then, the expected work of a pair of coupled hybrid paths at levels $\ell$ and $\ell-1$ is
\begin{align}\label{eq:costc}
\psi_\ell &:= C_1 \expt{\NMNRMKoneC(\ell)} + C_2 \expt{\NMNRMKtwoC(\ell)} + C_3 \expt{\NTLC(\ell)} \\ \nonumber
&+ \sum_{j=1}^J \expt{\int_{[0,T]} C_P(a_j(\bar X_\ell(s))\tau_{Ch}(\bar X_\ell(s),\delta_\ell))\indicator{TL}(\bar X_\ell(s)) ds}\\ \nonumber
 &+ \sum_{j=1}^J \expt{\int_{[0,T]} C_P(a_j(\bar X_{\ell-1}(s))\tau_{Ch}(\bar X_{\ell-1}(s),\delta_{\ell-1}))\indicator{TL}(\bar X_{\ell-1}(s)) ds}\COMMA
\end{align} 
where 
\begin{align*}
\NMNRMKoneC(\ell)&:=\NMNRMKone(\Delta t_\ell,\delta_\ell) + \NMNRMKone(\Delta t_{\ell-1},\delta_{\ell-1})\\
\NMNRMKtwoC(\ell)&:=\NMNRMKtwo(\Delta t_\ell,\delta_\ell)+\NMNRMKtwo(\Delta t_{\ell-1},\delta_{\ell-1})\\
\NTLC(\ell) &:=\NTL(\Delta t_\ell,\delta_\ell)+\NTL (\Delta t_{\ell-1},\delta_{\ell-1})\PERIOD
\end{align*} 



Now, recalling the definitions of the error decomposition, given at the beginning of Section \ref{calibration}, we have all the elements to formulate the work optimization problem.
Given a relative tolerance, $TOL{>}0$, we solve
\begin{align}
\left\{ \!
\begin{array}{l}\label{eq:minworkprob}
\min_{\{\Delta t_0, L,\seqof{M_{\ell},\delta_{\ell}}{\ell=0}{L}\}} \sum_{\ell =0}^{L}\psi_\ell M_{\ell}\\
\mbox{s.t.}\\
\EE_{E,L} +\EE_{I,L}+\EE_{S,L} \leq TOL\PERIOD
\end{array}
\right.
\end{align}

It is natural to consider the following family of auxiliary problems indexed on $L{\geq} 1$, where we assume for now that the double sequence, $(\Delta t_{\ell}, \delta_{\ell})_{\ell=0}^L$, is known: 
\begin{align}
\left\{ \!
\begin{array}{l}\label{eq:auxminprob}
\min_{\seqof{M_{\ell} \geq 1}{{\ell}=0}{L}} \sum_{{\ell}=0}^{L}\psi_\ell M_{\ell}\\
\mbox{s.t.}\\
\EE_{I,L}+ C_A \sqrt{\sum_{\ell=0}^L \frac{\mV{\ell}}{M_{\ell}}} \leq TOL{-}TOL^2\COMMA
\end{array}
\right.
\end{align}
where we have $C_A \geq 2$ to guarantee an asymptotic confidence level of at least 95\%.

Let us assume for now that  we know $\psi_\ell$, $\mV{\ell}$ and $\EE_{I,\ell}$, for $\ell=0,1,\ldots,L$. Let $L_0$ be the smallest value of $L$ such that $\EE_{I,L}{< }TOL{-}TOL^2$. 
This value exists and it is finite since the discretization error, $\EE_{I,L}$, tends to zero as 
$L$ goes to infinity.
For each $L\geq L_0$, define $w_{L}{:=}\sum_{{\ell}=0}^{L}\psi_\ell M^*_{\ell}$, where the sequence $\seqof{M^*_{\ell}}{\ell=0}{L}$ is the solution of the problem \eqref{eq:auxminprob}.
It is worth mentioning that $\seqof{M^*_{\ell}}{\ell=0}{L}$ is quickly obtained as the solution of  
the following Karush-Kuhn-Tucker problem (see, \eg, \cite{Luenberger}):
\begin{align}\label{eq:greedyKKT} 
\left\{ \!
\begin{array}{l}
\min_{\seqof{M_{\ell} \geq 1}{{\ell}=0}{L}}  \sum_{\ell=0}^L \psi_{\ell} M_{\ell}\\
\mbox{s.t.}\\
\sum_{\ell=0}^L \frac{\mV{\ell}}{M_{\ell}}\leq R
\end{array}
\right.\PERIOD
\end{align}
We do not develop here all the calculations, but a pseudo code is given in Algorithm \ref{alg:greedyoptKKT}.

Let us now analyze two extreme cases: i) for $L$ such that $\EE_{I,L}$ is less but very close to $TOL{-}TOL^2$, we have that  ${\sum_{\ell=0}^L {\mV{\ell}}/{M^*_{\ell}}}$ is a very small number. 
As a consequence, we obtain large values of $M^*_{\ell}$ and, hence, a large value of $w_{L}$. 
By adding one more level, \ie,  $L{\leftarrow}L{+}1$,  we expect a larger gap between $\EE_{I,L}$ and $TOL_0$; that means that we expect a larger value of ${\sum_{\ell=0}^L {\mV{\ell}}/{M^*_{\ell}}}$ that may lead to smaller values of $M^*_{\ell}$.  We observe that, in spite of adding one more term to $w_L$, this leads to a smaller value of $w_L$.
ii) At the other extreme, a large value of $L$ is associated with large values of $\psi_L$ and therefore with large values of $w_L$. 

This informal `extreme case analysis' has been confirmed by our numerical experiments (see, for instance, Figures \ref{fig:dec2-diag} and \ref{fig:gtt-diag}  (lower-right)), which allow us to conjecture that the sequence 
$\seqof{w_L}{L=L_0}{+\infty}$ is a convex function of $L$ and, hence, that it has a unique optimal value achieved at a certain $L^*$. 
A pseudo algorithm to find $L^*$ could be to start computing $w_{L_0}$ and  $w_{L_0+1}$. 
If  $w_{L_0+1}{\geq}w_{L_0}$, we accept $L^*{=}L_0$;
otherwise, we proceed to computing the next term of the sequence, $\seqof{w_L}{L=L_0}{+\infty}$. 
If, for some $p$, we have $w_{L_{p+1}}{\geq}w_{L_p}$, we accept $L^*{=}L_p$.
Of course, we can stop even if $w_{L_{p+1}}{<}w_{L_p}$, but the difference $\abs{w_{L_{p+1}}{-}w_{L_p}}$ is sufficiently small.
In this last case, we accept  $L^*{=}L_{p+1}$.

\subsubsection{Computational Complexity}\label{sec:compucomple}
\begin{figure}[h!]
\centering
\begin{minipage}{0.49\textwidth}
\includegraphics[scale=0.31]{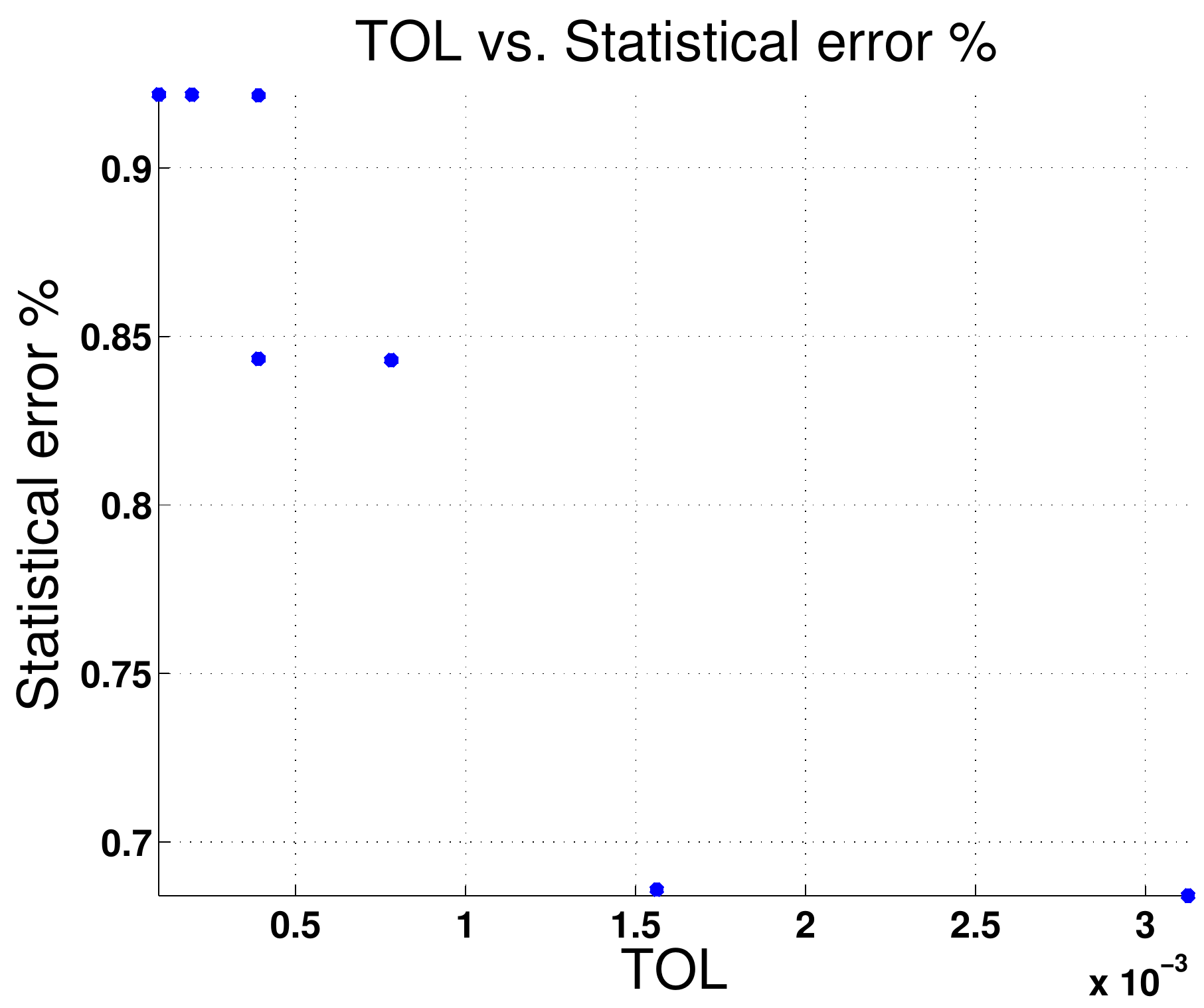}
\end{minipage}
\begin{minipage}{0.49\textwidth}
\includegraphics[scale=0.31]{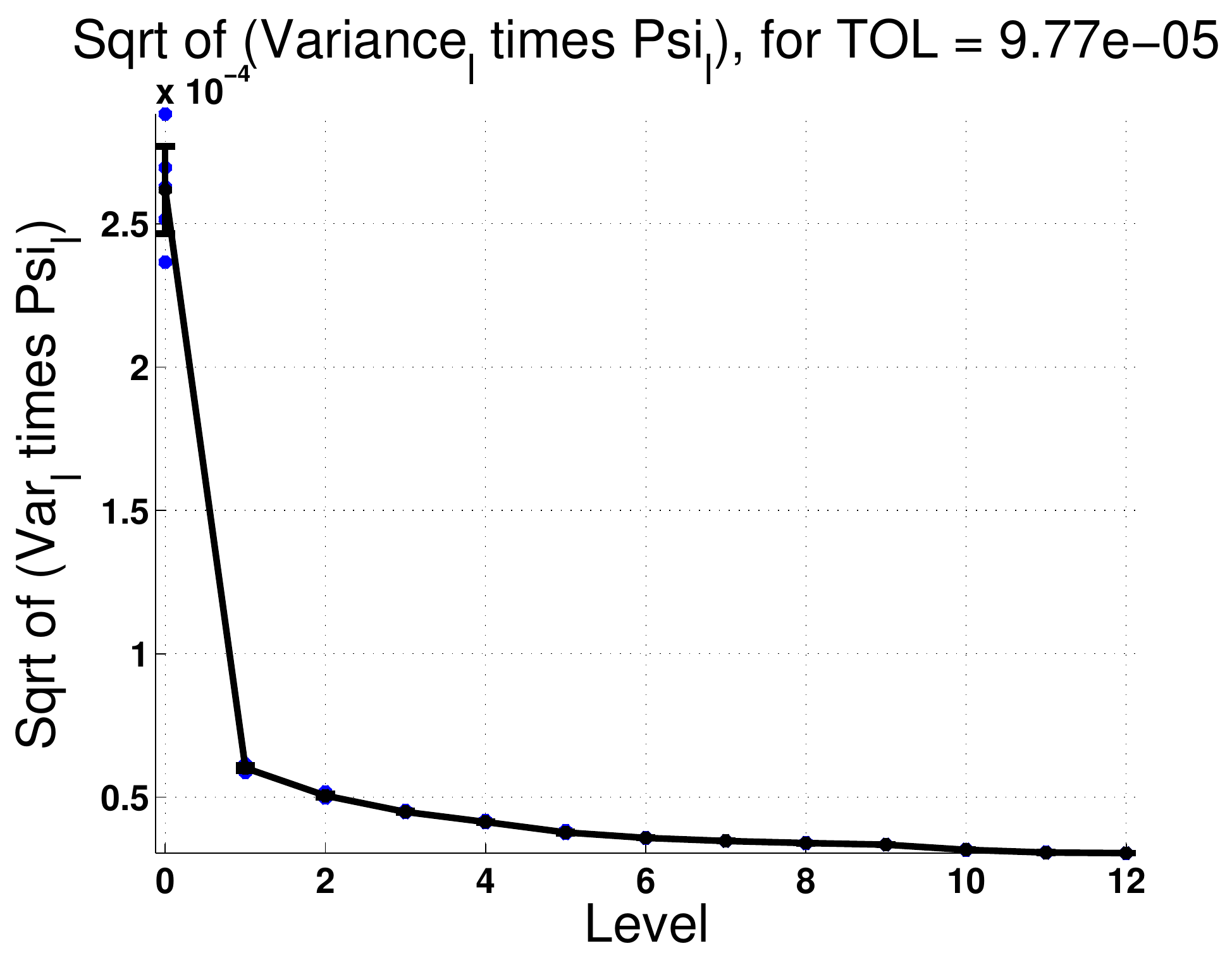}
\end{minipage}
\caption{Left: Percentage of the statistical error over the computational global error, for the simple decay model \eqref{ex:dec}. As mentioned in Section \ref{EEC}, it is well above $0.5$ for all the tolerances. Right: $\sqrt{\hmV{\ell}\hat{\psi}_\ell}$ as a function of $\ell$, for the smallest tolerance, which  decreases as the level increases. Observe that the contribution of level 0 is less than 50\% of the sum of the other levels.
}
\label{fig:statdec2}
\end{figure}

\begin{figure}[h!]
\centering
\begin{minipage}{0.49\textwidth}
\includegraphics[scale=0.31]{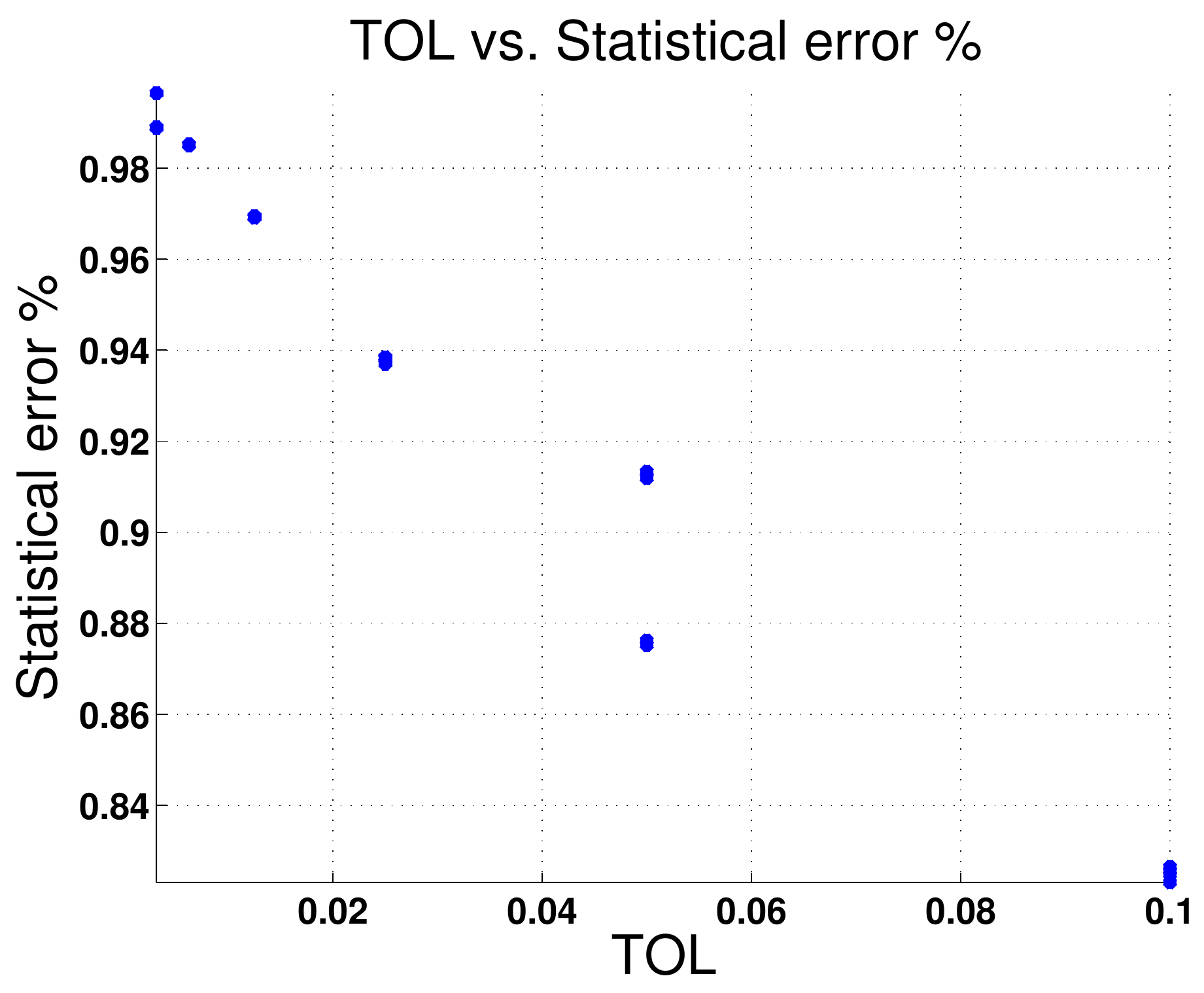}
\end{minipage}
\begin{minipage}{0.49\textwidth}
\includegraphics[scale=0.31]{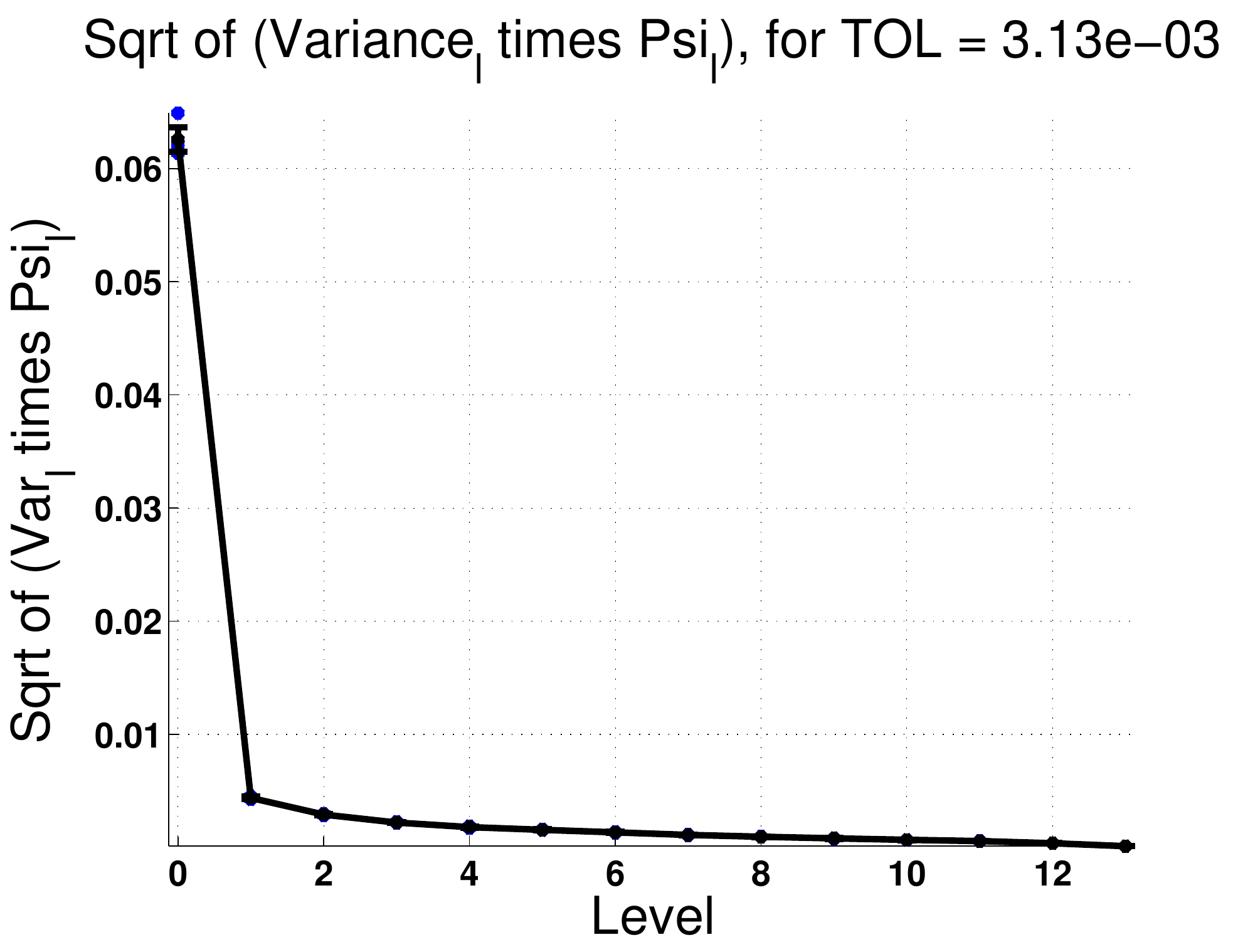}
\end{minipage}
\caption{Left: Percentage of the statistical error over the computational global error, for the gene transcription and translation model \eqref{ex:gtt}. As  mentioned in Section \ref{EEC}, it is well above $0.5$ for all the tolerances. Right: $\sqrt{\hmV{\ell}\hat{\psi}_\ell}$ as a function of $\ell$, for the smallest tolerance, which  decreases as the level increases. Observe that the contribution of level 0 is almost equal to the sum of the other levels.
}
\label{fig:statgtt}
\end{figure}

{
At this point, we have all the necessary elements to establish a key point of this work, the computational complexity of the multilevel hybrid Chernoff tau-leap method.

Let us now analyze the optimal amount of work at level $L$, $w_{L}$, as a function of the given relative tolerance, $TOL$. For simplicity, let us assume that $M^*_{\ell}{>}1,\, \ell{=}0,...,L$. In this case,  the optimal number of samples at level $\ell$ is given by
\begin{equation*}
M^*_{\ell} {=}  (C_A/\theta)^2  TOL^{-2}\sqrt{\mV{\ell}/ \psi_{\ell}} \sum_{\ell=0}^L \sqrt{\mV{\ell} \psi_{\ell}}\COMMA
\end{equation*} 
for some $\theta \in (0,1)$. In fact, $\theta$ is the proportion of the tolerance, $TOL$, that our computational cost optimization algorithm selects for the statistical error, $\EE_{S,L}$. In our algorithms, we impose 
$\theta\geq 0.5$; however, our numerical experiments always select a larger value (see Figures \ref{fig:statdec2} and \ref{fig:statgtt}).

By substituting $M^*_{\ell}$ into the total work formula, $w_{L}$, we conclude that the optimal expected work, conditional on $\theta$, is given by 
\begin{equation*}
\expt{w^*_{L} (TOL)\SEP \theta} = \LP\frac{C_A}{\theta}\sum_{\ell=0}^{L(\theta)} \sqrt{\mV{\ell} \psi_{\ell}}\RP^{2} TOL^{-2}\PERIOD
\end{equation*}
Due to the constraint, $\theta\geq 0.5$, we have that 
\begin{equation*}
w^*_{L} (TOL)  \leq \sup_L\left\{ \LP{2\,C_A}\sum_{\ell=0}^L \sqrt{\mV{\ell} \psi_{\ell}}\RP^{2}\right\} TOL^{-2}\PERIOD
\end{equation*}
Let us consider the series $\sum_{\ell=0}^{\infty} \sqrt{\mV{\ell} \psi_{\ell}}$. 
First, observe that the expected computational work per path at level $\ell$, $\psi_{\ell}$, is bounded by a multiple of the expected computational work of the MNRM (see Section \ref{sec:MNR}), \ie, $K \psi_{\text{MNRM}}$. In our numerical experiments, we observe that taking $K$ around $3$ is enough. 
Therefore, $\sum_{\ell=0}^{\infty} \sqrt{\mV{\ell} \psi_{\ell}} \leq \sqrt{K \psi_{\text{MNRM}}} \sum_{\ell=0}^{\infty}  \sqrt{\mV{\ell}}$.
Observe that, by construction, $\mV{\ell}\rightarrow 0$, superlinearly. More specifically, it satisfies the bound  $\mV{\ell}=\Ordo{\Delta t_\ell} \leq C \Delta t_0(1/2)^{\ell}$ for some positive constant $C$. 
Therefore, the series $\sum_{\ell =0}^{\infty} \sqrt{\mV{\ell}}$ is dominated by the geometric series $\sum_{\ell =0}^{\infty} (1/\sqrt{2})^{\ell}<\infty$ . We conclude that 
$\sup_{L}\{\sum_{\ell=0}^L \sqrt{\mV{\ell} \psi_{\ell}}\}$ is bounded and, therefore, the expected computational complexity of the multilevel hybrid Chernoff tau-leap method is $w^*_L(TOL){=}\Ordo{TOL^{-2}}$. 
}

\subsubsection{Some Comments on the Algorithms for Phase II}

In Algorithm \ref{alg:Cal}, we propose an iterative method to obtain an approximate solution to the problem
\eqref{eq:minworkprob}. Notice that we are assuming that there are at least two levels in the multilevel hierarchy, \ie, $L\geq 1$. 


To solve the problem \eqref{eq:minworkprob}, we bound the global exit error, 
$\EE_{E,L}$, by $TOL^2$. More specifically, we choose $\delta_L$ to be sufficiently small such that 
\begin{equation}\label{eq:lastdelta}
|\avg{g_{L}}{\cdot}|\,\, \delta_{L}\,\,  \avg{\NTL(\Delta t_L,\delta_L)}{\cdot} < TOL^2\PERIOD
\end{equation} 
At this point, it is crucial to observe that {if we impose the condition \eqref{eq:lastdelta} on any level $\ell{<}L$, then we are unnecessarily enforcing a dependence of $\delta_\ell$ on $TOL$. This dependence may result in very small values of $\delta_\ell$, which in turn may increase the expected number of exact steps and tau-leap steps at level $\ell$, implying a larger expected computational work at level $\ell$}. In the appendix of \cite{ourSL}, we proved that, when $\delta_\ell$ tends to zero, the expected values of the number of tau-leap steps at level $\ell$ go to zero, and therefore our hybrid MLMC strategy would converge to the SSA method without the desired reduction in  computational work. 
To avoid the dependence of $\seqof{\delta_\ell}{\ell=0}{L-1}$ on $TOL$, we adopt a different strategy based on the following decomposition:
\begin{align*}
\mV{\ell} = \var{g_\ell \indicator {A_{\ell}} - g_{\ell-1} \indicator{A_{\ell-1}}} &= \var{g_\ell - g_{\ell-1} \SEP 
A_{\ell}\cap A_{\ell-1}}\prob{A_{\ell}\cap A_{\ell-1}}\\
&+ \var{g_\ell\SEP A_{\ell}\cap A^c_{\ell-1}}\prob{A_{\ell}\cap A^c_{\ell-1}}\\
&+ \var{g_{\ell-1}\SEP A^c_{\ell}\cap A_{\ell-1}}\prob{A^c_{\ell}\cap A_{\ell-1}}\PERIOD
\end{align*}   
We impose that the first term of the right-hand side dominates the other two. This is because the conditional variances appearing in the last two terms are of order $\Ordo{1}$, while the conditional variance appearing in the first term is of order $\Ordo{\Delta t_{\ell}}$, and we make our computations with approximations of $\mV{\ell}$ assuming that $\prob{A_{\ell}\cap A_{\ell-1}}$ is close to one.  We proceed as follows: first, we approximate $\prob{A_{\ell}\cap A_{\ell-1}}$ by $\prob{A_{\ell}} \prob{A_{\ell-1}}$; then, we consider $1{-}\delta_\ell \avg{\NTL(\Delta t_\ell,\delta_\ell)}{\cdot}$ as an approximate upper bound for $\prob{A_{\ell}}$ when  $\delta_\ell \avg{\NTL(\Delta t_\ell,\delta_\ell)}{\cdot}{\ll} 1$. Those considerations lead us to impose 
\begin{align}\label{eq:deltas}
&\var{g_\ell - g_{\ell-1} \SEP A_{\ell}\cap A_{\ell-1}}
\LP 1{-}\delta_\ell \avg{\NTL(\Delta t_\ell,\delta_\ell)}{\cdot}\RP
\LP 1{-}\delta_{\ell-1} \avg{\NTL(\Delta t_{\ell-1},\delta_{\ell-1})}{\cdot}\RP>\\
&\var{g_\ell\SEP {A_{\ell}}\cap{A^c_{\ell-1}}}\delta_{\ell-1} \avg{\NTL(\Delta t_{\ell-1},\delta_{\ell-1})}{\cdot} + 
\var{g_{\ell-1}\SEP{A^c_{\ell}}\cap{A_{\ell-1}}}\delta_{\ell} \avg{\NTL(\Delta t_{\ell},\delta_{\ell})}{\cdot}\PERIOD\nonumber
\end{align}
To avoid simultaneous refinements on $\delta_\ell$ and $\delta_{\ell-1}$, based on \eqref{eq:deltas}, we impose on $\delta_\ell$ the following condition:
\begin{equation}\label{eq:deltasapp}
\hmV{\ell}\LP 1{-}\delta_{\ell} \avg{\NTL(\Delta t_{\ell},\delta_{\ell})}{\cdot}\RP^2 >
 2\,\, \svar{g}{\cdot}\delta_{\ell} \avg{\NTL(\Delta t_{\ell},\delta_{\ell})}{\cdot} \PERIOD
\end{equation}
 
Algorithms \ref{alg:pathdualscostSL} and \ref{alg:Cal} provide  $\avg{g_{\ell}}{\cdot}$, $\avg{\NTL}{\cdot}$ and the other required quantities. {Condition \eqref{eq:deltasapp} does not affect the telescoping sum property of our multilevel estimator, $\mathcal{M}_L$, defined in \eqref{MLMCest}, since each level, $\ell$, has its own $\delta_{\ell}$.}



\begin{rem}[Multilevel estimators used in Algorithm \ref{alg:Cal}]
\label{rem:algCal}
Although in algorithm \ref{alg:Cal} we show that the estimations of $\expt{g(X(T))}$ and $\var{g(X(T))}$ are computed using the information from the last level only, in fact we are computing them using a multilevel estimator. We omit the details in the algorithm for the sake of simplicity. For the case of $\expt{g(\bar X(T))}$, we use the standard mutilevel estimator, and, for the case of $\var{g(X(T))}$, we use the following telescopic decomposition:
\begin{align*}
\var{g(\bar X_l(T))} = \var{g(\bar X_0(T))} + \sum_{\ell=1}^{l}(\var{g(\bar X_\ell(T))}-\var{g(\bar X_{\ell-1}(T))})\COMMA
\end{align*}
where $l>1$ is a fixed level. Using the usual variance estimators for each level, we obtain an unbiased multilevel estimator of the variance of $g(\bar X)$. We refer to \cite{Bierig.Chernov} for details.
\end{rem}
\begin{rem}[Coupled paths exiting the lattice, $\latt$]
Algorithm \ref{alg:coupled} could compute four types of paths. It could happen that no approximate process (the coarse one, $X_{\ell-1}$, or the fine one, $X_\ell$) exits the lattice, which is the most common case. It could also happen that one of the approximate processes exits the lattice. And finally, both approximate processes could exit the lattice. The first case is the most common one and no further explanation is required. We now explain the case when one of processes exits the lattice. Suppose that the coarse one exits the lattice. In that case, until the fine process reaches time $T$ or exits the lattice, we still simulate the coupled process by simulating only the fine path using the single-level hybrid algorithm presented in \cite{ourSL}. If the fine path reaches $T$, we have that $\indicator{A_{\ell-1}}=0$, and $\indicator{A_{\ell}}=1$. Vice versa, if the fine process exits and the coarse one reaches $T$, we have  $\indicator{A_{\ell-1}}=1$ and $\indicator{A_{\ell}}=0$.
\end{rem}

\begin{rem}[Coupling with an exact path]
\label{rem:unbiased}
Algorithm \ref{alg:Cal} uses a computational-cost-based stopping criterion. That is, the algorithm stops refining the time mesh when the estimated total computational cost of the multilevel estimator, $\hat W_{\ML}{:=}\sum_{\ell =0}^{l}\hat \psi_\ell M_{\ell}$, at level $l$, is greater than the corresponding computational cost for level $l{-}1$, and only when condition $\hat \WE {<} TOL{-}TOL^{2}$ is already satisfied. In that case, $L^*{=}l{-}1$. The latter condition is required for obtaining a solution of the optimization problem \eqref{eq:greedyKKT}. In our numerical experiments, we observed that the computational cost of two hybrid  coupled paths, $\psi_\ell$, may be greater than the computational cost of ``hybrid-exact'' coupled paths; that is, the computational cost of a hybrid path at level $l{-}1$ coupled with an exact path at level $l$. That kind of path, used only at the last level, leads to the following unbiased multilevel estimator:
\begin{align*}
\mathcal{\tilde M}_L &:= \frac{1}{M_0} \sum_{m=1}^{M_0}g_0\indicator{A_0}(\omega_{m,0}) + \sum_{\ell=1}^{L-1} \frac{1}{M_{\ell}} \sum_{m=1}^{M_{\ell}} [g_{\ell} \indicator{A_\ell} - g_{\ell-1} \indicator{A_{\ell-1}}](\omega_{m,\ell}) \\
&+ \frac{1}{M_{L}} \sum_{m=1}^{M_L} [g(X(T)) - g_{\ell-1} \indicator{A_{L-1}}](\omega_{m,L})  \PERIOD
\end{align*}
Therefore, it is possible to add another stopping criterion to Algorithm \ref{alg:Cal} related to the comparison between the estimated computational cost of two hybrid coupled paths and the computational cost of hybrid-exact coupled paths. Please note that the condition $\delta_{L} \avg{N_{\TL,L}}{\cdot}\avg{g_{L}}{\cdot} {\leq} TOL^2$ trivially holds because $\avg{N_{\TL,L}}{\cdot}$ is zero in such a case. {In our numerical examples, there are no significant computational gains in the estimation phase from using that stopping rule and its corresponding estimator. This alternative hybrid unbiased estimator is inspired by the work of Anderson and Higham \cite{Anderson2012}}.
\end{rem}

\subsection{Phase III}
From Phase II, we found that, to compute our multilevel Monte Carlo estimator, $\mathcal{M}_L$, for a given tolerance, we have to run $M_{0}^*$ single hybrid paths with parameters $(\Delta t_0,\delta_0)$ and  $M_{\ell}^*$ coupled hybrid paths with parameters $(\Delta t_{\ell-1},\delta_{\ell-1})$ and $(\Delta t_\ell, \delta_\ell)$, for $\ell=1,2,\ldots,L^*$. But, we will follow a slightly different strategy: we run half of the required simulations and use them to update our estimations of  the sequences  
$\seqof{\EE_{I,\ell}}{\ell=0}{L^*}$, $\seqof{\mV{\ell}}{\ell=0}{L^*}$, and $\seqof{\psi_\ell}{\ell=0}{L^*}$. Then, we solve the problem \eqref{eq:auxminprob} again and re-calculate the values of $M_\ell^*$ for all $\ell$. We proceed iteratively until convergence. 
In this way, we take advantage of the information generated by new simulated paths and update the estimations of the sequences of weak errors, computational costs, and variances, obtaining more control over the total work of the method.



\section{Numerical Examples}
\label{sec:examples}

In this section, we present two examples to illustrate the performance of our proposed method, and we compare the results with the single-level approach given in \cite{ourSL}. For bench-marking purposes, we use Gillespie's Stochastic Simulation Algorithm (SSA) instead of the Modified Next Reaction Method (MNRM), because the former is widely used in the literature.

\subsection{A Simple Decay Model}

The classical radioactive decay model provides a simple and important example for the application of our  method. This model has only one species and one first-order reaction,
\begin{align}
\label{ex:dec}
X \xrightarrow{c} \emptyset\PERIOD
\end{align}
Its stoichiometric matrix, $\nu \in \rset$, and the propensity function, $a:\zset_+ \rightarrow \rset$, are given by 
\begin{align*}
\nu = -1 \mbox{   and   }  a(X) = c\,X \PERIOD
\end{align*}
Here, we choose $c = 1$, and define $g(x)=x$ as the scalar observable. In this particularly simple example, we have that $\expt{g(X(T))|X(t) = X_0} = X_0 \exp(-c(T{-}t))$.
Consider the initial condition $X_0{=}10^5$ and the final time $T{=}0.5$. In this case, the process starts relatively far from the boundary, \ie, it is a tau-leap dominated setting. 

\begin{figure}[h!]
\centering
\begin{minipage}{0.49\textwidth}
\includegraphics[scale=0.31]{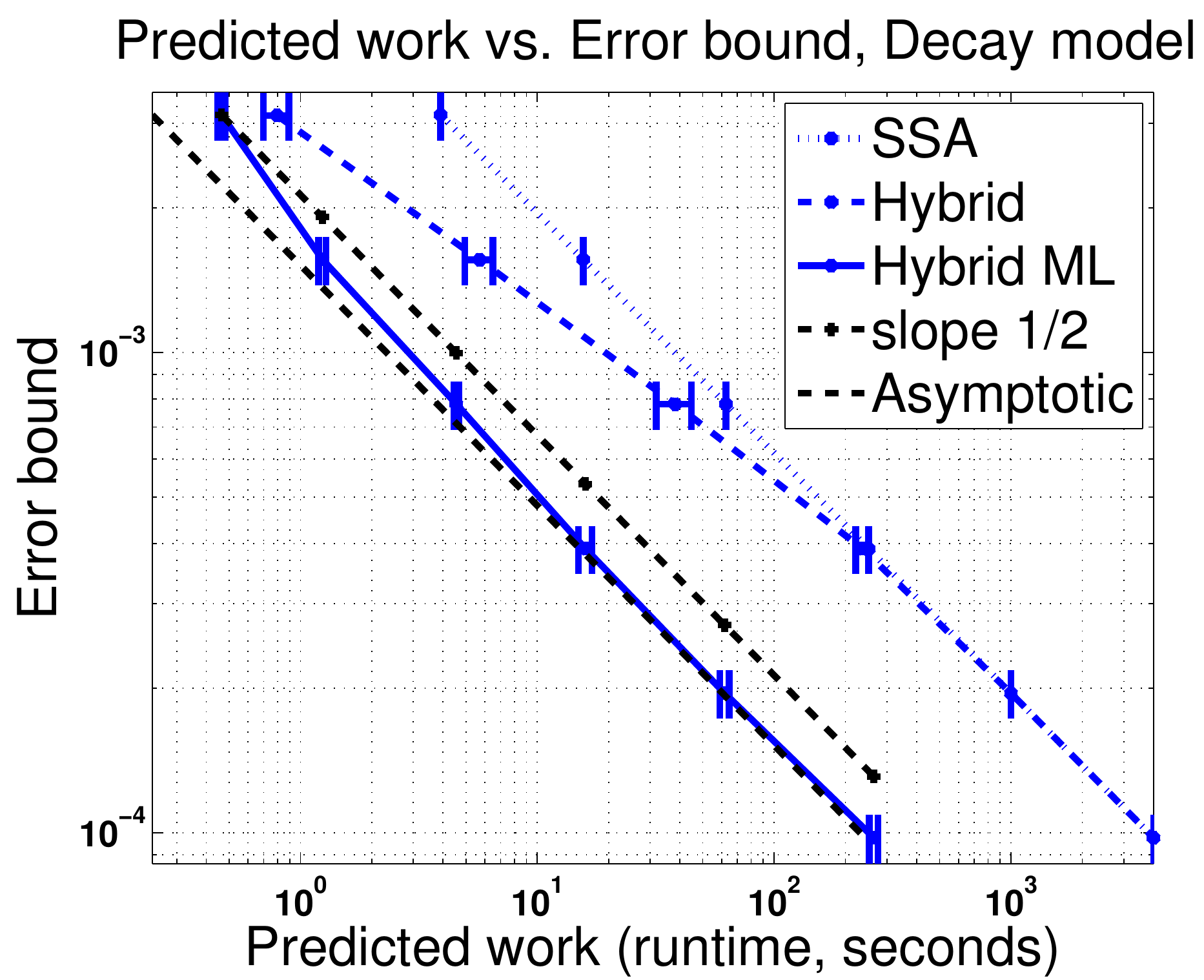}
\end{minipage}
\begin{minipage}{0.49\textwidth}
\includegraphics[scale=0.31]{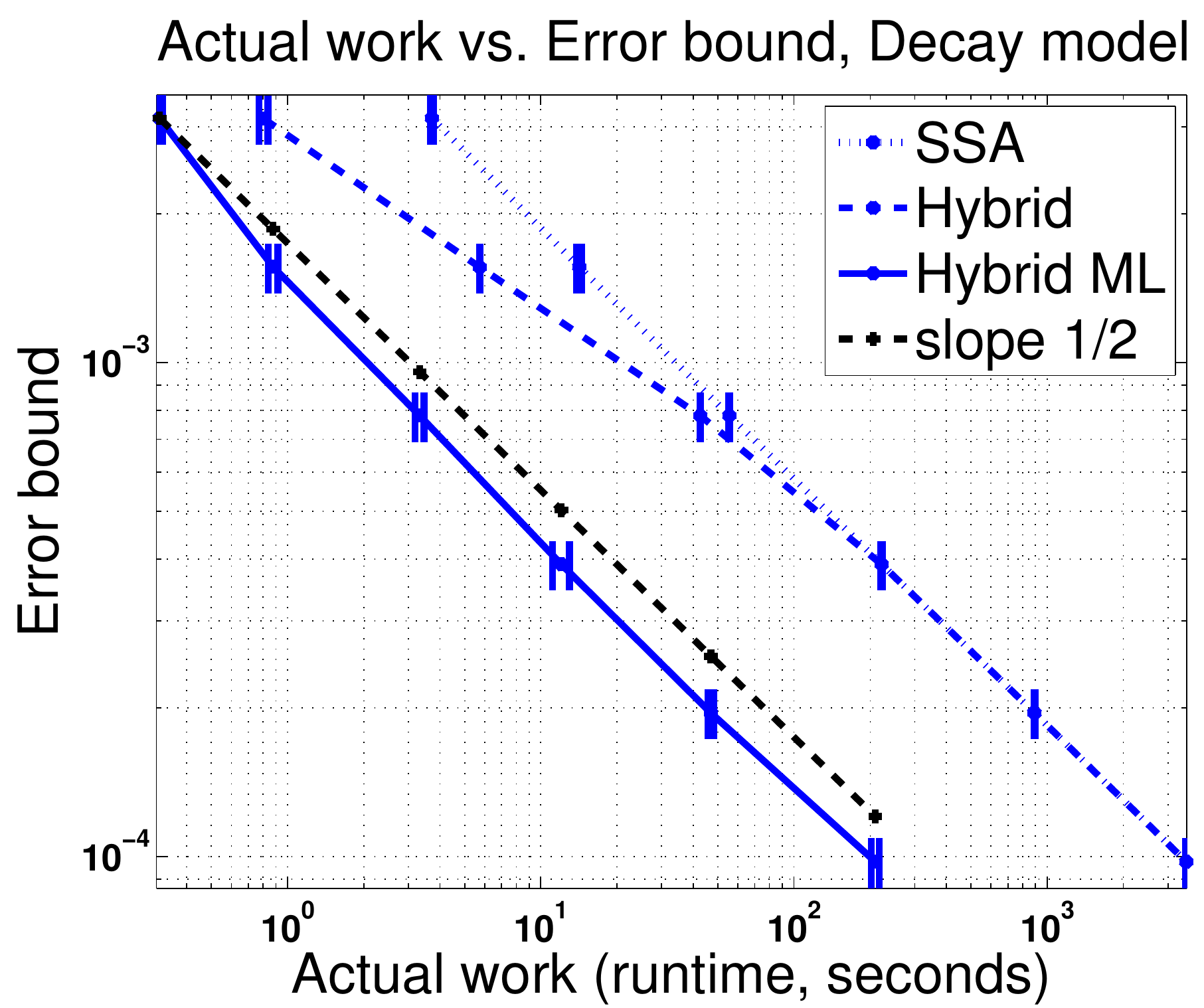}
\end{minipage}
\caption{Left: Predicted work (runtime) versus the estimated error bound for the simple decay model \eqref{ex:dec}, with $95\%$ confidence intervals. The multilevel hybrid method is preferred over the SSA and the single-level method for all the tolerances. Right: Actual computational work (runtime) versus the estimated error bound. Notice that the computational complexity has order $\Ordo{TOL^{-2}}$.}
\label{fig:dec2-worktime}
\end{figure}

We now analyze an ensemble of five independent runs of the calibration algorithm (Algorithm \ref{alg:Cal}), using different relative tolerances. In Figure \ref{fig:dec2-worktime}, we show, in the left panel, the total predicted work (runtime) for the single-level hybrid method, for the multilevel hybrid method and for the SSA method, versus the estimated error bound. The multilevel method is preferred over the SSA and the single-level hybrid method for all the tolerances. We also show the estimated asymptotic work of the multilevel method.
In the right panel, we show, for different tolerances, the actual work (runtime), using a 20 core Intel GLNXA64 architecture and MATLAB version R2014a. 

In Table \ref{tab:dec2}, we summarize an ensemble run of the calibration algorithm, where
$W_{\ML}$ is the average actual computational work of the multilevel estimator (the sum of all the seconds taken to compute the estimation) and $W_{\ssa}$ is the corresponding average actual work of the SSA. We compare those values with the corresponding estimations, $\hat W_{\ML}$ and $\hat W_{\ssa}$. 

\begin{figure}[h!]
\centering
\begin{minipage}{0.49\textwidth}
\includegraphics[scale=0.31]{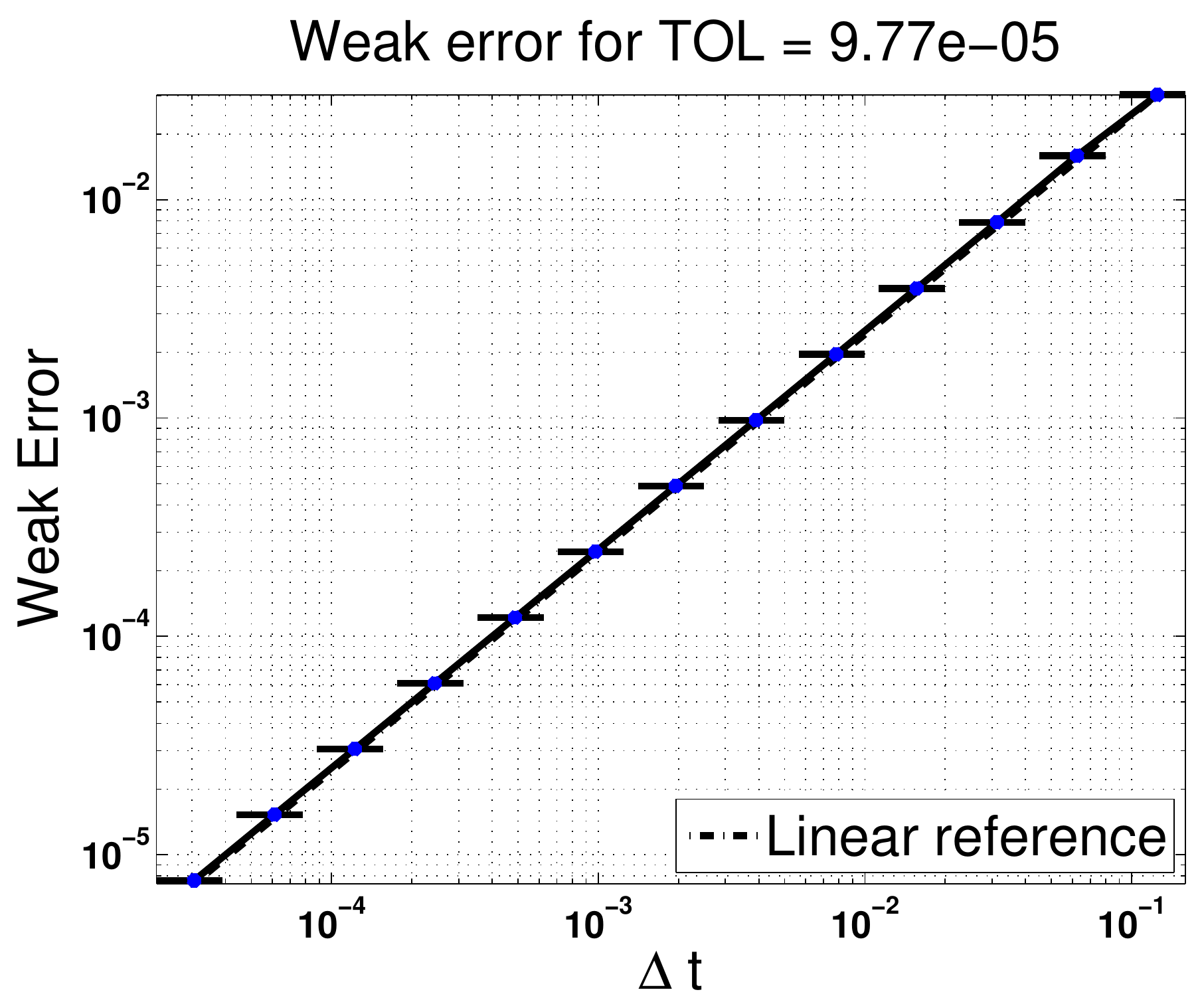}
\end{minipage}
\begin{minipage}{0.49\textwidth}
\includegraphics[scale=0.31]{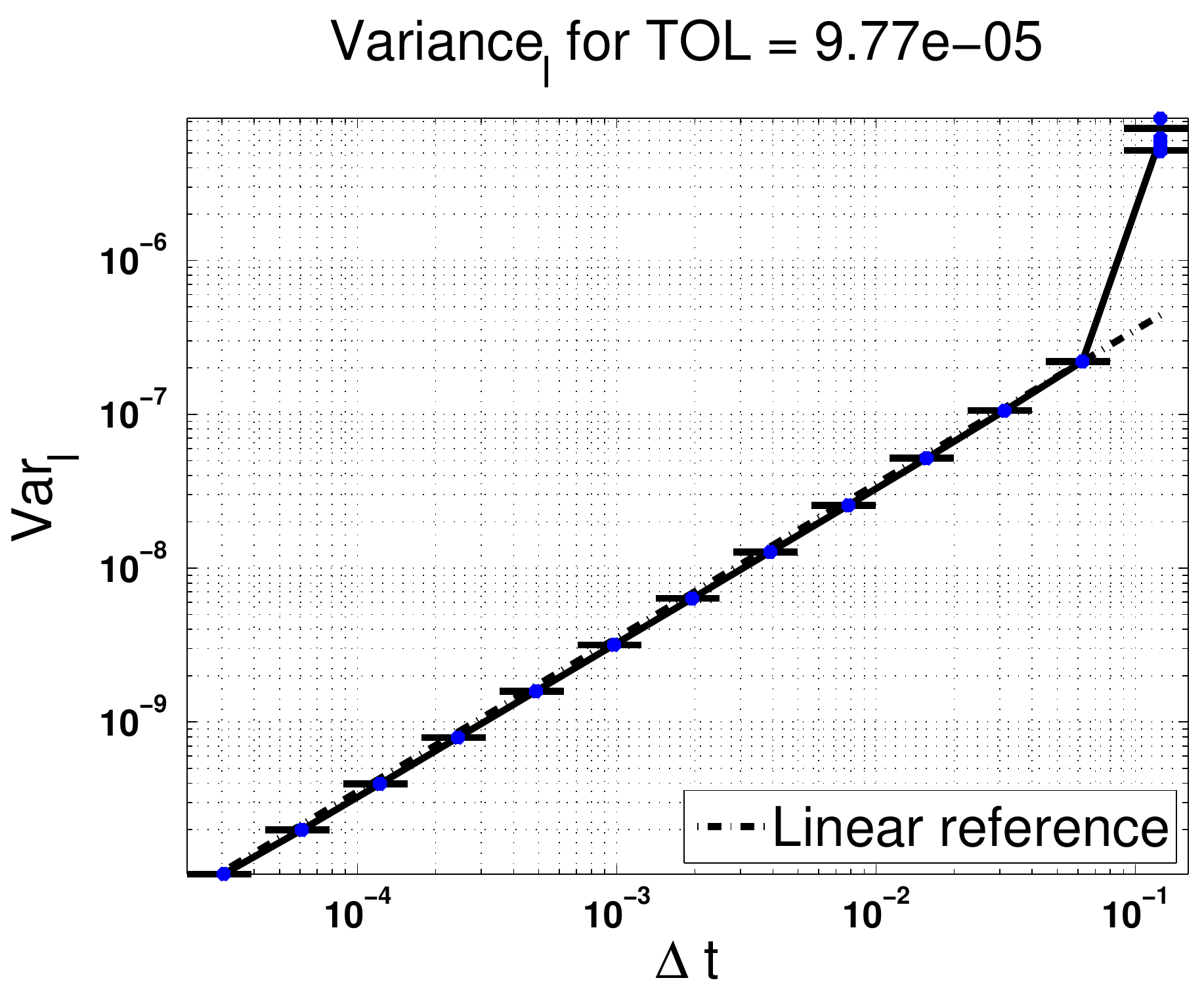}
\end{minipage}
\begin{minipage}{0.49\textwidth}
\includegraphics[scale=0.31]{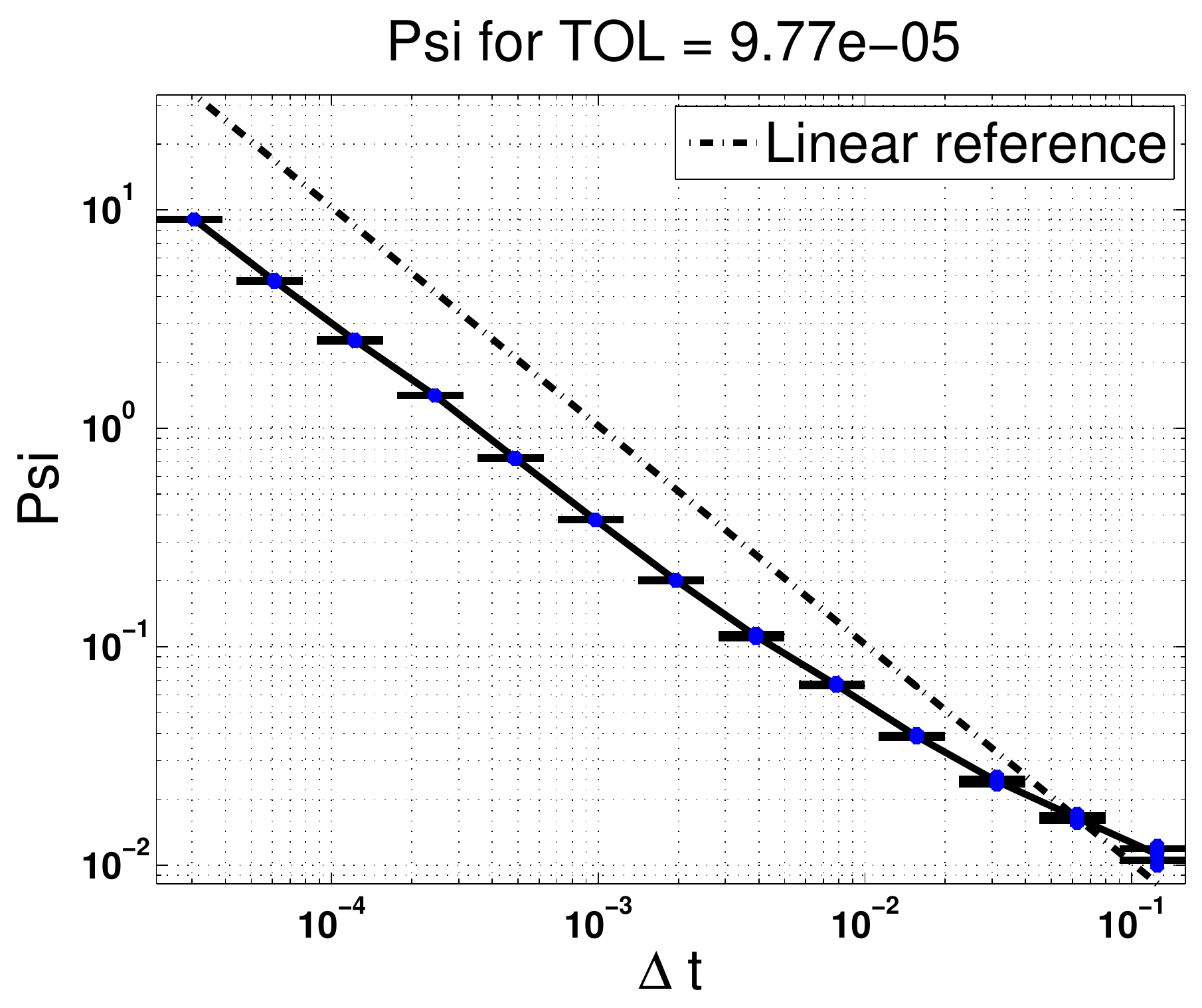}
\end{minipage}
\begin{minipage}{0.49\textwidth}
\includegraphics[scale=0.31]{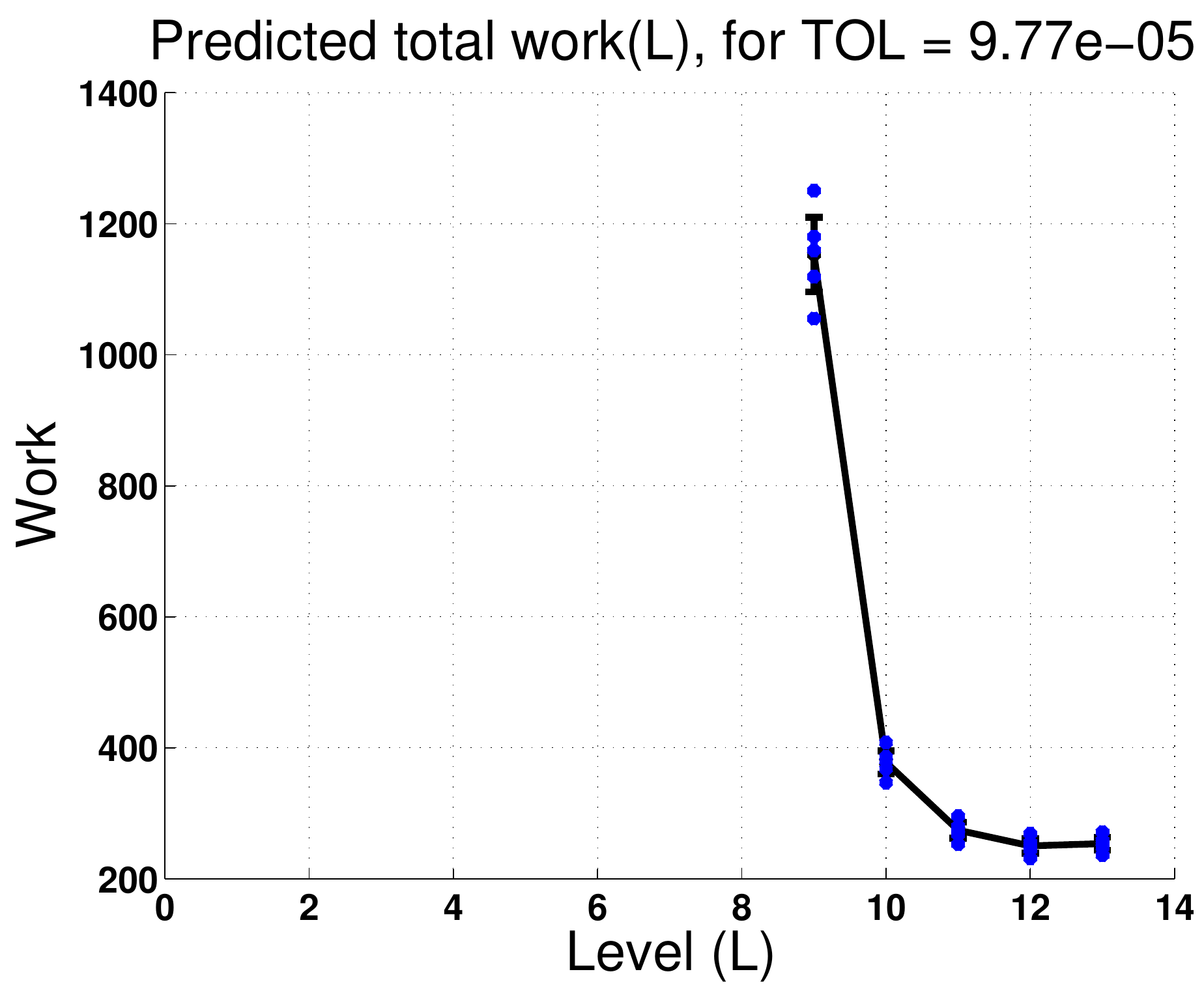}
\end{minipage}
\caption{Upper left: estimated weak error, $\WEH{\ell}$, as a function of the time mesh size, $\Delta t$, for the simple decay model \eqref{ex:dec}. Upper right: estimated variance of the difference between two consecutive levels, $\hmV{\ell}$, as a function of $\Delta t$. Lower left: estimated path work, $\hat \psi_\ell$, as a function of $\Delta t$. Lower right: estimated total computational work, $\sum_{l=0}^L \hat \psi_l M_l$, as a function of the level, $L$.
}
\label{fig:dec2-diag}
\end{figure}

In Figure \ref{fig:dec2-diag}, we can observe how the estimated weak error, $\WEH{\ell}$, and the estimated variance of the difference of the functional between two consecutive levels, $\hmV{\ell}$, decrease linearly as we refine the time mesh. This corresponds to the pure tau-leap case since the process, $X$, remains far from the boundary in 
$[0,T]$.  As expected, the linear relationship for the variance starts at level $1$. 
The estimated total path work, $\hat \psi_{\ell}$, increases as we refine the mesh. Observe that it increases more slowly than linearly. This is because the work needed for generating Poisson random variables becomes less as we refine the time mesh.
In the lower right panel, we show the total computational work, only in the cases in which $\WEH{\ell} < TOL{-}TOL^2$. 
\newcommand{\TOLDEC}{9.77e-05}

\begin{figure}[h!]
\centering
\begin{minipage}{0.49\textwidth}
\includegraphics[scale=0.31]{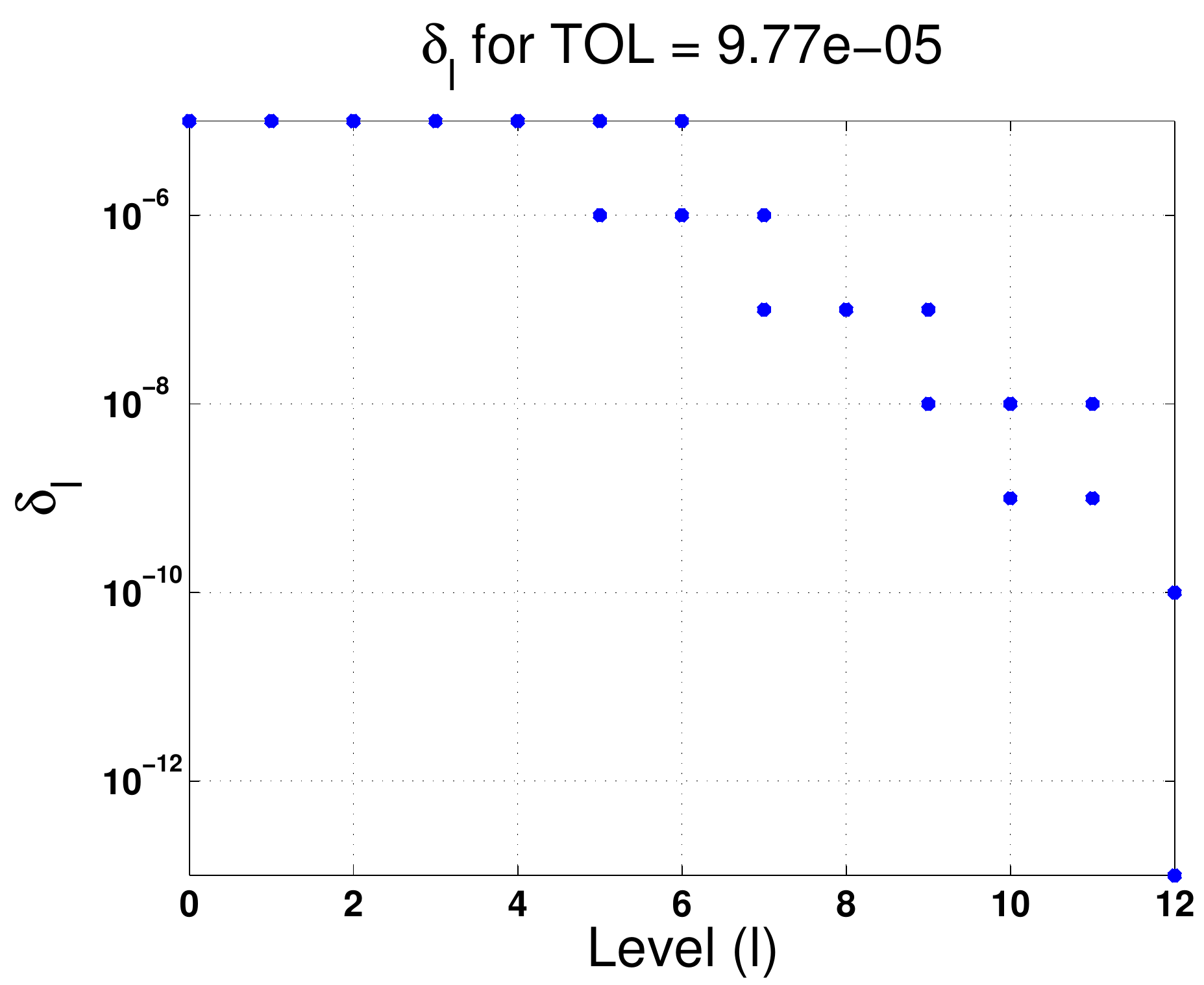}
\end{minipage}
\begin{minipage}{0.49\textwidth}
\includegraphics[scale=0.31]{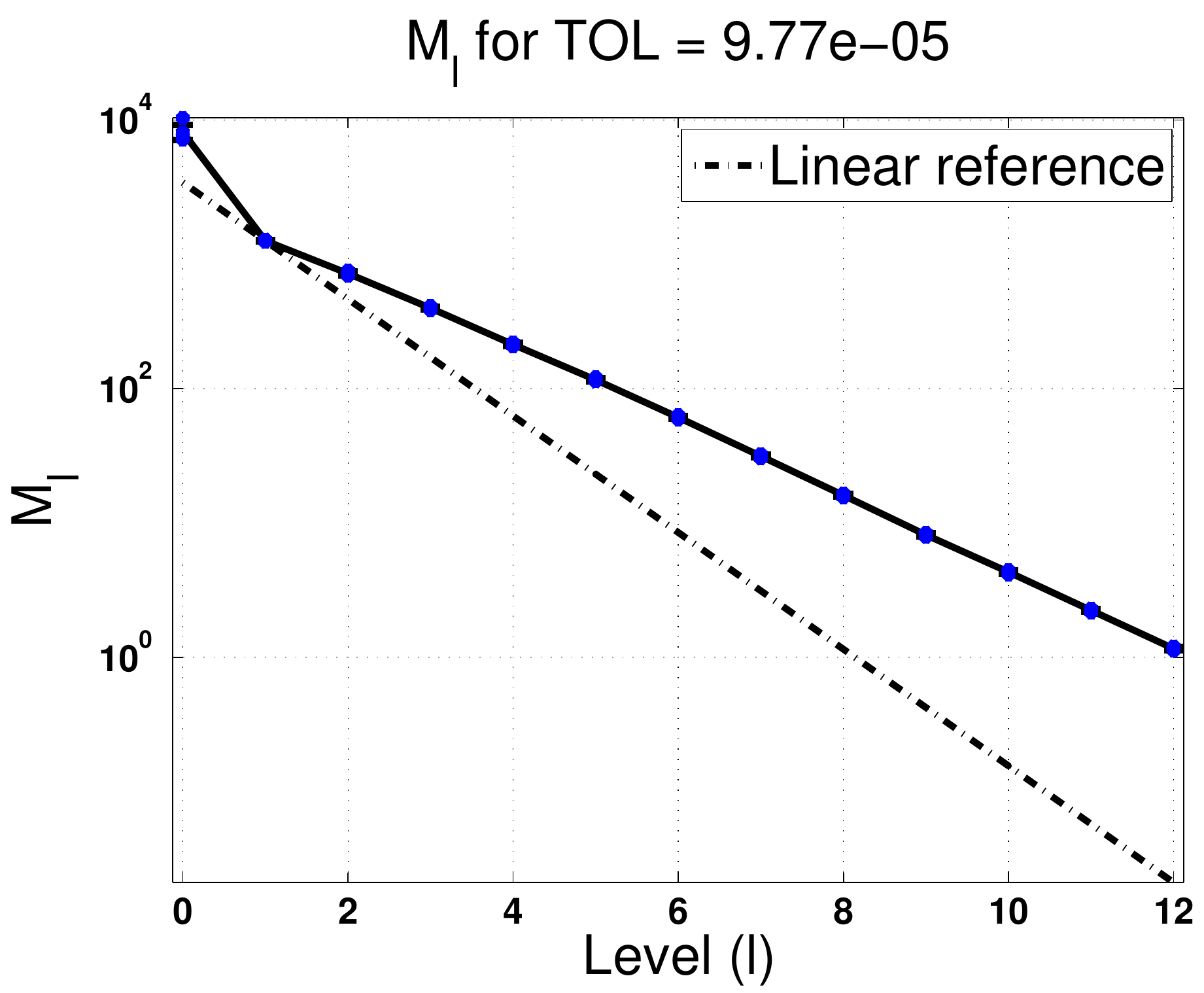}
\end{minipage}
\caption{One-step exit probability bound, $\delta_\ell$, and $M_\ell$ for $\ell{=}0,1,...,L^*$, for the smallest tolerance, for the simple decay model \eqref{ex:dec}.
}
\label{fig:dec2-out}
\end{figure}

In Figure \ref{fig:dec2-out}, we show the main outputs of Algorithm \ref{alg:Cal}, $\delta_\ell$ and $M_\ell$ for $\ell=0,...,L^*$, for the smallest considered tolerance. In this case, $L^*$ is 12. We observe that the number of realizations decreases slower than linearly, from levels $1$ to $L^*{-}1$, until it reaches $M_{L^*}{=}1$.

\begin{table}[h!]
\centering
\begin{tabular}{l|lll|lll|lll}
$TOL$ & $L^*$ & Min & Max & $\frac{\hat{W}_{ML}}{\hat{W}_{\ssa}}$ & Min & Max & $\frac{W_{ML}}{W_{\ssa}}$ & Min & Max \\ \noalign{\smallskip} \hline\noalign{\smallskip} 
3.13e-03 & 5 & 5 &5 & 0.03 &0.02 &0.04 & 0.03 &0.02 &0.05 \\ 
1.56e-03 & 6 & 6 &6 & 0.04 &0.02 &0.10 & 0.04 &0.02 &0.13 \\ 
7.81e-04 & 8 & 8 &8 & 0.03 &0.02 &0.05 & 0.03 &0.02 &0.06 \\ 
3.91e-04 & 9.2 & 9 &10 & 0.02 &0.02 &0.03 & 0.02 &0.01 &0.03 \\ 
1.95e-04 & 11 & 11 &11 & 0.02 &0.02 &0.03 & 0.02 &0.02 &0.04 \\ 
9.77e-05 & 12 & 12 &12 & 0.03 &0.02 &0.03 & 0.03 &0.02 &0.03 \\ 
\noalign{\smallskip}\hline 
\end{tabular}

\bigskip
\caption{Details of the ensemble run of Algorithm \ref{alg:Cal} for the simple decay model \eqref{ex:dec}.
As an example, the second row of the table indicates that, for a tolerance $TOL{=}1.56\cdot 10^{-3}$, six levels are needed. 
The predicted work of the multilevel hybrid method is, on average, $4\%$ of the predicted work of the SSA method, which coincides with the actual work. {Observed minimum and maximum values in the ensemble are also provided}.
}
\label{tab:dec2}
\end{table}

\begin{figure}[h!]
\centering
\begin{minipage}{0.49\textwidth}
\includegraphics[scale=0.31]{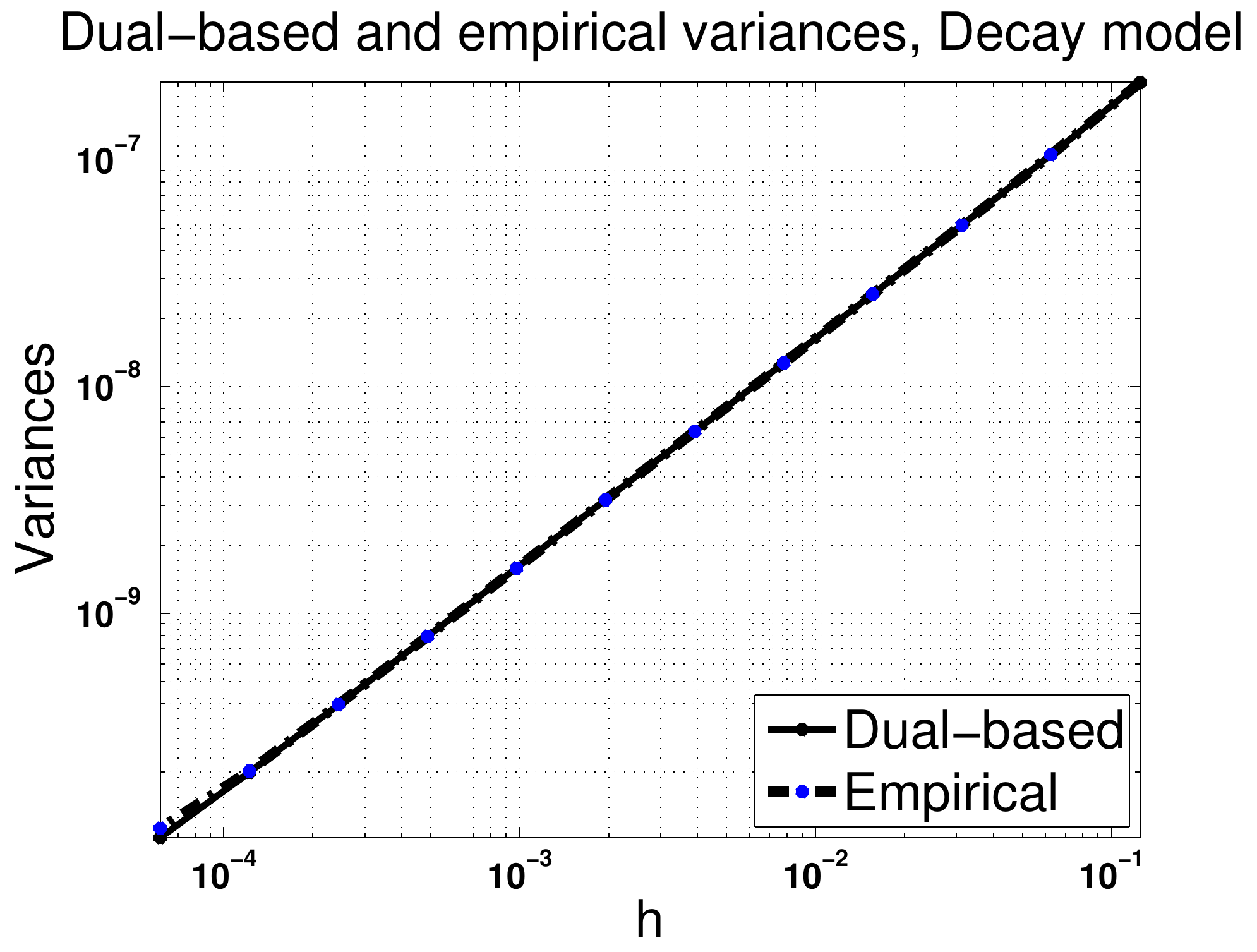}
\end{minipage}
\begin{minipage}{0.49\textwidth}
\includegraphics[scale=0.31]{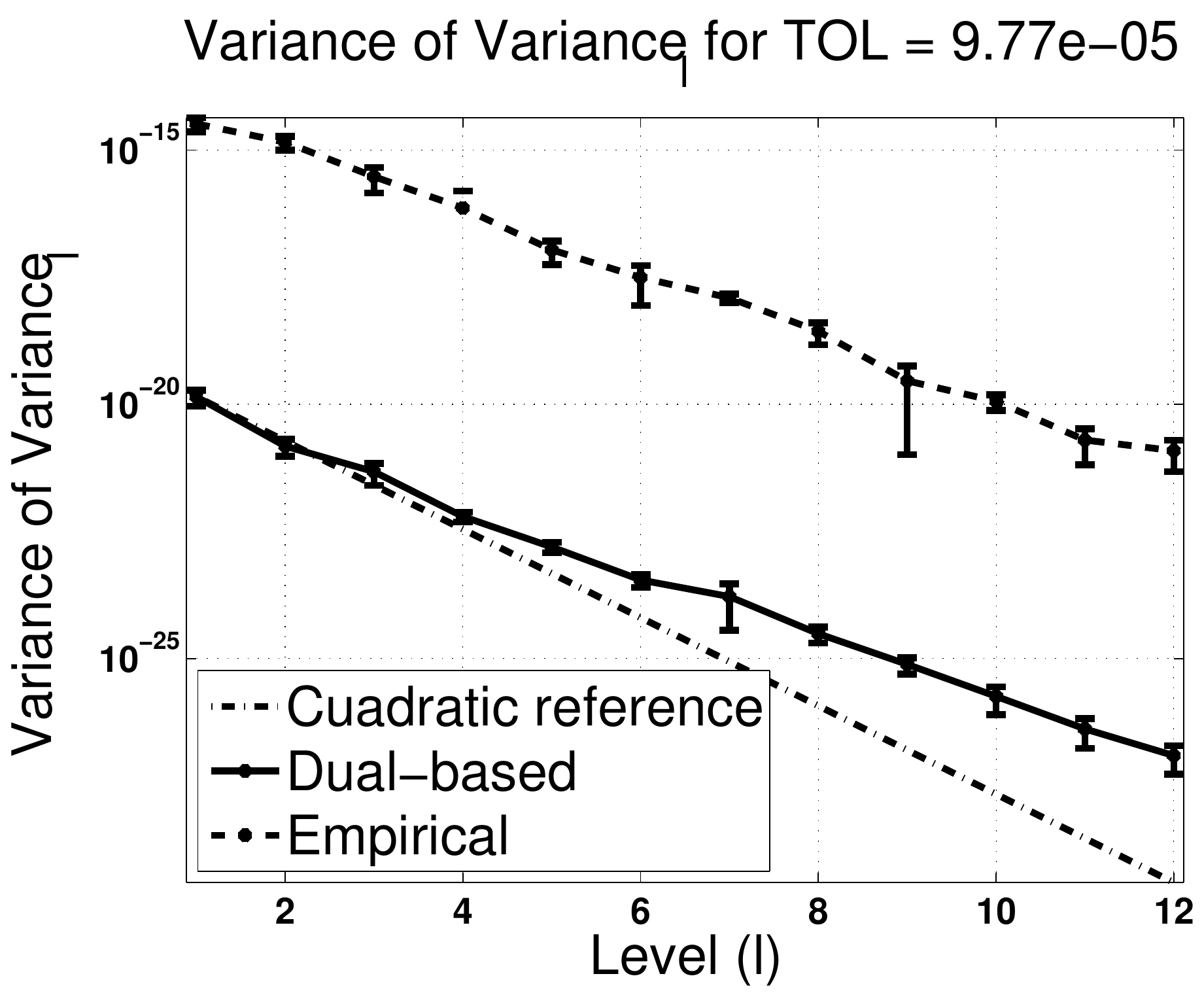}
\end{minipage}
\caption{Left: performance of the formula \eqref{eq:varhatestimated} as a strong error estimate, for the simple decay model \eqref{ex:dec}. Here, $h{=}\Delta t$. Right: estimated variance of $\mV{\ell}$ with 95\% confidence intervals.
}
\label{fig:effdec2}
\end{figure}

In the left panel of Figure \ref{fig:effdec2}, we show the performance of  formula \eqref{eq:varhatestimated}, implemented in Algorithm \ref{alg:varggl}, used to estimate the strong error, $\mV{\ell}$, {defined in Section \ref{calibration}}. The quotient of $\hmV{\ell}$ over a standard Monte Carlo estimate of $\mV{\ell}$ is almost 1 for the first ten levels. At levels 11 and 12, we obtain  0.99 and  0.91, respectively. Both quantities are estimated using a coefficient of variation less than 5\%, but there is a remarkable difference in terms of computational work in favor of our dual-weighted estimator. 
In the right panel of the same figure, we show the estimated variance of $\mV{\ell}$, computed by dual-weighted estimation \eqref{eq:varhatestimated} and computed by direct sampling. Observe that, in this case, the computational savings may be up to order $\Ordo{10^5}$.

\begin{figure}[h!]
\centering
\begin{minipage}{0.49\textwidth}
\includegraphics[scale=0.31]{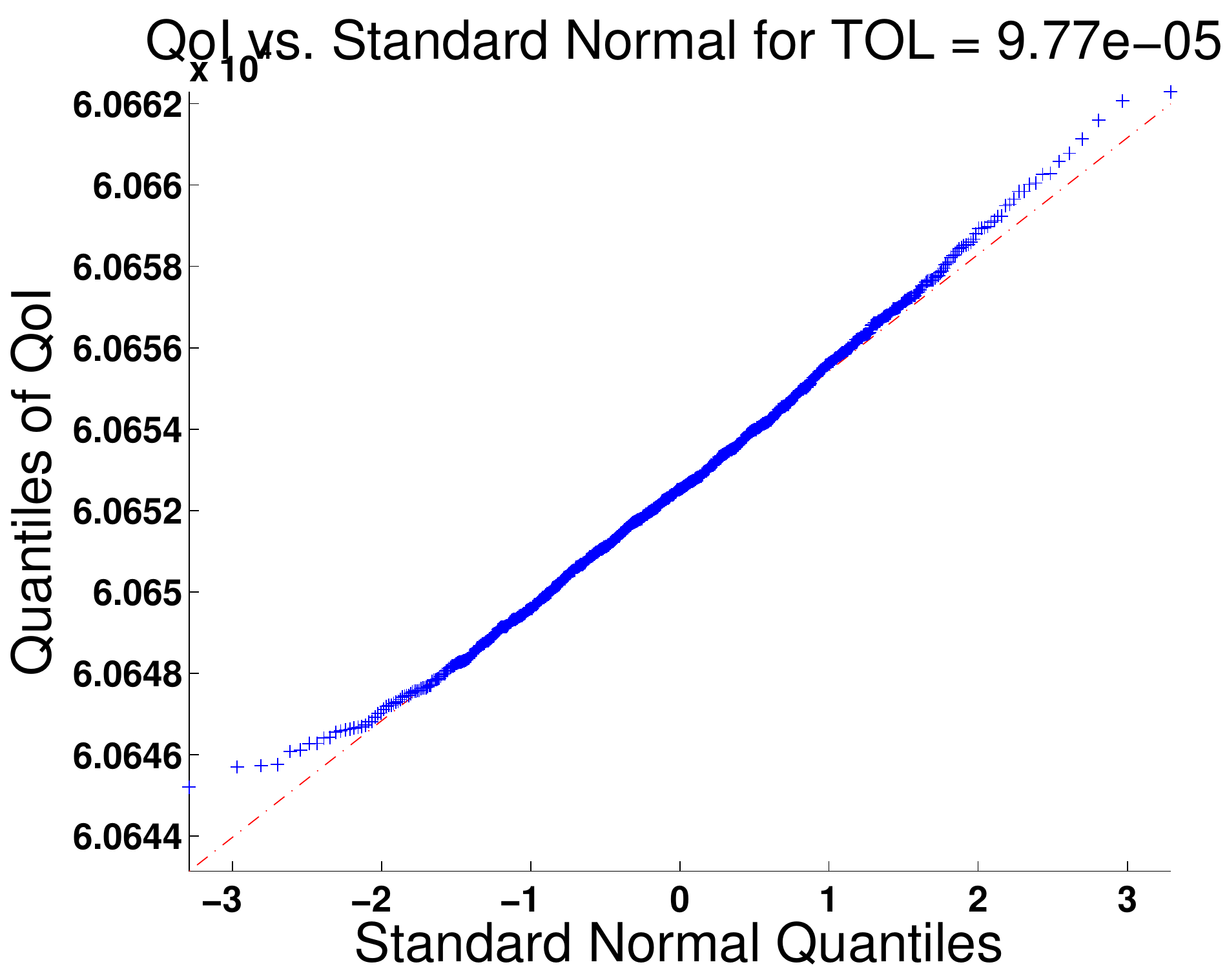}
\end{minipage}
\begin{minipage}{0.49\textwidth}
\includegraphics[scale=0.31]{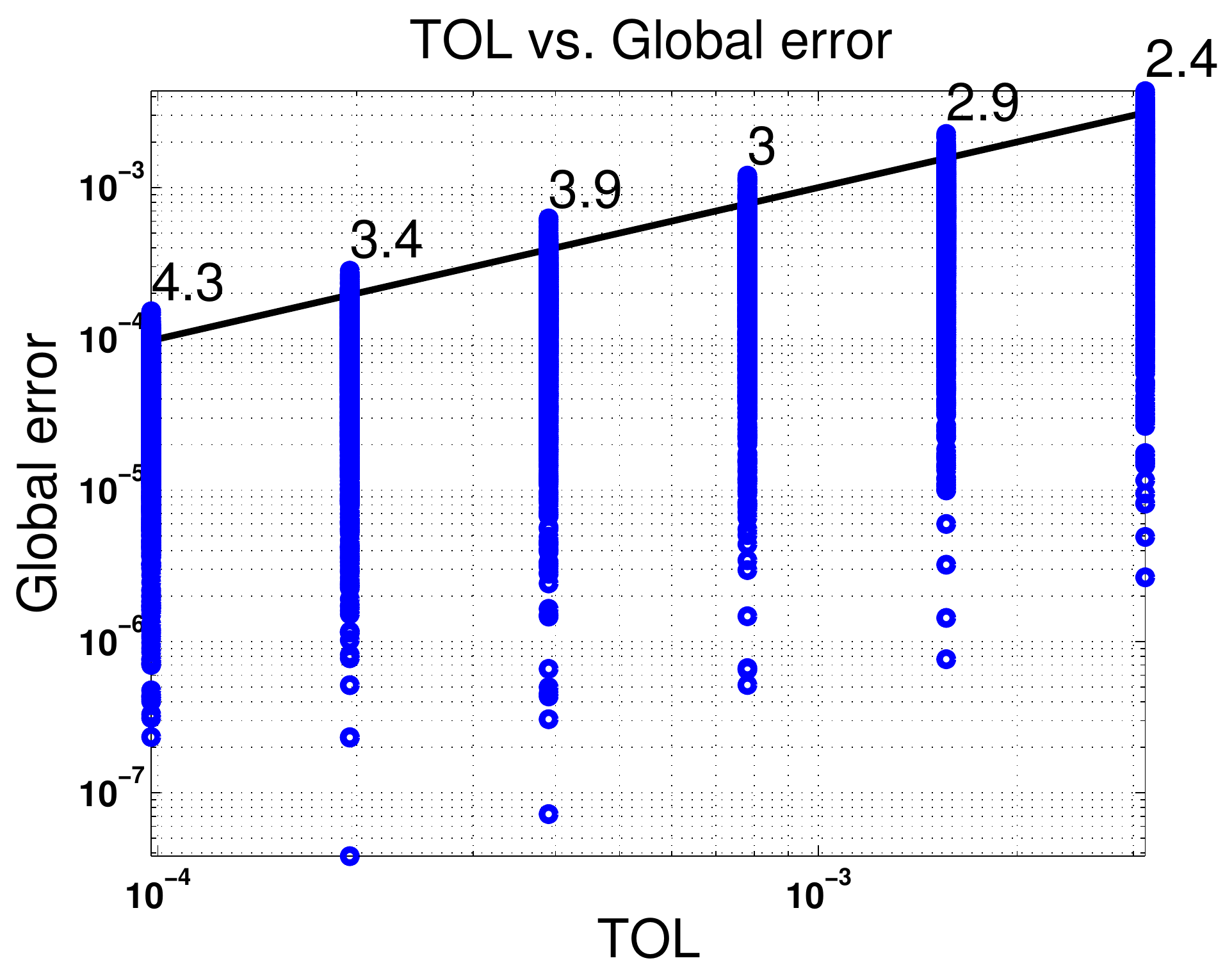}
\end{minipage}
\caption{Left: QQ-plot for the hybrid Chernoff MLMC estimates, $\mathcal{M}_L$, in the simple decay model \eqref{ex:dec}. {Also, we performed a Shapiro-Wilk normality test, and we obtained a p-value of $0.0105$}. Right: $TOL$ versus the actual computational error. The numbers above the straight line show the percentage of runs that had errors larger than the required tolerance. We observe that in all cases, except for the smallest tolerance, the computational error follows the imposed tolerance with the expected confidence of 95\%.}
\label{fig:qqdec2}
\end{figure}

In the simulations, we observed that, as we refine $TOL$, the optimal number of levels approximately increases logarithmically, which is a desirable feature. We fit the model $L^*=a+b\log(TOL^{-1})$, obtaining $b{=}{-}2.11$ and $a{=}{-}7.3$.

The QQ-plot in Figure \ref{fig:qqdec2} shows, for the smallest considered $TOL$, $10^3$ independent realizations of the multilevel estimator, $\mathcal{M}_L$ (defined by \eqref{MLMCest}). Those $10^3$ points are generated using 5 sets of parameters given by an independent run of the calibration algorithm (Algorithm \ref{alg:Cal}). This plot, complemented with a Shapiro-Wilk normality test, validates our assumption about the Gaussian distribution of the statistical error. Observe that the estimates are concentrated around the theoretical value $X_0 \exp(-c(T{-}t)) = 10^5 \exp(-0.5)\approx 6.0653e+04$. In the same figure, we also show $TOL$ versus the actual computational error. It can be seen that the prescribed tolerance is achieved with the required confidence of 95\%, in all the tolerances.

\subsection{Gene Transcription and Translation \cite{Anderson2012}} This model has five reactions,
\begin{align}
\label{ex:gtt}
\nonumber \emptyset \xrightarrow{c_1} R,& \ \ R \xrightarrow{c_2} R+P \\
2P \xrightarrow{c_3}  D,& \ \ R \xrightarrow{c_4} \emptyset \\
\nonumber P \xrightarrow{c_5} \emptyset 
\end{align}
described respectively by the stoichiometric matrix and the propensity function
\begin{align*}
\nu = \left( 
 \begin{array}{cccc}     1    & 0 &    0 \\
     0    & 1 &    0 \\
     0   & -2 &    1 \\
    -1    & 0 &    0 \\
     0   & -1 &    0 
 \end{array} 
 \right) \mbox{   and   }  a(X) = \left( \begin{array}{c}  c_1 \\ c_2 R \\ c_3 P(P{-}1)\\ c_4R \\ c_5 P  \end{array} \right)\COMMA
\end{align*}
where $X(t)=(R(t),P(t),D(t))$, and $c_1 {=} 25$, $c_2 {=}10^3$, $c_3{=}0.001$, $c_4{=}0.1$, and $c_5{=}1$. In the simulations, the initial condition is $(0, 0,0)$ and the final time is $T{=}1$. The observable is given by $g(X)=D$.
We observe that 
the abundance of the mRNA species, represented by $R$, is close to zero for $t \in [0,T]$. However, as we point out in \cite{ourSL}, the reduced abundance of one of the species is not enough to ensure that the SSA method should be used. 

\begin{figure}[h!]
\centering
\begin{minipage}{0.49\textwidth}
\includegraphics[scale=0.31]{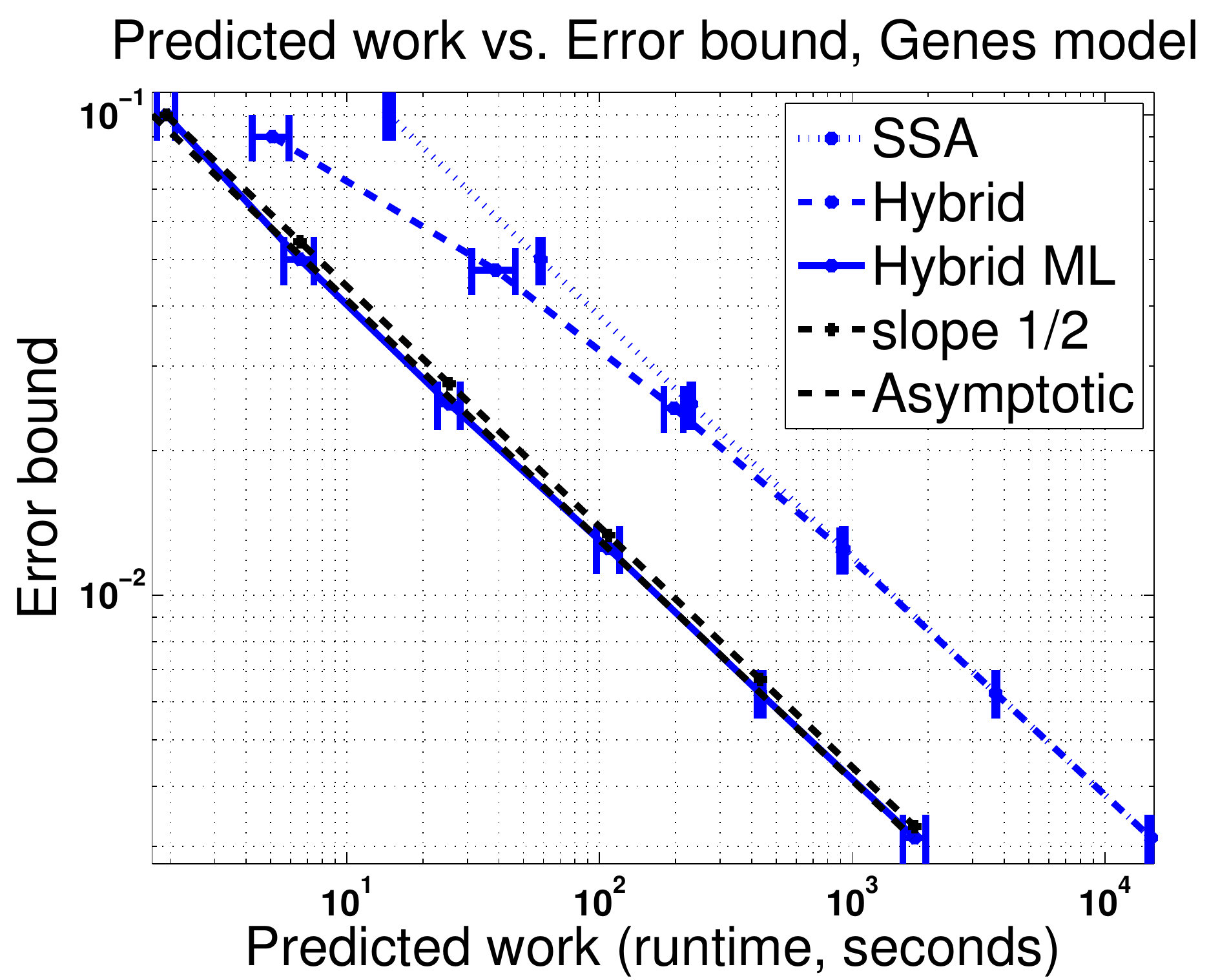}
\end{minipage}
\begin{minipage}{0.49\textwidth}
\includegraphics[scale=0.31]{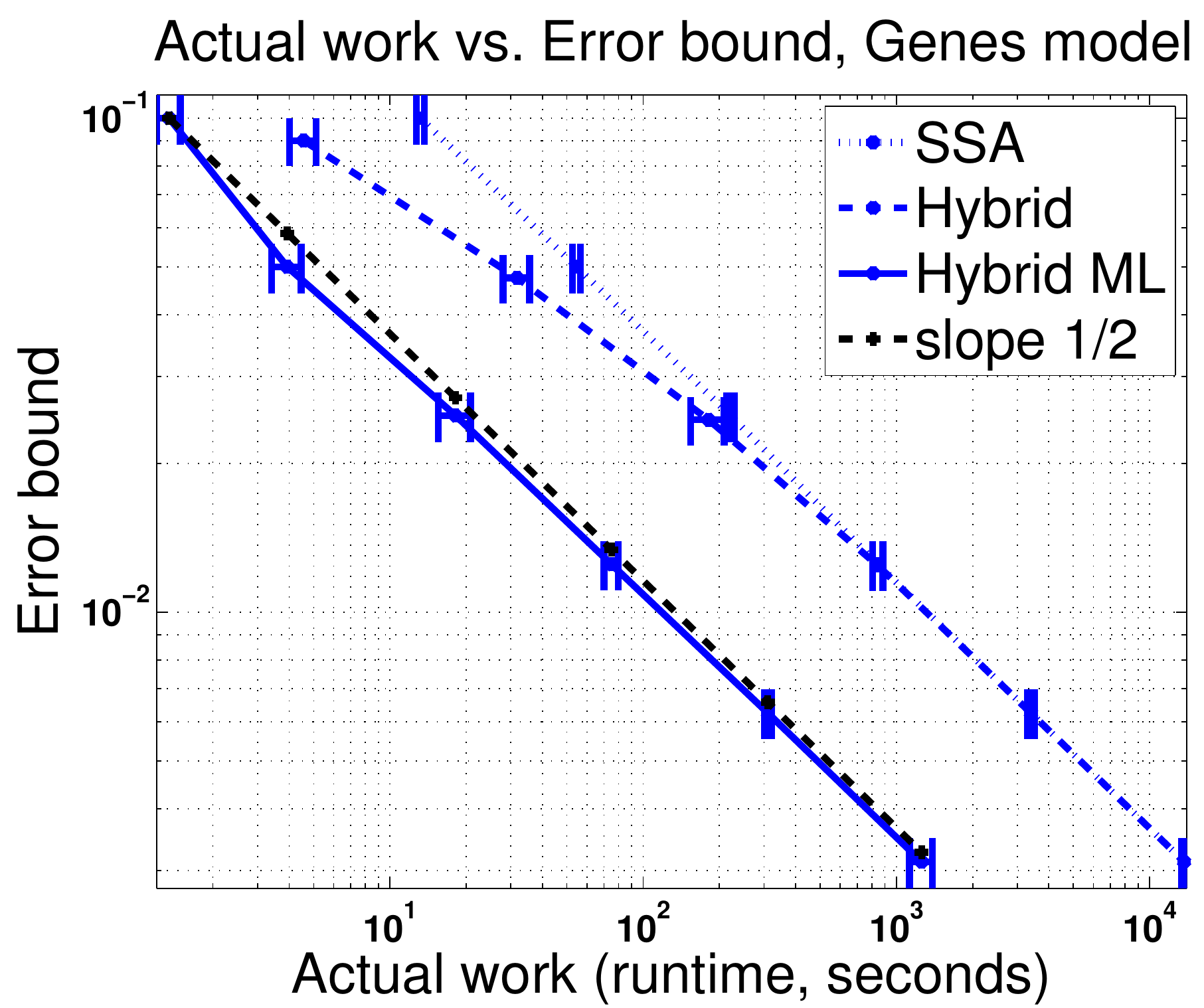}
\end{minipage}
\caption{Left: Predicted work (runtime) versus the estimated error bound for the gene transcription and translation model \eqref{ex:gtt}. The hybrid method is preferred over the SSA for the first three tolerances only. 
The multilevel hybrid method is preferred over the SSA and the single-level method 
for all the tolerances. Right: Actual work (runtime) versus the estimated error bound.}
\label{fig:gtt-worktimes}
\end{figure}

We now analyze an ensemble of five independent runs of the calibration algorithm (Algorithm \ref{alg:Cal}), using different relative tolerances. In Figure \ref{fig:gtt-worktimes}, we show, in the left panel, the total predicted work (runtime) for the single-level hybrid method, for the multilevel hybrid method and for the SSA method, versus the estimated error bound. 
We also show the estimated asymptotic work of the multilevel method.
Again, the multilevel hybrid method outperforms the others and we remark that the observed computational work of the multilevel method is of order $\Ordo{TOL^{-2}}$.

\begin{figure}[h!]
\centering
\begin{minipage}{0.49\textwidth}
\includegraphics[scale=0.31]{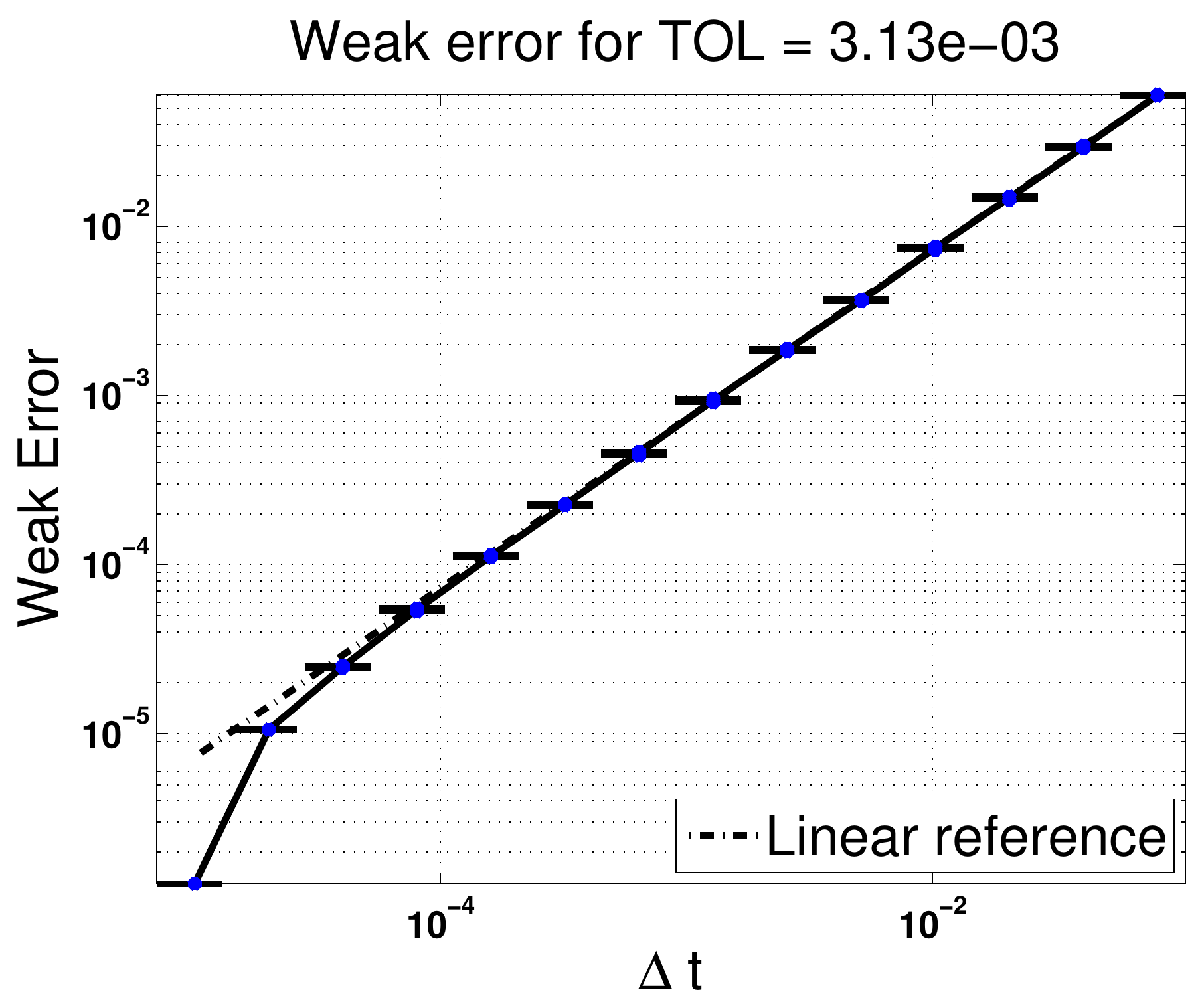}
\end{minipage}
\begin{minipage}{0.49\textwidth}
\includegraphics[scale=0.31]{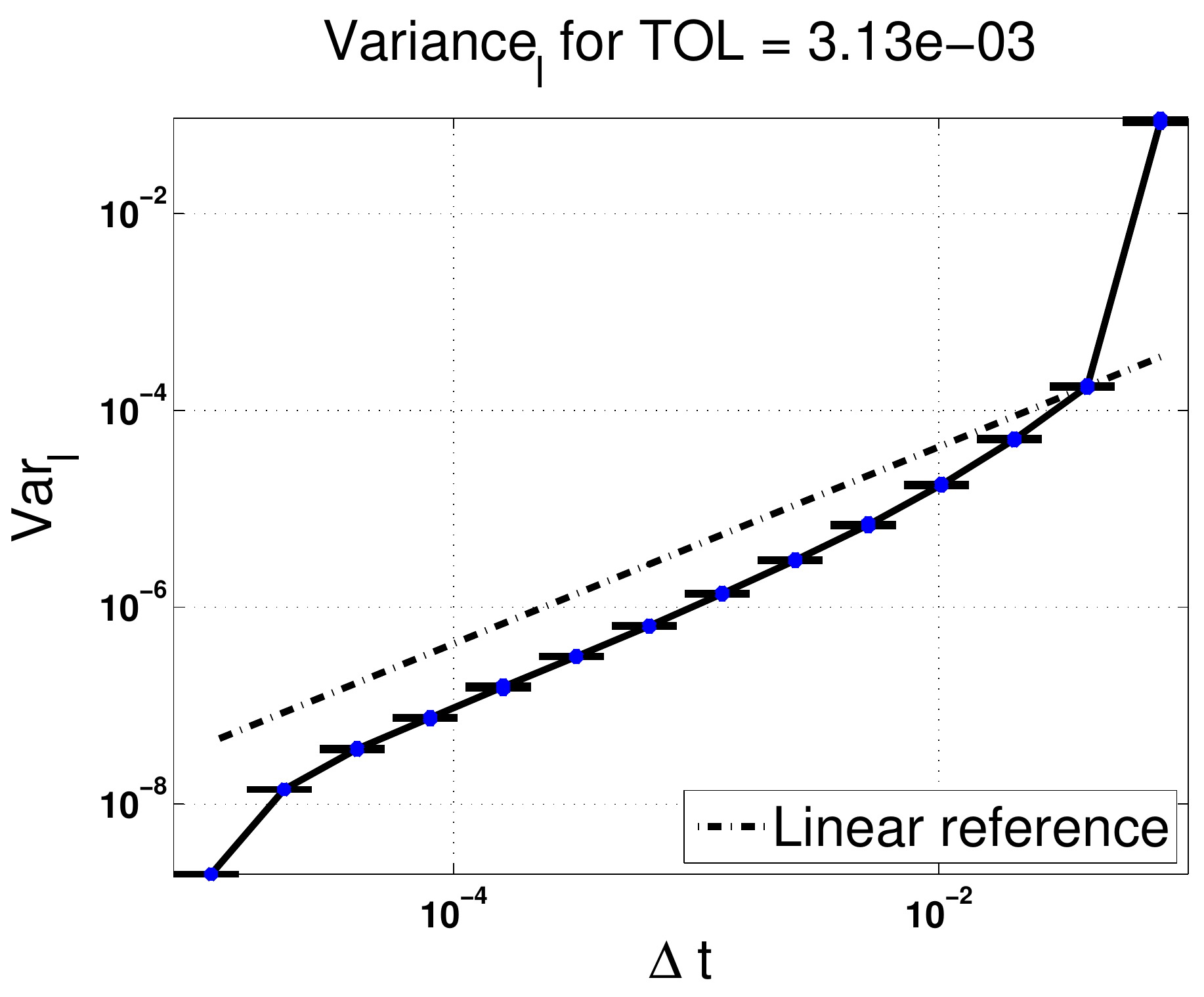}
\end{minipage}
\begin{minipage}{0.49\textwidth}
\includegraphics[scale=0.31]{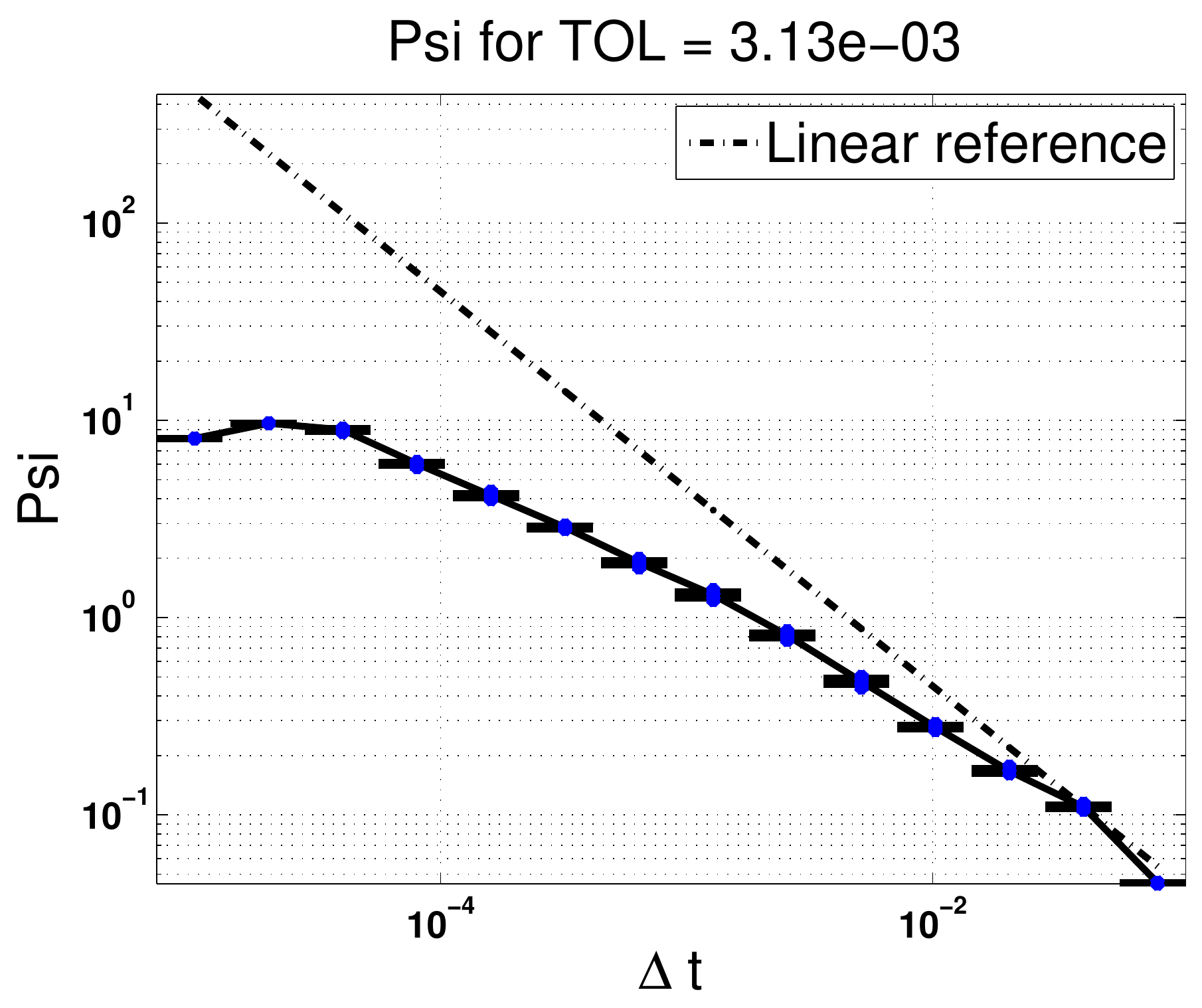}
\end{minipage}
\begin{minipage}{0.49\textwidth}
\includegraphics[scale=0.31]{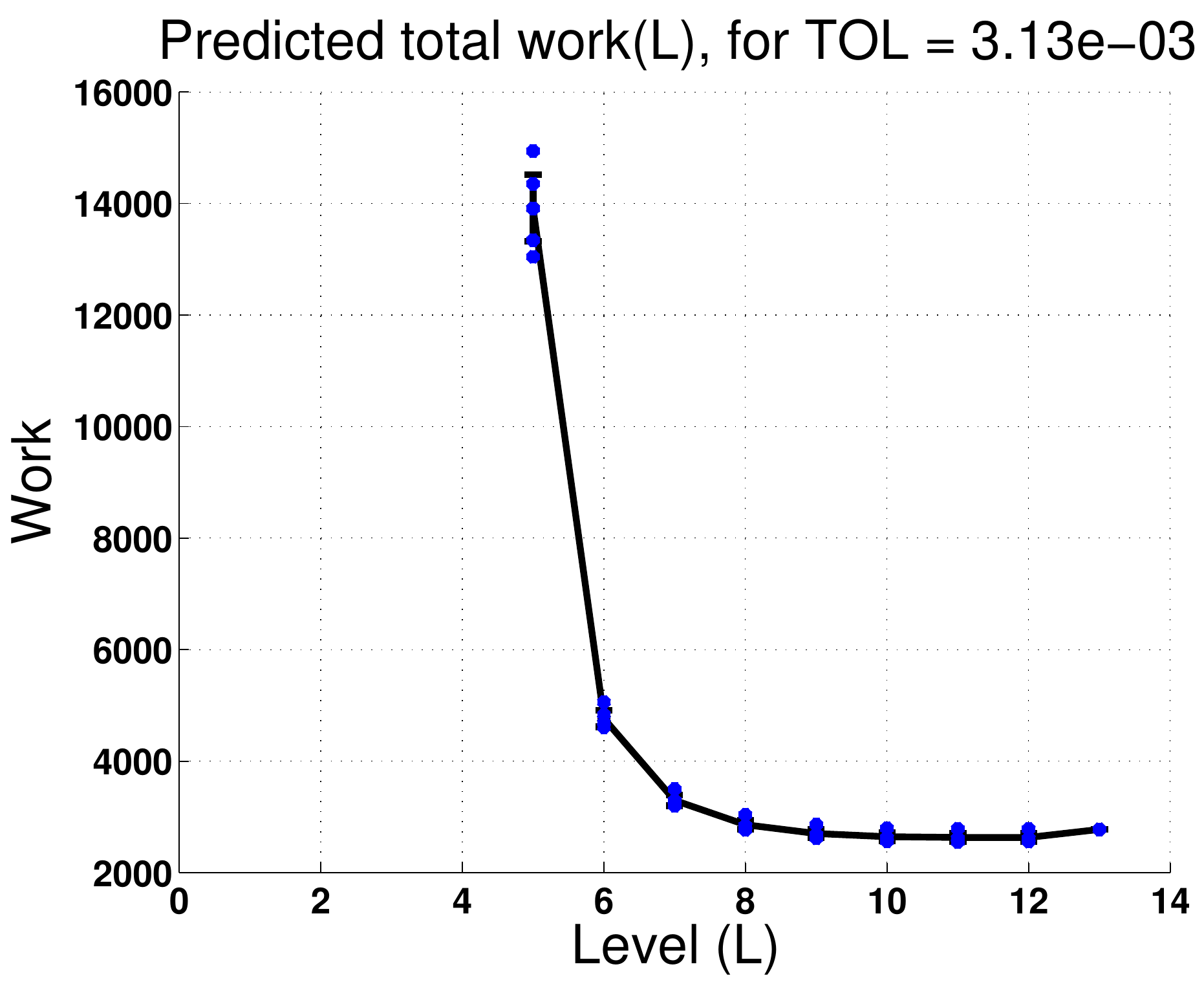}
\end{minipage}
\caption{Upper left: estimated weak error, $\WEH{\ell}$, as a function of the time mesh size, $\Delta t$, for the gene transcription and translation model \eqref{ex:gtt}. Upper right: estimated variance of the difference between two consecutive levels, $\hmV{\ell}$, as a function of $\Delta t$. Lower left: estimated path work, $\hat \psi_\ell$, as a function of $\Delta t$. Lower right: estimated total computational work, $\sum_{l=0}^L \hat \psi_l M_l$, as a function of the level, $L$.
}
\label{fig:gtt-diag}
\end{figure}

In Figure \ref{fig:gtt-diag}, we can observe how the estimated weak error decreases linearly for the coarser time meshes, but, as we continue refining the time mesh, it quickly decreases towards zero. In the case of the estimated variance, $\hmV{\ell}$, it decreases faster than linearly, and it also quickly decreases towards zero afterwards.
This is a consequence of the transition from a hybrid regime to a pure exact one.
The estimated total path work, $\hat \psi_{\ell}$, increases sublinearly as we refine the mesh. 
Note that $\hat \psi_{\ell}$ reaches a maximum, which corresponds to a SSA-dominant regime.
In the lower right panel, we show the total computational work only in the cases in which $\WEH{\ell} < TOL{-}TOL^2$. 

\newcommand{\TOLGTT}{3.13e-03}

\begin{figure}[h!]
\centering
\begin{minipage}{0.49\textwidth}
\includegraphics[scale=0.31]{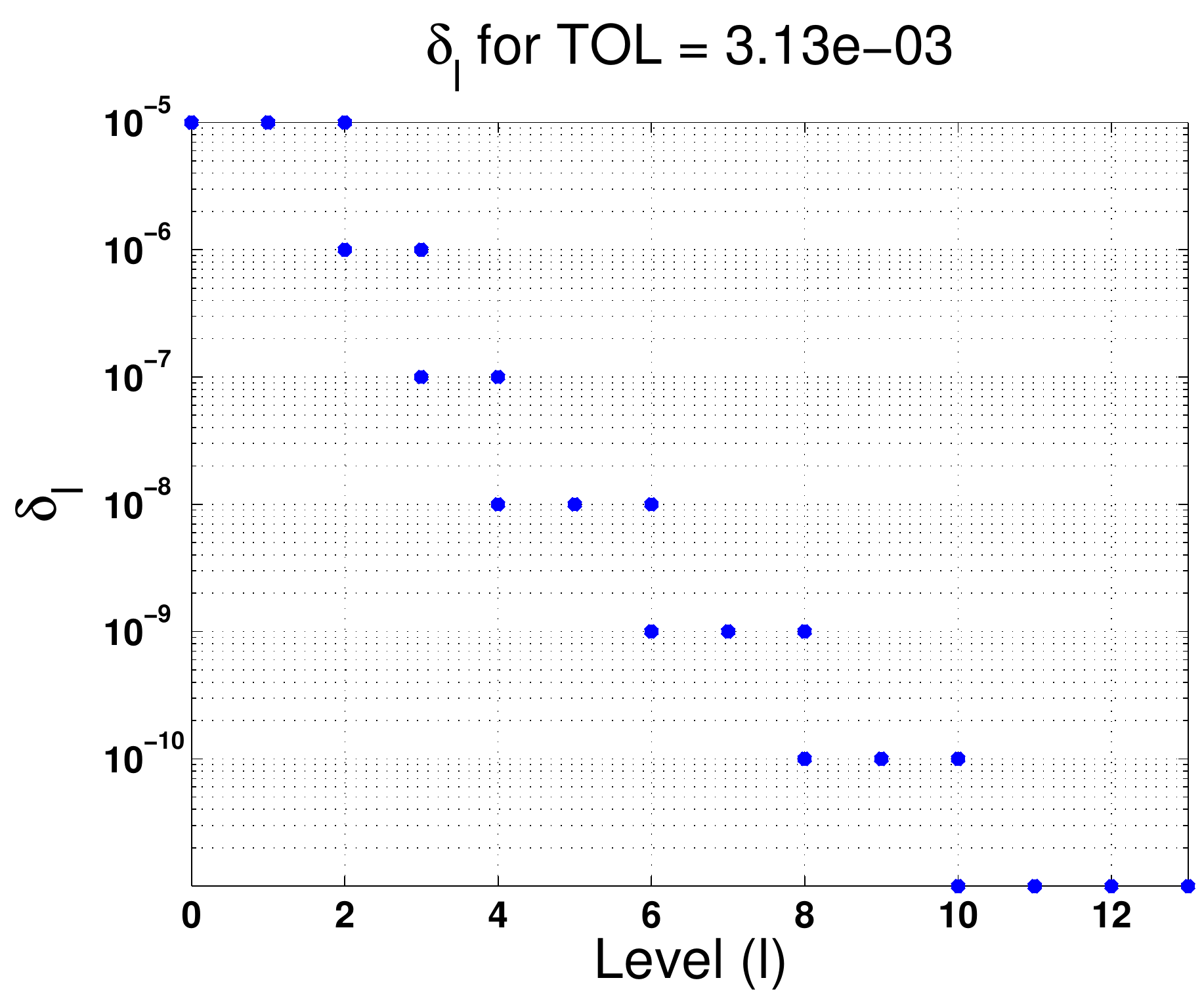}
\end{minipage}
\begin{minipage}{0.49\textwidth}
\includegraphics[scale=0.31]{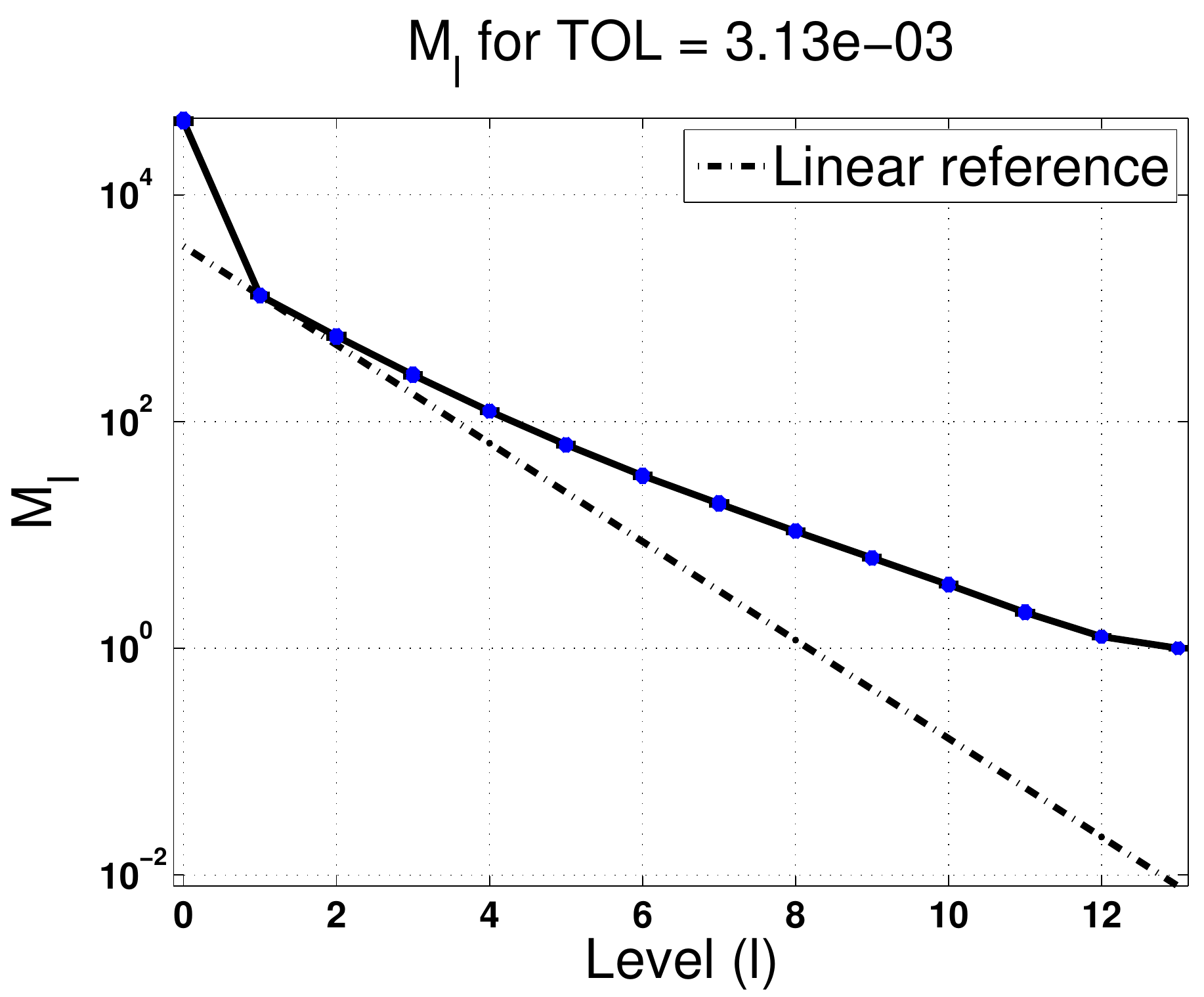}
\end{minipage}
\caption{The one-step exit probability bound, $\delta_\ell$, and $M_\ell$ for $\ell{=}0,1,...,L^*$, for the smallest tolerance in the gene transcription and translation model \eqref{ex:gtt}.
}
\label{fig:gtt-out}
\end{figure}

In Figure \ref{fig:gtt-out}, we show the main outputs of Algorithm \ref{alg:Cal}, $\delta_\ell$ and $M_\ell$ for $\ell=0,...,L^*$, for the smallest tolerance. We observe that the number of realizations decreases slower than linearly from levels $1$ to $12$.

\begin{table}[h!]
\centering
\begin{tabular}{l|lll|lll|lll}
$TOL$ & $L^*$ & Min & Max & $\frac{\hat{W}_{ML}}{\hat{W}_{\ssa}}$ & Min & Max & $\frac{W_{ML}}{W_{\ssa}}$ & Min & Max \\ \noalign{\smallskip} \hline\noalign{\smallskip} 
1.00e-01 & 3 & 3 &3 & 0.04 &0.04 &0.04 & 0.06 &0.05 &0.07 \\ 
5.00e-02 & 4.6 & 4 &5 & 0.04 &0.03 &0.04 & 0.05 &0.05 &0.05 \\ 
2.50e-02 & 6 & 6 &6 & 0.03 &0.03 &0.04 & 0.05 &0.04 &0.05 \\ 
1.25e-02 & 8 & 8 &8 & 0.03 &0.03 &0.03 & 0.05 &0.05 &0.06 \\ 
6.25e-03 & 10 & 10 &10 & 0.03 &0.03 &0.03 & 0.05 &0.04 &0.05 \\ 
3.13e-03 & 11.4 & 11 &13 & 0.03 &0.03 &0.03 & 0.05 &0.04 &0.05 \\ 
\noalign{\smallskip}\hline 
\end{tabular}
\bigskip
\caption{Details for the ensemble run of Algorithm \ref{alg:Cal} for the gene transcription and translation model \eqref{ex:gtt}.}
\label{tab:gtt}
\end{table}

\begin{figure}[h!]
\centering
\begin{minipage}{0.49\textwidth}
\includegraphics[scale=0.31]{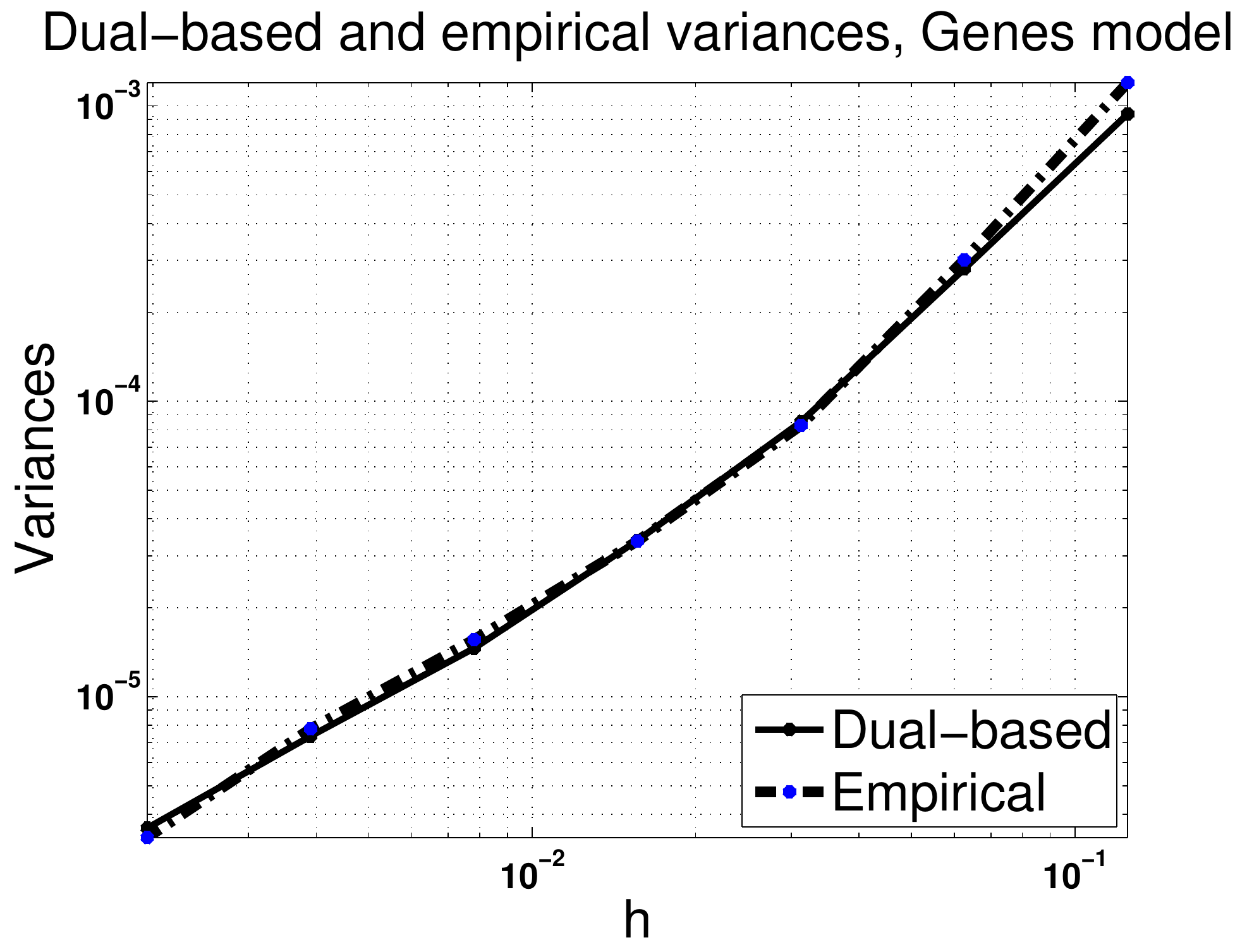}
\end{minipage}
\begin{minipage}{0.49\textwidth}
\includegraphics[scale=0.31]{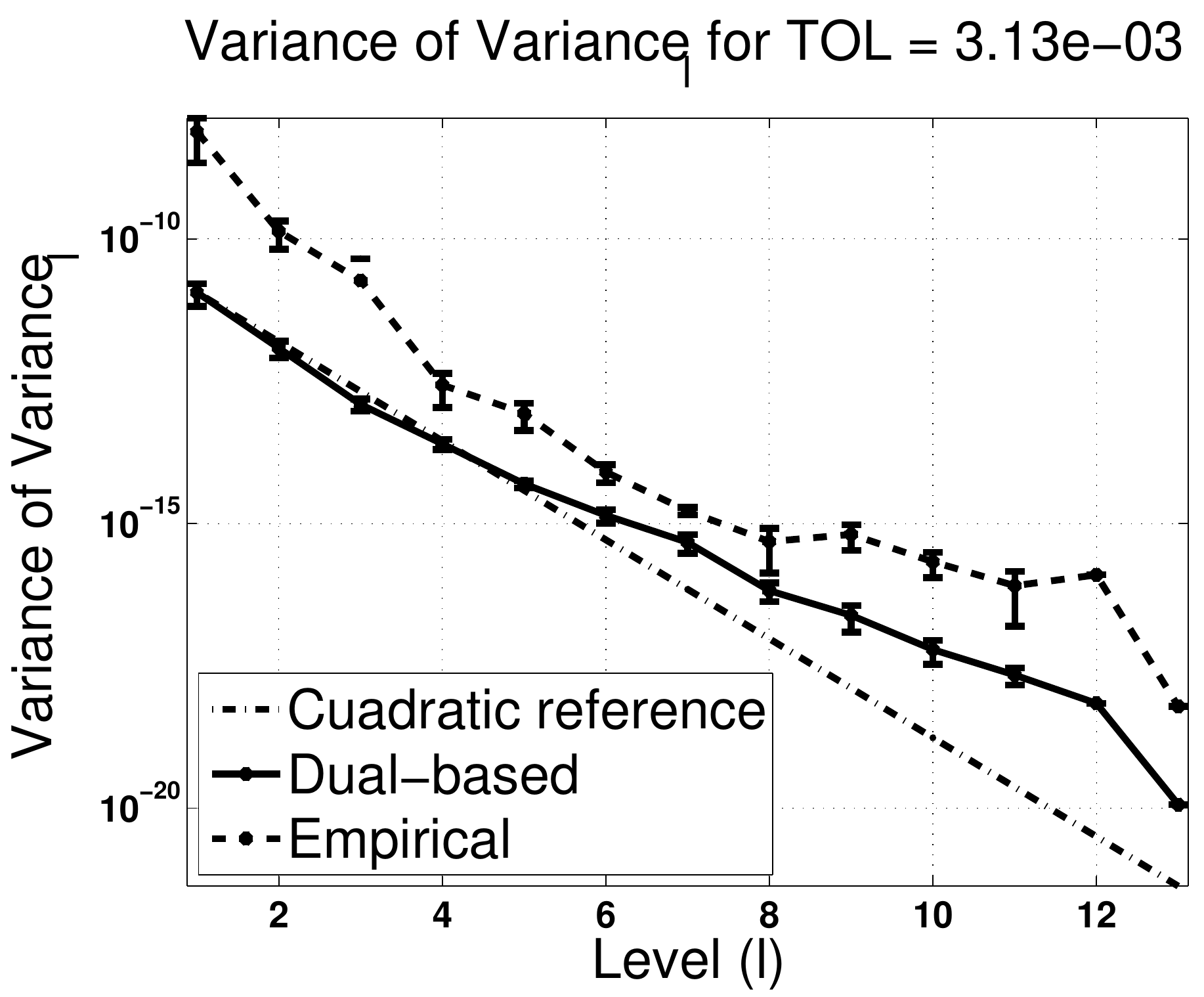}
\end{minipage}
\caption{Left: performance of formula \eqref{eq:varhatestimated} as a strong error estimate for the gene transcription and translation model \eqref{ex:gtt}. Here, $h{=}\Delta t$. Right: estimated variance of $\mV{\ell}$ with 95\% confidence intervals.}
\label{fig:effgtt}
\end{figure}

In Figure \ref{fig:effgtt}, we see that our dual-weighted estimator of the strong error, $\mV{\ell}$, gives essentially the same results as the standard Monte Carlo estimator, but with much less computational work. 
In this case, an accurately empirical estimate of $\mV{7}$ took almost 48 hours, but the dual-based computation of $\hmV{7}$ just took few minutes.

\begin{figure}[h!]
\centering
\begin{minipage}{0.49\textwidth}
\includegraphics[scale=0.31]{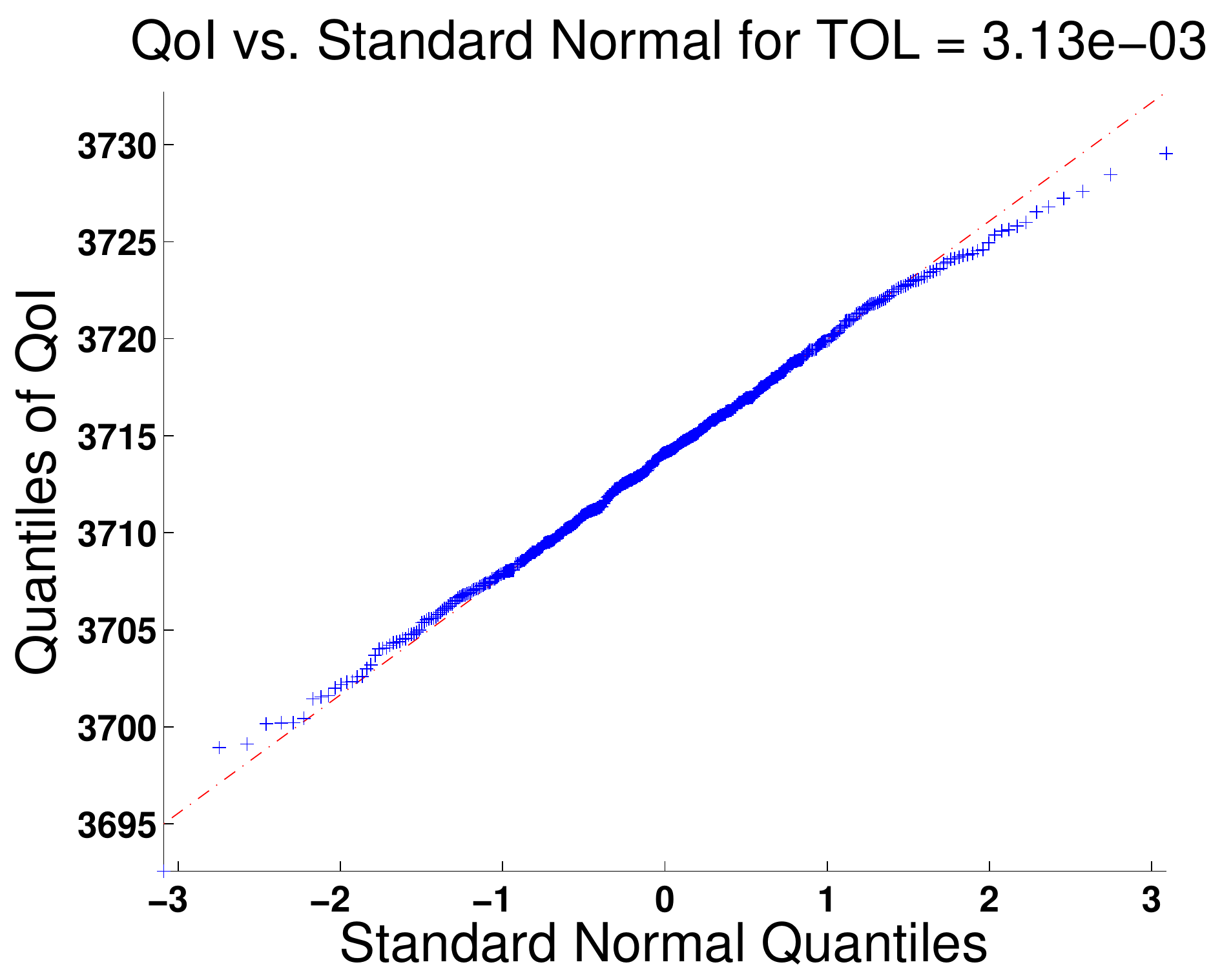}
\end{minipage}
\begin{minipage}{0.49\textwidth}
\includegraphics[scale=0.31]{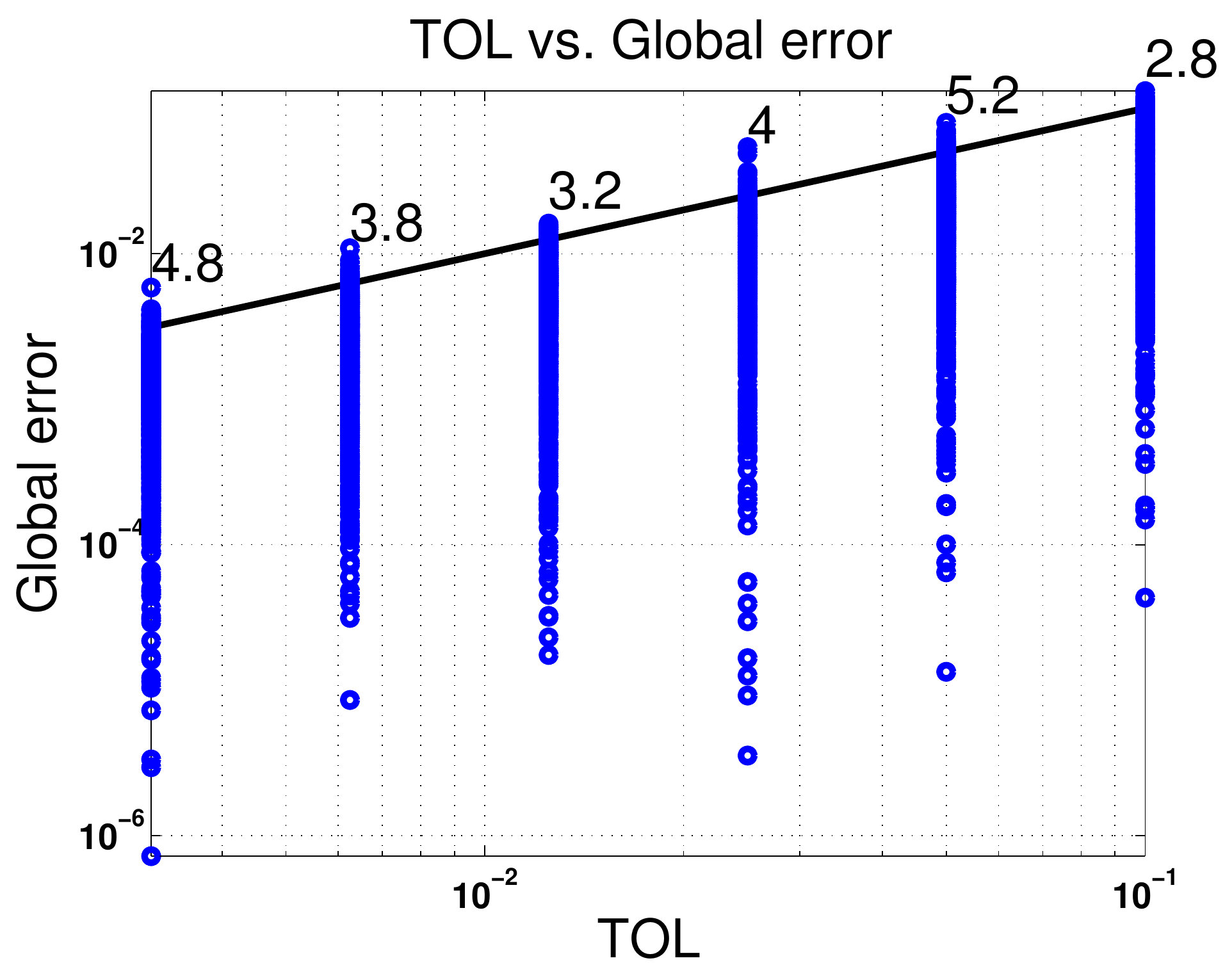}
\end{minipage}
\caption{Left: QQ-plot based on $\mathcal{M}_L$ estimates for the gene transcription and translation model \eqref{ex:gtt}. {Also, we performed a Shapiro-Wilk normality test and we obtained a p-value of $0.6$}. Right: $TOL$ versus the actual global computational error. The numbers above the straight line show the percentage of runs that had errors larger than the required tolerance. We observe that in all cases (except the second for a very small margin) the computational error follows the imposed tolerance with the expected confidence of 95\%.}
\label{fig:qqgtt}
\end{figure}

In the simulations, we observed that, as we refine $TOL$, the optimal number of levels approximately increases logarithmically, which is a desirable feature. We fit the model $L^*=a+b\log(TOL^{-1})$, obtaining $b{=}{-}2.48$ and $a{=}{-}2.85$.

The QQ-plot in the Figure \ref{fig:qqgtt}, computed in the same way as in the previous example, together with a Shapiro-Wilk normality test, shows the validity of the Gaussian assumption for the statistical errors. In the same figure, we also show $TOL$ versus the actual global computational error. It can be seen that the prescribed tolerance is achieved, except for the second smallest tolerance, with the required confidence of 95\%, since $C_A{=}1.96$.

\subsubsection*{MLMC Hybrid-Path Analysis}
\begin{figure}[h!]
\centering
\includegraphics[scale=0.5,angle=270]{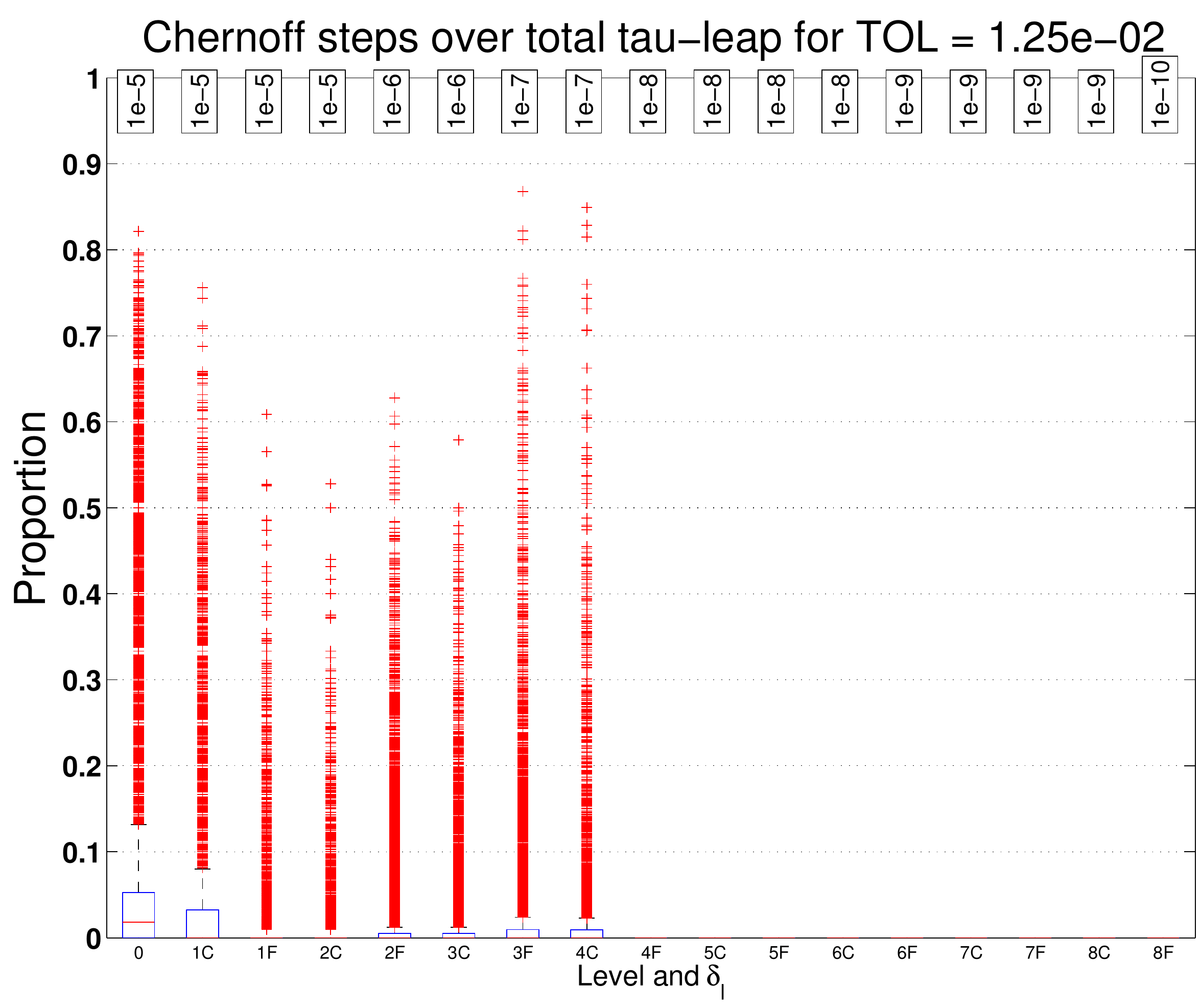}
\caption{Proportion of the number of Chernoff tau-leap steps over the total number of tau-leap steps for the gene transcription and translation model \eqref{ex:gtt}. In the x-axis, we show the corresponding level (starting from level 0) and, subsequently, the coarse (C) and fine (F) level. Below the title, we show the corresponding $\delta_\ell$ of each level. We observe a small increase in the proportion of the number of Chernoff steps from levels 1F/2C to levels 3F/4C (strictly speaking, a shift in the median and the third quartile). This is due to consecutive refinements in the values of $\delta$, from $1e{-}5$ to $1e{-}7$, producing smaller and smaller values of $\tau_{Ch}$.}
\label{fig:box1}
\end{figure}

\begin{figure}[h!]
\centering
\includegraphics[scale=0.5,angle=270]{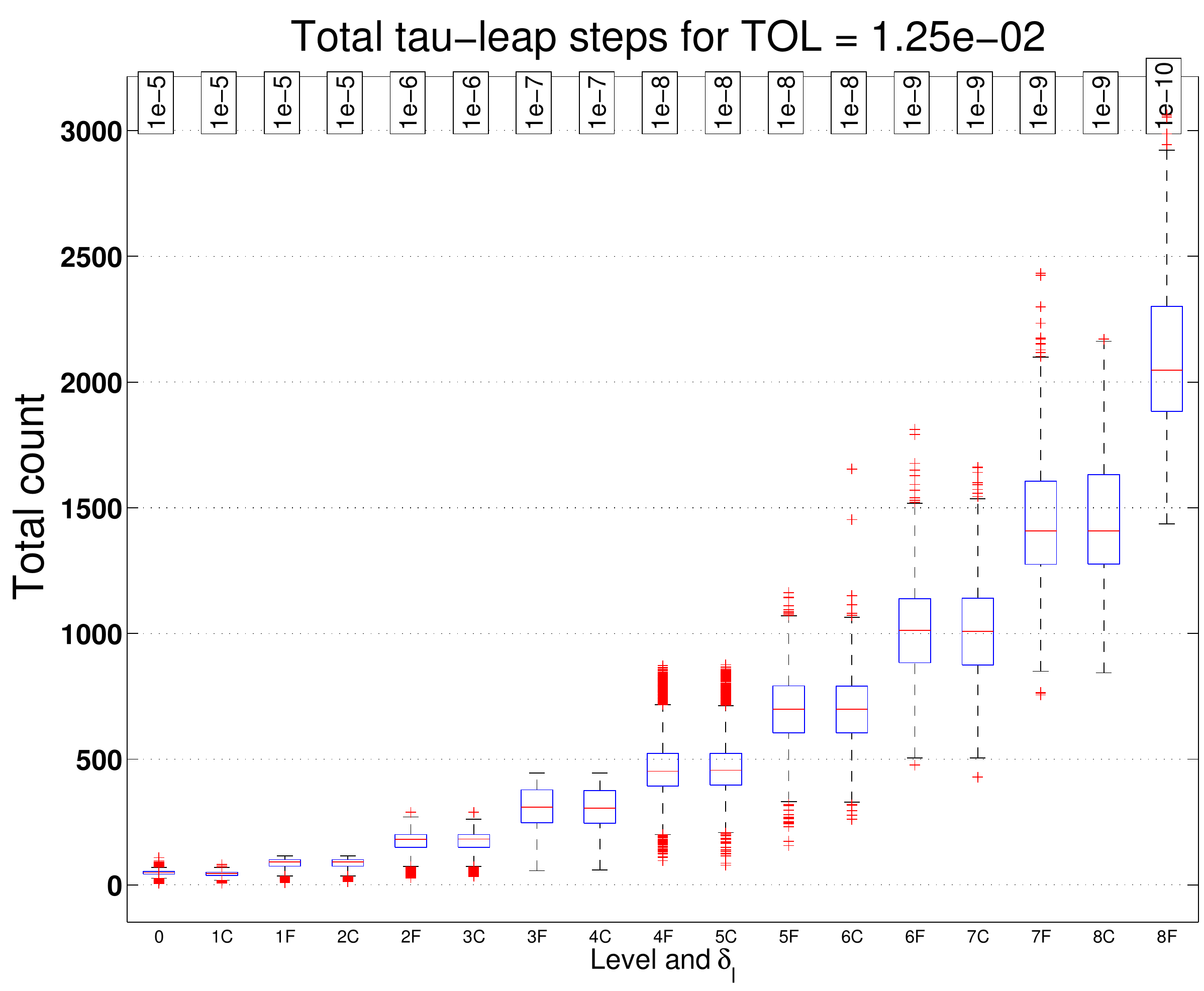}
\caption{Total number of tau-leap steps per path for the gene transcription and translation model \eqref{ex:gtt}. In the x-axis, we show the corresponding pairings of two consecutive levels (starting from level 0) and, subsequently, the coarse (C) and fine (F) meshes for two consecutive levels. Below the title, we show the corresponding $\delta_\ell$ of each level. The domain $I_{\TL}$ of the tau-leap method decreases with refinements, but, since the size of time mesh halves form by passing from one level to the next one, we see an increasing number of tau-leap steps until, at a certain level, there are no more tau-leap steps due to the relative computational cost of the tau-leap method.}
\label{fig:box2}
\end{figure}

\begin{figure}[h!]
\centering
\includegraphics[scale=0.5,angle=270]{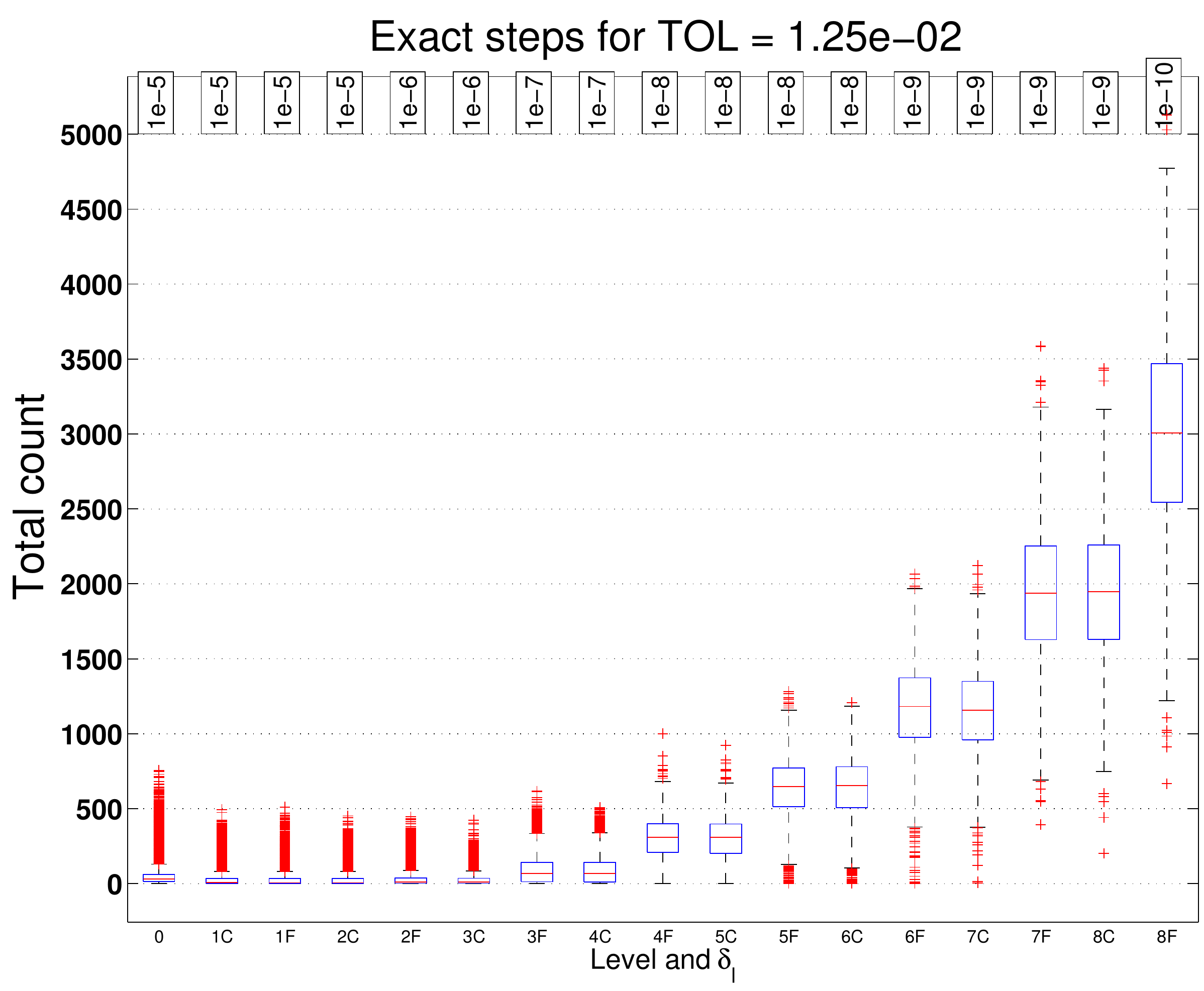}
\caption{Total number of exact steps per path for the gene transcription and translation model \eqref{ex:gtt}. In the x-axis, we show the corresponding pairings of two consecutive levels (starting from level 0) and, subsequently, the coarse (C) and fine (F) meshes for two consecutive levels. Below the title, we show the corresponding $\delta_\ell$ of each level.
The domain $I_{\MNRM}$ of the exact method is monotonically increasing with refinements of the time mesh and the one-step exit probability bound. As a consequence, we expect the total count of exact paths to be a monotonically increasing function of the level, $\ell$.}
\label{fig:box3}
\end{figure}

\begin{figure}[h!]
\centering
\begin{minipage}{0.32\textwidth}
\includegraphics[scale=0.24]{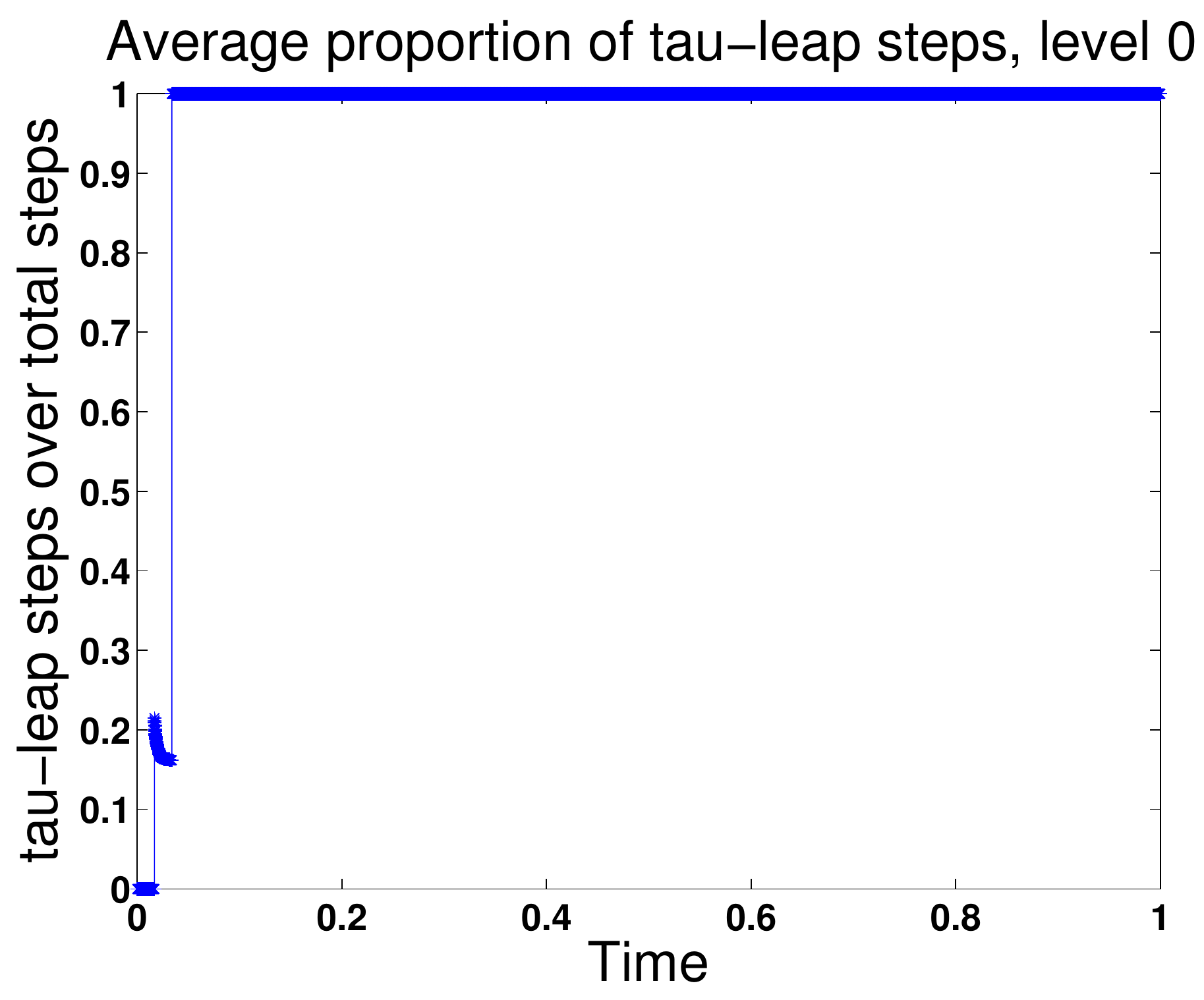}
\end{minipage}
\hfill
\begin{minipage}{0.32\textwidth}
\includegraphics[scale=0.24]{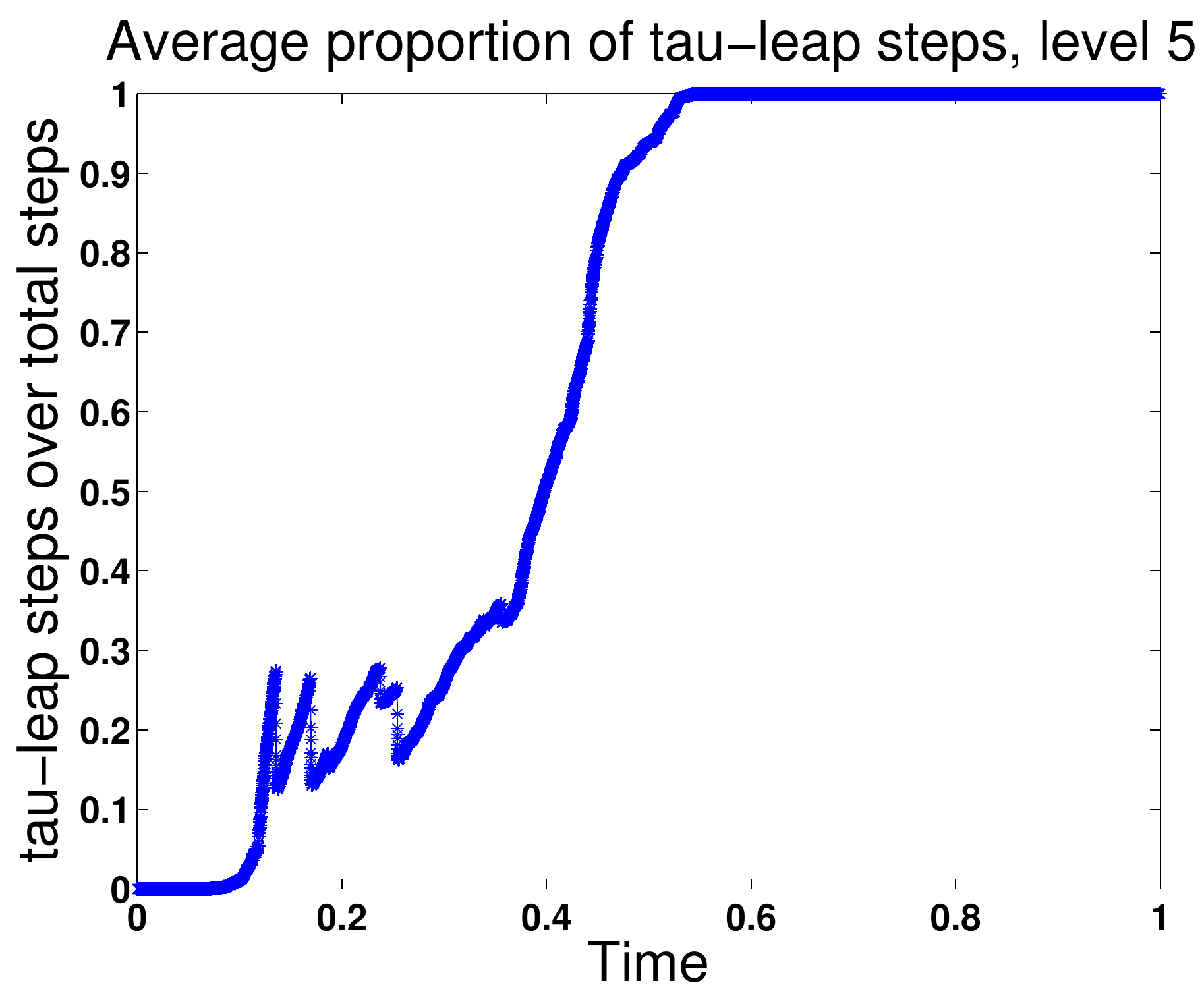}
\end{minipage}
\hfill
\begin{minipage}{0.32\textwidth}
\includegraphics[scale=0.24]{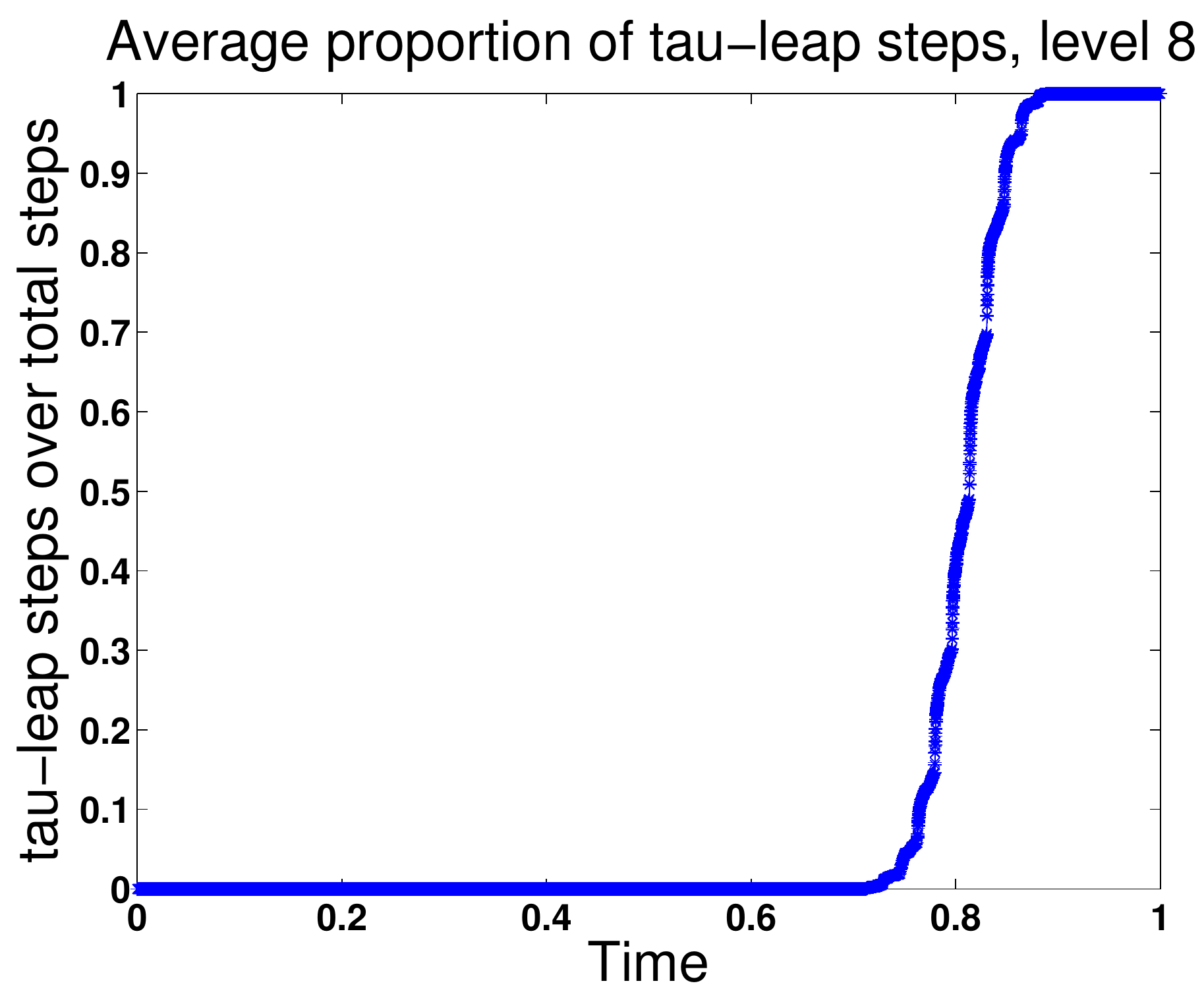}
\end{minipage}
\caption{This figure depicts the `blending' effect produced by our hybrid path-simulation algorithm. Here, we can see the proportion of the tau-leap steps adaptively taken based on expected work optimization. We see how the presence of the tau-leap decreases when we move to the deepest levels. We observe that, for the chosen tolerance, to couple with an exact path at the last level is not optimal.}
\label{fig:blending}
\end{figure}

{We now analyze an ensemble of $10^3$ independent runs of the multilevel estimator, $\mathcal{M}_L$, for $TOL=1.25e{-}2$. 
In this case, $L^*=8$.
In Figures \ref{fig:box1}, \ref{fig:box2} and \ref{fig:box3}, we show boxplots corresponding to that ensemble. In each one, we indicate the coupling pair (on the x-axis) and the value of $\delta_{\ell}$ (below the title of the plot). 
In each figure, the first boxplot starting from the left, corresponds to single-level hybrid simulations at the coarsest level, $\ell{=}0$, with a time mesh of size $\Delta t_0$, and with an exit bound for the one-step exit probability, $\delta_0{=}1e{-}5$. Next, we show the boxplots corresponding to coupled hybrid paths, at levels $\ell{=}0$ and $\ell{=}1$, generated using time meshes of size $\Delta t_0$ and $\Delta t_1$, respectively, and exit probability bounds, $\delta_0$ and $\delta_1$, respectively. This is indicated under the boxplots with the symbols 1C and 1F, which stand for `Coarse' and `Fine' levels in the first coupling, respectively. We proceed in the same fashion until the final level, $L^*$.
At this point, it is crucial to observe that the probability law for the samples in the boxplots indicated by $kF$ and $(k{+}1)C$ should be the same for any $k\in\{0,1,\ldots,L^*{-}1\}$ (in the single-level case, we interpret the symbol 0 as 0F). This is because both samples are generated using the same time mesh of size $\Delta t_k$ and the same one-step exit probability bound, $\delta_k$, see Remark \ref{rem:telescoping}.

Figure \ref{fig:box1} shows the total proportion of Chernoff tau-leap steps over the total number of tau-leap steps.
Here, we understand that a Chernoff tau-leap step is taken when the size of $\tau_{Ch}$ (see Section \ref{hybrid_algo}) is strictly smaller than the distance from the current time to the next mesh point and, therefore, the Chernoff bound is acting as an actual constraint for the size of the tau-leap step. We can see how the Chernoff steps are present in the first levels but not in the final ones, where exact steps are preferred according to our computational work criterion.
We observe a small increase of the proportion of the number of Chernoff steps from levels 1F/2C to levels 3F/4C (strictly speaking, a shift in the median and the third quartile). This is due to consecutive refinements in the values of $\delta$, from $1e{-}5$ to $1e{-}7$, producing smaller and smaller values of $\tau_{Ch}$. This is also because the Chernoff step size is, for some time points, still smaller than the grid size and also because the cost of reaching the time horizon using Chernoff steps is still preferred over the cost of using exact steps. 
The abundance of outliers at all levels up to $\ell{=}3$ indicates that the Chernoff bound is actively controlling the size of the tau-leap steps.  

Figures \ref{fig:box2} and \ref{fig:box3} show the total count of tau-leap and exact steps, respectively. These plots are intended as diagnostic plots with two main objectives: i) checking the telescoping sum property as stated in Remark \ref{rem:telescoping}, and ii) understanding the `blending' phenomena in our simulated hybrid paths, that is, the presence of both methods, tau-leap and exact. 
It could be useful to think in terms of the domain of each method:
given a time mesh, $\Delta t$, and a value of the one-step exit probability bound, $\delta$, we could decompose the interval $[0,T]$ into two domains, $I_{\TL}$ and $I_{\MNRM}$, for the tau-leap and exact methods, respectively. 
The domain, $I_{\MNRM}$, should be monotonically increasing with refinements of the time mesh and $\delta$, since when the size of the time mesh, $\Delta t_\ell$ or $\delta_\ell$, goes to zero, the expected number of tau-leap steps also goes to zero, see \cite{ourSL}, Appendix A. As a consequence, we expect the total count of exact paths to be a monotonically increasing function of the level, $\ell$. On the other hand, the domain $I_{\TL}$ decreases, but, since the size of the time mesh halves form by passing from one level to the next one, we expect to see also an increasing number of tau-leap steps, at least for no very deep levels.
The blending effect of the hybrid decision rules in Algorithm \ref{alg:coupled} are depicted in Figure \ref{fig:blending}, where the proportion of the tau-leap steps over the total number of steps is shown for levels $\ell\in\{0,5,8\}$. In the left panel, we can see that the number of tau-leap steps dominates except close to the origin, where the coarse time-mesh is finer. Remember that in our methodology, our initial mesh can be nonuniform. We then see how the domain, $I_{\MNRM}$, increases until it occupies almost 80\% of the time interval $[0,T]$.
}

\begin{rem}
The savings in computational work when generating Poisson random variables heavily depend on MATLAB's performance capabilities. For example, we do not generate the random variates in batches, as in \cite{Anderson2012}, and that could have an impact on the results. In fact, we should expect better results from our method if we implement our algorithms in more performance-oriented languages or if we sample Poisson random variables in batches.
\end{rem}

{
\begin{rem}(Level 0 time mesh)
In this example, we use an adaptive mesh at level 0. This is because this example is mildly stiff. Using a uniform time mesh at level 0 imposes a small time step size requirement for all time which is not needed. Moreover, this issue is propagated to the finer levels. In all our numerical examples, at level 0, we use the coarsest possible time mesh such that the Forward Euler method is numerically stable.
\end{rem}
}

\section{Conclusions}
\label{sec:conclusions}
In this work, we developed a multilevel Monte Carlo version for the single-level hybrid Chernoff tau-leap algorithm presented in \cite{ourSL}. We showed that the computational complexity of this method is of order $\Ordo{TOL^{-2}}$ and, therefore, that it can be seen as a variance reduction of the SSA method, which has the same complexity.
This represents an important advantage of the hybrid tau-leap with respect to the pure tau-leap in the multilevel context.
In our numerical examples, we obtained substantial gains with respect to both the SSA and the single-level hybrid Chernoff tau-leap.
The present approach, like the one in \cite{ourSL}, also provides an approximation of $\expt{g(X(T))}$ with  prescribed accuracy and  confidence level, with nearly optimal computational work.
For reaching this optimality, we derived novel formulas based on dual-weighted residual estimations for computing the variance of the difference of the observables between two consecutive levels in coupled hybrid paths {and also the bias of the deepest level (see \eqref{eq:EILw} and \eqref{eq:varhatestimated}). 
These formulas are particularly relevant in the present context of Stochastic Reaction Networks due to the fact that alternative standard sample estimators become too costly at deep levels because of the presence of large kurtosis.}

Future extensions may involve better hybridization techniques as well as implicit and higher-order versions of the hybrid MLMC.



\section*{Acknowledgments}
{The authors would like to thank two anonymous reviewers for their constructive comments that helped us to improve our manuscript. We also would like to thank Prof. Mike Giles for very enlightening discussions.}
The authors are members of the KAUST SRI Center for
Uncertainty Quantification in the Computer, Electrical and Mathematical Sciences and Engineering Division at King Abdullah University of Science and Technology (KAUST). This work was supported by KAUST.

\clearpage
\appendix

\begin{algorithm}
\caption{Coupled hybrid path. Inputs: the initial state, $X(0)$, the final time, $T$, the propensity functions, $a{=}\seqof{a_j}{j=1}{J}$, the stoichiometric vectors, $\nu{=}\seqof{\nu_j}{j=1}{J}$, two one-step exit probability bounds; one for the coarse level, $\bar{\delta}$, and another for the fine level, $\bar{\bar{\delta}}$, and  
two time meshes, one coarse $\seqof{t_k}{k=0}{K}$, such that $t_K {=} T$ and a finer one, 
$\seqof{s_l}{l=0}{K'}$, such that $s_0 {=} t_0$, $s_{K'} {=} t_K$, and $\seqof{t_k}{k=0}{K} \subset \seqof{s_l}{l=0}{K'}$.
Outputs: a sequence of states evaluated at the coarse grid, $\seqof{\bar{X}(t_k)}{k=0}{K} \subset \latt$, 
such that $t_K \leq T$, 
a sequence of states evaluated at the fine grid $\seqof{\dbar{X}(s_l)}{l=0}{K'} \subset \latt$,
such that $\bar{X}(t_K) \in \latt$ or $\dbar{X}(s_{K'}) \in \latt$.
If $t_K<T$, both paths exited the $\latt$ lattice before the final time, $T$. 
It also returns the number of times the tau-leap method was successfully applied at the fine level and at the coarse level, and the number of exact steps at the fine level and at the coarse level. For the sake of simplicity, we omit sentences involving the recording of $\bar{X}(t_k)$ and $\dbar{X}(s_l)$ from the current state variables $\bar{X}$ and $\dbar{X}$, respectively, the counting of the number of steps, and the return sentence.}
\label{alg:coupled}
\begin{algorithmic}[1]
\STATE $t \leftarrow 0$, $\bar{X} \leftarrow X(0)$, $\dbar{X} \leftarrow X(0)$ 
\STATE $\bar{t} \leftarrow$ next grid point in the coarse grid larger than $t$ 
\STATE $(\bar{H},\bar{m},\bar{a})  \leftarrow$ Algorithm \ref{alg:timehor} with ($\bar{X}$,$t$,$\bar{t}$,$T$,$\bar{\delta}$,$a$)

\STATE $\dbar{t} \leftarrow$ next grid point in the fine grid larger than $t$ 
\STATE $(\dbar{H},\dbar{m},\dbar{a}) \leftarrow$ Algorithm \ref{alg:timehor} with ($\dbar{X}$,$t$,$\dbar{t}$,$T$,$\dbar{\delta}$,$a$)

\WHILE{$t<T$} 
	\STATE $H \leftarrow \min\{\bar H, \dbar H\}$
		\IF {$\bar m = \TL$ \AND $\dbar m = \TL$}
			\STATE $S \leftarrow$ Algorithm \ref{alg:auxilco2} with ($\bar a,\dbar{a}$)
			\STATE $\Lambda \leftarrow \mathcal{P}(S {\cdot} (H{-}t))$ (generate Poisson random variates)
			\STATE $\dbar{X} \leftarrow \dbar{X} + (\Lambda_1 {+} \Lambda_2)\nu$
			\STATE $\bar{X} \leftarrow \bar{X} + (\Lambda_1 {+} \Lambda_3)\nu$
			\STATE $t \leftarrow H$		
		\ELSE
			\STATE Initialize internal clocks $R, P$ if needed (see Algorithm \ref{alg:mnr})
			\WHILE{$t<H$}
				\IF{$\dbar{m}=\MNRM$}
					\STATE $\dbar{a} \leftarrow a(\dbar{X})$
				\ENDIF
				\IF{$\bar{m}=\MNRM$}
					\STATE $\bar{a} \leftarrow a(\bar{X})$
				\ENDIF
				
				\STATE $S \leftarrow$ Algorithm \ref{alg:auxilco2} with ($\bar{a},\dbar a$)
				\STATE $(t, \bar{X}, \dbar{X}, R, P) \leftarrow$ Algorithm \ref{alg:auxilco} with $(t,H,\bar{X},\dbar{X},R,P,S)$
			\ENDWHILE			
		\ENDIF
\IF{$t<T$}
	\IF{$H = \bar H$}
		\STATE $\bar{t} \leftarrow$ next grid point in the coarse grid larger than $t$ 
		\STATE $(\bar{H},\bar{m},\bar{a}) \leftarrow$ Algorithm \ref{alg:timehor} with ($\bar{X}$,$t$,$\bar{t}$,$T$,$\bar{\delta}$,$a$)
	\ENDIF
	\IF {$H = \dbar H$}
		\STATE $\dbar{t} \leftarrow$ next grid point in the fine grid  larger than $t$ 
		\STATE $(\dbar{H},\dbar{m},\dbar{a}) \leftarrow$ Algorithm \ref{alg:timehor} with ($\dbar{X}$,$t$,$\dbar{t}$,$T$,$\dbar{\delta}$,$a$)
	\ENDIF
\ENDIF
\ENDWHILE

\end{algorithmic}
\end{algorithm}

\begin{algorithm}
\caption{Compute next time horizon. Inputs: the current state, $\tilde{X}$, the current time, $t$, the next grid point, $\tilde{t}$, the final time, $T$, the one step exit probability bound, $\tilde{\delta}$, and the propensity functions, $a{=}\seqof{a_j}{j=1}{J}$.
Outputs: the next horizon $H$, the selected method $m$, current propensity values $\tilde{a}$.}
\label{alg:timehor}
\begin{algorithmic}[1]
\STATE $\tilde{a} \leftarrow a(\tilde{X})$
\STATE $(m,\tilde{\tau}) \leftarrow$ Algorithm \ref{alg:sel} with ($\tilde{X}$,$t$,$\tilde{a}$,$\tilde{\delta}$,$\tilde{t}$)
\IF{$ m= \TL$}
	\STATE $H \leftarrow  \min\{\tilde{t},t{+}\tilde{\tau},T\}$
\ELSE
	\STATE $H \leftarrow  \min\{t{+}\tilde{\tau},T\}$
\ENDIF
\RETURN $(H,m,\tilde a )$
\end{algorithmic}
\end{algorithm}

\begin{algorithm}[h!]
\caption{Auxiliary function used in Algorithm \ref{alg:coupled}. Inputs: the current time, $t$, the current time horizon, $H$, the current system state at coarser level, $\bar{X}$, and finer level, $\dbar{X}$, the internal clocks $R_i$, $P_i$, $i{=}1,2,3$, and the values, $S_i$, $i{=}1,2,3$ (see Section \ref{sec:MNR} for more information on these values). Outputs: updated time, $t$, updated system states, $\bar{X}$, $\dbar{X}$, and updated internal clocks $R_i$, $P_i$, $i{=}1,2,3$.
}
\label{alg:auxilco}
\begin{algorithmic}[1]
\STATE $\Delta t_i \leftarrow (P_i - R_i)/S_i$, for $i{=}1,2,3$
\STATE $\Delta \leftarrow \min_i\{ \Delta t_i \}$
\STATE ${\mu} \leftarrow \text{argmin}_i \{\Delta t_i\}$
					\IF{$t+\Delta > H$}
						\STATE $R \leftarrow R + S {\cdot} (H{-}t)$
						\STATE $t \leftarrow H$
					\ELSE
						\STATE update $\bar{X}$ and $\dbar{X}$
						\STATE $R \leftarrow R + S \Delta$
						\STATE $r \leftarrow$ uniform$(0,1)$
						\STATE $P_\mu \leftarrow P_\mu + \log (1/r)$
						\STATE $t \leftarrow t + \Delta$
					\ENDIF
\RETURN $(t,\bar X,\dbar X,R,P)$
\end{algorithmic}
\end{algorithm}

\begin{algorithm}[h!]
\caption{Auxiliary function used in Algorithm \ref{alg:coupled}. Inputs: the propensity values at the coarse and fine grid, $\bar{a},\dbar{a}$. Output: $S_i$, $i{=}1,2,3$.
}
\label{alg:auxilco2}
\begin{algorithmic}[1]
			\STATE $S_1 \leftarrow \min(\bar{a},\dbar{a})$
			\STATE $S_2 \leftarrow \bar{a}-S_1$
			\STATE $S_3 \leftarrow \dbar{a} -S_1$
			\RETURN $S$
\end{algorithmic}
\end{algorithm}

%


\label{sec:algorithms}
\begin{algorithm}
\caption{\small
Multilevel calibration and error estimation.  Inputs: same as Algorithm \ref{alg:coupled} plus the observable, $g$, and the prescribed tolerance, $TOL{>}0$. Outputs: $\seqof{M_{\ell}}{\ell=0}{L}$, $\seqof{\delta_\ell}{\ell=0}{L}$, $\seqof{ \seqof{t_{n,\ell}}{n=0}{N_\ell}}{\ell=0}{L}$, the estimated computational work of the multilevel estimator, $\hat W_{\ML}$, and the estimated computational work of the SSA method, $\hat W_{\ssa}$. We denote by $g_l \equiv g(\bar X_{l}(T;\bar{\omega}))$, and $g_{l{+1}}{-}g_l \equiv g(\bar X_{{l{+}1}}(T;\bar{\omega})){-g(\bar X_{{l}}(T;\bar{\omega}))}$. Here, $C^*$ is the unitary cost of a pure SSA step, and $c$ is the factor of refinement of $\delta$ (in our experiments $c{=}10$). 
See also Remark \ref{rem:algCal} regarding the estimators of $\var{g(X(T))}$ and $\expt{g(X(T))}$, and Remark \ref{rem:unbiased}. 
}
\label{alg:Cal}
\begin{algorithmic}[1]
\begin{small}
\STATE $l \leftarrow 0$, $\delta_l \leftarrow 0.01$, $\hat W_{ML}^{(a)} \leftarrow \infty$ 
\STATE Set initial meshes $\seqof{t_k}{k=0}{K}$ and $\seqof{s_l}{l=0}{K'}$
\STATE fin-delta $\leftarrow$ \FALSE
\WHILE {\NOT fin-delta}
	\STATE $(\hat \psi_0, \svar{g_l}{\cdot}, \avg{\{g_{l},\WE,\NSSAP, \NTL\}}{\cdot}) \leftarrow$ Algorithm \ref{alg:pathdualscostSL}
	\IF{$\hat{\mathcal{V}}_l (1{-}\delta_l \avg{\NTL}{\cdot})^2 \geq 2 \svar{g_l}{\cdot} \delta_l  \avg{\NTL}{\cdot}$ \AND $\delta_l \avg{\NTL}{\cdot} < 0.1$}
		\STATE fin-delta $\leftarrow$ \TRUE
		\STATE  Refine $\delta_l$ by a factor of $c$
	\ENDIF
\ENDWHILE
\STATE $\delta_{l+1} \leftarrow \delta_l$
\STATE fin $\leftarrow$ \FALSE
\WHILE {\NOT fin}
	\STATE fin-delta $\leftarrow$ \FALSE
	\WHILE {\NOT fin-delta}
		\STATE $(\hat{\psi}_{l+1}, \hat{\mathcal{V}}_{l+1}, \avg{\{g_{l{+}1},\NSSAP,\WE,N_{TL,l+1}\}}{\cdot},\svar{g_{l+1}}{\cdot}) \leftarrow$ Algorithm \ref{alg:pathdualscost}
		
		\IF{$\hat{\mathcal{V}}_{l+1} (1{-}\delta_{l+1} \avg{N_{TL,l+1}}{\cdot})^2 \geq 2 \svar{g_{l+1}}{\cdot} \delta_{l+1}  \avg{N_{TL,l+1}}{\cdot}$ \\ \AND $\delta_{l+1} \avg{N_{TL,l+1}}{\cdot} < 0.1$}
			\STATE fin-delta $\leftarrow$ \TRUE
			\STATE $\delta_l \leftarrow \delta_{l+1}$
		\ELSE
			\STATE Refine $\delta_{l+1}$ by a factor of $c$
		\ENDIF
	\ENDWHILE
	\STATE $M_{\ssa} \leftarrow {C_A^2 \svar{g_{l{+}1}}{\cdot}}/{TOL^2}$
	\STATE $\hat W_{\ssa} \leftarrow C^* M_{\ssa} \avg{\NSSAP}{\cdot}$
	\IF {$\hat \WE < TOL{-}TOL^{2}$}
		\STATE $\seqof{M_\ell}{\ell=0}{l{+}1} \leftarrow$ Algorithm \ref{alg:greedyoptKKT} with $(\seqof{\hat{\psi}_\ell}{\ell=0}{l{+}1},\seqof{\hat{\mathcal{V}}_\ell}{\ell=0}{l+1},TOL,\hat \WE)$
		\STATE $\hat W_{\ML} \leftarrow \sum_{\ell=0}^{l+1}\hat \psi_\ell M_\ell$
	\ELSE
		\STATE $\hat W_{\ML} \leftarrow \infty$
	\ENDIF
	\IF {($\hat W_{\ML}^{(a)} > \hat W_{\ML}$ \OR $\hat \WE > TOL{-}TOL^{2})$ \AND $\avg{N_{TL,l+1}}{\cdot} > 0$}
		\STATE $l \leftarrow l+1$
		\STATE $\hat W_{\ML}^{(a)} \leftarrow \hat W_{\ML}$
		\STATE Refine meshes $\seqof{t_k}{k=0}{K}$ and $\seqof{s_l}{l=0}{K'}$
	\ELSE
		\STATE fin $\leftarrow \delta_{l+1} \avg{N_{\TL,l+1}}{\cdot}\avg{g_{l{+}1}}{\cdot} \leq TOL^2$
		\IF {\NOT fin}
			\STATE $\delta_{l+1} \leftarrow c^{\lfloor (log_c(TOL^2/(\avg{g_{l{+}1}}{\cdot} \cdot \avg{N_{\TL,l+1}}{\cdot}) \rfloor}$
			\WHILE {\NOT fin}
				\STATE $(\avg{\{g_{l{+}1},N_{TL,l+1}\}}{\cdot}) \leftarrow$ Algorithm \ref{alg:pathdualscost}
				\STATE fin $\leftarrow \delta_{l+1} \avg{N_{\TL,l+1}}{\cdot}\avg{g_{l{+}1}}{\cdot} \leq TOL^2$
				\IF {\NOT fin}
					\STATE Refine $\delta_{l+1}$ by a factor of $c$
				\ENDIF
			\ENDWHILE
		\ENDIF
		\IF {$\avg{N_{\TL,l+1}}{\cdot}=0$}
			\STATE $l \leftarrow l+1$
			\STATE $\seqof{M_\ell}{\ell=0}{l} \leftarrow$ Algorithm \ref{alg:greedyoptKKT} with $(\seqof{\hat{\psi}_\ell}{\ell=0}{l},\seqof{\hat{\mathcal{V}}_\ell}{\ell=0}{l},TOL,0)$
		\STATE $\hat W_{\ML} \leftarrow \sum_{\ell=0}^{l}\hat \psi_\ell M_\ell$
		\ENDIF	
	\ENDIF
\ENDWHILE
\end{small}
\end{algorithmic}
\end{algorithm}

\begin{algorithm}
\caption{Auxiliary function for Algorithm \ref{alg:Cal}. Inputs: same as Algorithm \ref{alg:coupled}. 
Outputs: the estimated runtime of the coupled path, $\hat \psi$, 
an estimate of $\var{g(\bar X(T)){-}g(\dbar X(T))}$, $\hat{\mathcal{V}}$, 
an estimate of $\expt{g(X(T))}$, $\avg{g(\dbar X(T)}{\cdot}$, an estimate of the expected number of steps needed by the SSA method, $\avg{\NSSAP}{\cdot})$, an estimate of $\expt{\WE}$, $\avg{\WE}{\cdot}$, an estimate of the expected number of tau-leap steps taken at the fine level, $\avg{\NTL }{\cdot}$, and an estimate of $\var{g(X(T))}$, $\svar{g(\dbar X(T)}{\cdot}$. Here, $\dbar X(t)$ refers to the approximated process using a finer grid than the approximated process, $\bar X(t)$. Moreover, ($\bar X(t)$,$\dbar X(t)$) are two coupled paths. Here, $\indicator{TL} (k)=1$ if and only if the decision at time $t_k$ was tau-leap.
Set appropriate values for $M_0$ and $CV_0$. For the sake of simplicity, we omit the arguments of the algorithms  when there is no risk of confusion. See also Remark \ref{rem:algCal} regarding the estimators of $\var{g(X(T))}$ and $\expt{g(X(T))}$.
}
\label{alg:pathdualscost}
\begin{algorithmic}[1]
\STATE $M \leftarrow M_0$, $cv \leftarrow \infty$, $M_f \leftarrow 0$ 
\WHILE{$cv > CV_0$} 
	\FOR{$m \leftarrow 1$ to $M$ }
		\STATE Generate two coupled paths: $(\bar{X}(s_l;\dbar {\omega}_m))_{l=0}^{K'},(\dbar{X}(s_l;\dbar {\omega}_m))_{l=0}^{K'},   \leftarrow$ Algorithm \ref{alg:coupled}
		\IF {the path does not exit $\latt$}
			\STATE $M_f \leftarrow M_f + 1$
			\STATE $(S_e(\dbar {\omega}_m),S_v(\dbar {\omega}_m)) \leftarrow$ Algorithm \ref{alg:varggl} with $(\bar{X}(t_k;\dbar {\omega}_m))_{k=0}^{K}$
			\STATE $\WE(\dbar {\omega}_m) \leftarrow$ Algorithm \ref{alg:weakerror} with $\seqof{\dbar{X}(s_l;\dbar{\omega}_m)}{l=0}{K'}$
			\STATE Estimate $\NSSAP (\dbar {\omega}_m)$, using  $\int_0^T a_0 (\dbar{X}(s))ds$ (see \cite{ourSL})	
			\STATE $C_{Poi}(\bar{\bar{\omega}}_m) \leftarrow 
\sum_{j=1}^J \sum_{l=0}^{K'} C_P(a_j(\bar X(s_l))(s_{l+1}{-}s_{l})\indicator{TL}(l)$
			\STATE \quad \quad \quad \quad \quad $+\, \sum_{j=1}^J  \sum_{l=0}^{K'} C_P(a_j(\bar{ \bar X} (s_l))(s_{l+1}{-}s_l)\indicator{TL}(l)$
			\STATE Compute $\NMNRMKoneC$, $\NMNRMKtwoC$, $\NTLC$, and $N_{\TL}$
		\ENDIF
	\ENDFOR
	\STATE $\hat{\mathcal{V}} \leftarrow \svar{S_e}{M_f}+\avg{S_v}{M_f}$
	\STATE Compute the coefficient of variation $cv_\mathcal{V}$ and $cv_{\WE}$ of $\hat{\mathcal{V}}$ and $\avg{\WE}{\cdot}$, respectively.\!\!\!\!\!\!
	
	\STATE $cv \leftarrow \max\{cv_\mathcal{V},cv_{\WE}\}$ 
	\STATE $\hat{\psi} \leftarrow C_1\avg{\NMNRMKoneC}{M_f} {+} C_2\avg{\NMNRMKtwoC}{M_f} {+} C_3 \avg{\NTLC}{M_f} {+} \avg{C_{Poi}}{M_f}$
	\STATE $M \leftarrow 2 M$
\ENDWHILE
\RETURN $\LP \hat{\psi}, \hat{\mathcal{V}}, \avg{\{g(\dbar{X}(T)),\NSSAP,\WE,N_{\TL}\}}{M_f},\svar{g(\dbar{X}(T))}{M_f}\RP$
\end{algorithmic}
\end{algorithm}

\begin{algorithm}[h!]
\caption{Compute the discretization error of a given approximated path. 
Inputs: $(\bar X(t_k))_{k=0}^K$. 
Here, $\indicator{TL} (k)=1$ if and only if the decision at time $t_k$ was tau-leap, and $Id$ is the $d\times d$ identity matrix
Output: $\WE$. Notes: $x_k \equiv \bar X(t_k)$.}
\label{alg:weakerror}
\begin{algorithmic}[1]
\STATE $\EE_I \leftarrow 0$
\STATE Compute $\varphi_K \leftarrow \nabla g(x_k)$ 
\FOR {$k \leftarrow K{-}1$ to $1$}		 
	\STATE  $\Delta t_k \leftarrow t_{k+1} - t_k$
	\STATE Compute $\JAC_a = [\partial_i a_j(x_k)]_{j,i}$
	\STATE  $\varphi_k \leftarrow \LP Id+{\Delta t_k}\, \JAC_a^T \,\nu^T\RP \varphi_{k+1}$
	\STATE $\Delta a_k \leftarrow a(x_{k+1}) - a(x_k)$
	\STATE $\EE_I \leftarrow \EE_I  + \frac{\Delta t_k}{2} ( \Delta a_k \,\indicator{TL} (k)\,\nu^T)\,\varphi_k$	
\ENDFOR
\RETURN  $\EE_I$
\end{algorithmic}
\end{algorithm}


\begin{algorithm}
\caption{
Compute $S_e\equiv S_e(\bar \omega)$ and $S_v\equiv S_v(\bar \omega)$ defined in \eqref{eq:varhat}. Inputs: $(\bar X(t_k))_{k=0}^K$ and a positive constant $c$.
Outputs: $S_e$ and $S_v$. Notes: if $a$ is a vector, then, $\text{diag}(a)$ is a diagonal matrix with main diagonal $a$. Here, $\indicator{TL} (k)=1$ if and only if the decision at time $t_k$ was tau-leap, $Id$ is the $d\times d$ identity matrix, $x_k \equiv \bar X(t_k)$, and $\Phi(x)$ is the cumulative distribution function of a Gaussian  random variable.
}
\label{alg:varggl}
\begin{algorithmic}[1]
\STATE $S_e \leftarrow 0$
\STATE $S_v \leftarrow 0$
\STATE Compute $\varphi_K \leftarrow \nabla g(x_k)$ 
\FOR {$k \leftarrow K{-}1$ to $1$}		 
	\STATE  $\Delta t_k \leftarrow t_{k+1} - t_k$
	\STATE Compute $\JAC_a = [\partial_i a_j(x_k)]_{j,i}$
	\STATE $\varphi_k \leftarrow \LP Id+{\Delta t_k}\, \JAC_a^T \,\nu^T\RP \varphi_{k+1}$
	\STATE $\nu_{\varphi} \leftarrow \nu^T \varphi_k$
	\STATE $\nu_{a}\leftarrow (\JAC_a\,\nu)^T$
	\STATE $\mu_j \leftarrow \frac{\Delta t_k}{2} \sum_i (\nabla a_j(x_k){\cdot}\nu_i)\,a_i(x_k)$ 
	\STATE $\bar{\mu}_j \leftarrow \frac{\Delta t_k}{2} \sum_i |(\nabla a_j(x_k){\cdot}\nu_i)|\,a_i(x_k)$ 
	\STATE $\sigma^2_j \leftarrow \frac{\Delta t_k}{2} \sum_i (\nabla a_j(x_k){\cdot}\nu_i)^2\,a_i(x_k)$ 
	\STATE $S_e {\leftarrow} S_e {+} \indicator{\TL} (k)  \frac{\Delta t_k}{2} \mu\,\nu_{\varphi}$
	\STATE $\text{aux}_1 \leftarrow \frac{(\Delta t_k)^3}{8} (\nu_a \nu_{\varphi})^T \text{diag}(a)(\nu_a \nu_{\varphi})$
	\STATE $\text{aux}_2 \leftarrow \frac{\Delta t_k}{2}
	 \sum_j (\varphi_k{\cdot}\nu_j)^2 \indicator{\left\{\frac{\Delta t_k}{2} a_j(x_k)>c \right\}}	 
	 \left(\mu_j(1-2\Phi(-\frac{\mu_j}{\sigma_j}))+\sqrt{\frac{2}{\pi}}\,\sigma_j\,\exp(-\frac{1}{2} (\frac{\mu_j}{\sigma_j})^2)\right)$ 	
	\STATE $\text{aux}_3 \leftarrow \frac{\Delta t_k}{2}
	 \sum_j (\varphi_k{\cdot}\nu_j)^2 \indicator{\left\{\frac{\Delta t_k}{2} a_j(x_k)<c\right\}}  
	 \min\left\{\bar{\mu}_{j},\sqrt{\mu_j^2+\sigma^2_j}\right\}$
 	\STATE $S_v {\leftarrow} S_v {+} \indicator{\TL} (k)(\text{aux}_1+\text{aux}_2+\text{aux}_3)$ 
\ENDFOR
\RETURN  $(S_e,S_v)$
\end{algorithmic}
\end{algorithm}

\begin{algorithm}
\caption{Solve the optimization problem \eqref{eq:greedyKKT} using a greedy scheme. Inputs:  the estimations of the coupled path cost for all the levels, $\seqof{\hat{\psi}_\ell}{{\ell}=0}{L}$, the estimation of the variance of the quantity of interest at level 0, $\hat{\mathcal{V}}_{0}$, the estimations of the differences  of the quantity of interest for all the coupled levels, $\seqof{\hat{\mathcal{V}}_{\ell}}{{\ell}=1}{L}$, the prescribed tolerance, $TOL$, and the weak error estimation for level $L$, $\WE$. Output: the number of realizations needed for each level, $\seqof{M}{{\ell}=0}{L}$. }
\label{alg:greedyoptKKT}
\begin{algorithmic}[1]
\item[]  Define $q_k :=  \dfrac{\sum_{{\ell}=0}^{L-k} \sqrt{\hat{\psi}_\ell \hat{\mathcal{V}}_{\ell}}}{RHS{-}\sum_{{\ell}=L-k+1}^{L}\hat{\mathcal{V}}_{\ell}}$\STATE $RHS \leftarrow ((TOL{-}TOL^{2}-\WE)/C_A)^2$
\STATE fin $\leftarrow \FALSE$
\STATE $k \leftarrow 0$
\WHILE {\NOT fin \AND $k \leq L$}
\IF {$\hat{\psi}_{L{-}k}-q^2_k \,\,\hat{\mathcal{V}}_{L{-}k} < 0$}
\STATE fin $\leftarrow \TRUE$
\STATE $ \seqof{M_{\ell}}{\ell=0}{L{-}k} \leftarrow q_k \sqrt{\hat{\mathcal{V}}_{\ell}/\hat{\psi}_{\ell}}$
\ELSE
\STATE $M_{L{-}k} \leftarrow 1$
\STATE $k \leftarrow k+1$
\ENDIF
\ENDWHILE
\RETURN $\seqof{M_{\ell}}{\ell=0}{L}$
\end{algorithmic}
\end{algorithm}

\begin{algorithm}[h!]
\caption{Auxiliary function for Algorithm \ref{alg:Cal}. Inputs: same as Algorithm \ref{alg:coupled}. 
Outputs: the estimated runtime of the hybrid path at level 0, $\hat \psi_0$, 
an estimate of $\var{g(\bar X(T))}$, $\svar{g(\bar X(T)}{\cdot}$, an estimate of $\expt{g(X(T))}$, $\avg{g(\bar X(T)}{\cdot}$, an estimate of $\expt{\WE}$, $\avg{\WE}{\cdot}$, an estimate of the expected number of steps needed by the SSA method, $\avg{\NSSAP}{\cdot})$ and $\avg{\NTL }{\cdot}$. Here, $\indicator{TL} (k)=1$ if and only if the decision at time $t_k$ was tau-leap.
Notes: the values $C_1$, $C_2$ and $C_3$ are described in Section \ref{EEC}. Set appropriate values for $M_0$ and $CV_0$. For the sake of simplicity, we omit the arguments of the algorithms  when there is no risk of confusion.}
\label{alg:pathdualscostSL}
\begin{algorithmic}[1]
\STATE $M \leftarrow M_0$, $cv \leftarrow \infty$, $M_f \leftarrow 0$ 
\WHILE{$cv > CV_0$} 
	\FOR{$m \leftarrow 1$ to $M$ }
	\STATE $((\bar{X}(t_k))_{k=0}^K,\NTL,\NSSAKone,\NSSAKtwo) \leftarrow$ generate one hybrid path (see \cite{ourSL}) 
		\IF {the path does not exit $\latt$}	
			\STATE $M_f \leftarrow M_f+1$
			\STATE Compute $g(\bar X(T;\bar {\omega}_m))$
			\STATE $\WE \leftarrow$ Algorithm \ref{alg:weakerror} 
with $(\bar{X}(t_k))_{k=0}^K$
			\STATE $(S_e(\dbar {\omega}_m),S_v(\dbar {\omega}_m)) \leftarrow$ Algorithm \ref{alg:varggl} with $(\bar{X}(t_k))_{k=0}^{K}$
		\STATE Estimate $\NSSAP (\dbar {\omega}_m)$, using  $\int_0^T a_0 (\dbar{X}(s))ds$ (see \cite{ourSL})	
			\STATE $C_{Poi}(\bar \omega_m) \leftarrow \sum_{j=1}^J \sum_{k=0}^K C_P(a_j(\bar X(t_k))(t_{k+1}{-}t_k))\indicator{TL}(k)$
		\ENDIF
	\ENDFOR
	\STATE $\hat{\mathcal{V}} \leftarrow \svar{S_e}{M_f}+\avg{S_v}{M_f}$
	\STATE Estimate the coefficients of variation $cv_\mathcal{V}$, $cv_g$ and $cv_{\WE}$ of the estimators of $\var{g(\bar X(T)){-}g(\dbar X(T))}$, $\var{g(\bar X(T))}$ and $\expt{\WE}$, respectively. 
	\STATE $cv\leftarrow \max\{cv_\mathcal{V},cv_g,cv_{\WE}\}$
	\STATE $\hat \psi_0 {\leftarrow} C_1\avg{\NSSAKone}{M_f} {+} C_2\avg{\NSSAKtwo}{M_f} {+} C_3 \avg{\NTL}{M_f} {+} \avg{C_{Poi}}{M_f}$\!\!\!
	
	\STATE $M \leftarrow 2 M$
\ENDWHILE
\RETURN $(\hat \psi_0, \svar{g(\bar X(T))}{M_f}, \avg{\{g(\bar X(T)),\WE,\NSSAP, \NTL\}}{M_f})$
\end{algorithmic}
\end{algorithm}



\newpage


\end{document}